
\ifx\shlhetal\undefinedcontrolsequence\let\shlhetal\relax\fi

\input amstex
\expandafter\ifx\csname mathdefs.tex\endcsname\relax
  \expandafter\gdef\csname mathdefs.tex\endcsname{}
\else \message{Hey!  Apparently you were trying to
  \string\input{mathdefs.tex} twice.   This does not make sense.} 
\errmessage{Please edit your file (probably \jobname.tex) and remove
any duplicate ``\string\input'' lines}\endinput\fi




\catcode`\X=12\catcode`\@=11

\def\n@wcount{\alloc@0\count\countdef\insc@unt}
\def\n@wwrite{\alloc@7\write\chardef\sixt@@n}
\def\n@wread{\alloc@6\read\chardef\sixt@@n}
\def\r@s@t{\relax}\def\v@idline{\par}\def\@mputate#1/{#1}
\def\l@c@l#1X{\firstpart.#1}\def\gl@b@l#1X{#1}\def\t@d@l#1X{{}}

\def\crossrefs#1{\ifx\all#1\let\tr@ce=\all\else\def\tr@ce{#1,}\fi
   \n@wwrite\cit@tionsout\openout\cit@tionsout=\jobname.cit 
   \write\cit@tionsout{\tr@ce}\expandafter\setfl@gs\tr@ce,}
\def\setfl@gs#1,{\def\@{#1}\ifx\@\empty\let\next=\relax
   \else\let\next=\setfl@gs\expandafter\xdef
   \csname#1tr@cetrue\endcsname{}\fi\next}
\def\m@ketag#1#2{\expandafter\n@wcount\csname#2tagno\endcsname
     \csname#2tagno\endcsname=0\let\tail=\all\xdef\all{\tail#2,}
   \ifx#1\l@c@l\let\tail=\r@s@t\xdef\r@s@t{\csname#2tagno\endcsname=0\tail}\fi
   \expandafter\gdef\csname#2cite\endcsname##1{\expandafter
     \ifx\csname#2tag##1\endcsname\relax?\else\csname#2tag##1\endcsname\fi
     \expandafter\ifx\csname#2tr@cetrue\endcsname\relax\else
     \write\cit@tionsout{#2tag ##1 cited on page \folio.}\fi}
   \expandafter\gdef\csname#2page\endcsname##1{\expandafter
     \ifx\csname#2page##1\endcsname\relax?\else\csname#2page##1\endcsname\fi
     \expandafter\ifx\csname#2tr@cetrue\endcsname\relax\else
     \write\cit@tionsout{#2tag ##1 cited on page \folio.}\fi}
   \expandafter\gdef\csname#2tag\endcsname##1{\expandafter
      \ifx\csname#2check##1\endcsname\relax
      \expandafter\xdef\csname#2check##1\endcsname{}%
      \else\immediate\write16{Warning: #2tag ##1 used more than once.}\fi
      \multit@g{#1}{#2}##1/X%
      \write\t@gsout{#2tag ##1 assigned number \csname#2tag##1\endcsname\space
      on page \number\count0.}%
   \csname#2tag##1\endcsname}}

\def\multit@g#1#2#3/#4X{\def\t@mp{#4}\ifx\t@mp\empty%
      \global\advance\csname#2tagno\endcsname by 1 
      \expandafter\xdef\csname#2tag#3\endcsname
      {#1\number\csname#2tagno\endcsnameX}%
   \else\expandafter\ifx\csname#2last#3\endcsname\relax
      \expandafter\n@wcount\csname#2last#3\endcsname
      \global\advance\csname#2tagno\endcsname by 1 
      \expandafter\xdef\csname#2tag#3\endcsname
      {#1\number\csname#2tagno\endcsnameX}
      \write\t@gsout{#2tag #3 assigned number \csname#2tag#3\endcsname\space
      on page \number\count0.}\fi
   \global\advance\csname#2last#3\endcsname by 1
   \def\t@mp{\expandafter\xdef\csname#2tag#3/}%
   \expandafter\t@mp\@mputate#4\endcsname
   {\csname#2tag#3\endcsname\lastpart{\csname#2last#3\endcsname}}\fi}
\def\t@gs#1{\def\all{}\m@ketag#1e\m@ketag#1s\m@ketag\t@d@l p
\let\realscite\scite
\let\realstag\stag
   \m@ketag\gl@b@l r \n@wread\t@gsin
   \openin\t@gsin=\jobname.tgs \re@der \closein\t@gsin
   \n@wwrite\t@gsout\openout\t@gsout=\jobname.tgs }
\outer\def\localtags{\t@gs\l@c@l}
\outer\def\globaltags{\t@gs\gl@b@l}
\outer\def\newlocaltag#1{\m@ketag\l@c@l{#1}}
\outer\def\newglobaltag#1{\m@ketag\gl@b@l{#1}}

\newif\ifpr@ 
\def\m@kecs #1tag #2 assigned number #3 on page #4.%
   {\expandafter\gdef\csname#1tag#2\endcsname{#3}
   \expandafter\gdef\csname#1page#2\endcsname{#4}
   \ifpr@\expandafter\xdef\csname#1check#2\endcsname{}\fi}
\def\re@der{\ifeof\t@gsin\let\next=\relax\else
   \read\t@gsin to\t@gline\ifx\t@gline\v@idline\else
   \expandafter\m@kecs \t@gline\fi\let \next=\re@der\fi\next}
\def\pretags#1{\pr@true\pret@gs#1,,}
\def\pret@gs#1,{\def\@{#1}\ifx\@\empty\let\n@xtfile=\relax
   \else\let\n@xtfile=\pret@gs \openin\t@gsin=#1.tgs \message{#1} \re@der 
   \closein\t@gsin\fi \n@xtfile}

\newcount\sectno\sectno=0\newcount\subsectno\subsectno=0
\newif\ifultr@local \def\ultralocal{\ultr@localtrue}
\def\firstpart{\number\sectno}
\def\lastpart#1{\ifcase#1 \or a\or b\or c\or d\or e\or f\or g\or h\or 
   i\or k\or l\or m\or n\or o\or p\or q\or r\or s\or t\or u\or v\or w\or 
   x\or y\or z \fi}

\def\resetall{\global\advance\sectno by 1\subsectno=0
   \gdef\firstpart{\number\sectno}\r@s@t}
\def\resetsub{\global\advance\subsectno by 1
   \gdef\firstpart{\number\sectno.\number\subsectno}\r@s@t}
\def\newsection#1\par{\resetall\vskip0pt plus.3\vsize\penalty-250
   \vskip0pt plus-.3\vsize\bigskip\bigskip
   \message{#1}\leftline{\bf#1}\nobreak\bigskip}
\def\subsection#1\par{\ifultr@local\resetsub\fi
   \vskip0pt plus.2\vsize\penalty-250\vskip0pt plus-.2\vsize
   \bigskip\smallskip\message{#1}\leftline{\bf#1}\nobreak\medskip}


\newdimen\marginshift

\newdimen\margindelta
\newdimen\marginmax
\newdimen\marginmin

\def\margininit{       
\marginmax=3 true cm                  
				      
\margindelta=0.1 true cm              
\marginmin=0.1true cm                 
\marginshift=\marginmin
}    

\def\t@gsjj#1,{\def\@{#1}\ifx\@\empty\let\next=\relax\else\let\next=\t@gsjj
   \def\@@{p}\ifx\@\@@\else
   \expandafter\gdef\csname#1cite\endcsname##1{\citejj{##1}}
   \expandafter\gdef\csname#1page\endcsname##1{?}
   \expandafter\gdef\csname#1tag\endcsname##1{\tagjj{##1}}\fi\fi\next}
\newif\ifshowstuffinmargin
\showstuffinmarginfalse
\def\jjtags{\ifx\shlhetal\relax 
  \else
\ifx\shlhetal\undefinedcontrolseq
\else
\showstuffinmargintrue
\ifx\all\relax\else\expandafter\t@gsjj\all,\fi\fi \fi
}

\def\tagjj#1{\realstag{#1}\oldmginpar{\zeigen{#1}}}
\def\citejj#1{\rechnen{#1}\mginpar{\zeigen{#1}}}     

\def\rechnen#1{\expandafter\ifx\csname stag#1\endcsname\relax ??\else
                           \csname stag#1\endcsname\fi}

\newdimen\theight

\def\marginfont{\sevenrm}

\def\trymarginbox#1{\setbox0=\hbox{\marginfont\hskip\marginshift #1}%
		\global\marginshift\wd0 
		\global\advance\marginshift\margindelta}

\def \oldmginpar#1{%
\ifvmode\setbox0\hbox to \hsize{\hfill\rlap{\marginfont\quad#1}}%
\ht0 0cm
\dp0 0cm
\box0\vskip-\baselineskip
\else 
             \vadjust{\trymarginbox{#1}%
		\ifdim\marginshift>\marginmax \global\marginshift\marginmin
			\trymarginbox{#1}%
                \fi
             \theight=\ht0
             \advance\theight by \dp0    \advance\theight by \lineskip
             \kern -\theight \vbox to \theight{\rightline{\rlap{\box0}}%
\vss}}\fi}

\newdimen\upordown
\global\upordown=8pt
\font\tinyfont=cmtt8 
\def\mginpar#1{\smash{\hbox to 0cm{\kern-10pt\raise7pt\hbox{\tinyfont #1}\hss}}}
\def\mginpar#1{{\hbox to 0cm{\kern-10pt\raise\upordown\hbox{\tinyfont #1}\hss}}\global\upordown-\upordown}


\def\t@gsoff#1,{\def\@{#1}\ifx\@\empty\let\next=\relax\else\let\next=\t@gsoff
   \def\@@{p}\ifx\@\@@\else
   \expandafter\gdef\csname#1cite\endcsname##1{\zeigen{##1}}
   \expandafter\gdef\csname#1page\endcsname##1{?}
   \expandafter\gdef\csname#1tag\endcsname##1{\zeigen{##1}}\fi\fi\next}
\def\verbatimtags{\showstuffinmarginfalse
\ifx\all\relax\else\expandafter\t@gsoff\all,\fi}
\def\zeigen#1{\hbox{$\scriptstyle\langle$}#1\hbox{$\scriptstyle\rangle$}}


\def\margintag#1{\ifshowstuffinmargin\oldmginpar{\zeigen{#1}}\fi}

\def\(#1){\edef\dot@g{\ifmmode\ifinner(\hbox{\noexpand\etag{#1}})
   \else\noexpand\eqno(\hbox{\noexpand\etag{#1}})\fi
   \else(\noexpand\ecite{#1})\fi}\dot@g}

\newif\ifbr@ck
\def\eat#1{}
\def\[#1]{\br@cktrue[\br@cket#1'X]}
\def\br@cket#1'#2X{\def\temp{#2}\ifx\temp\empty\let\next\eat
   \else\let\next\br@cket\fi
   \ifbr@ck\br@ckfalse\br@ck@t#1,X\else\br@cktrue#1\fi\next#2X}
\def\br@ck@t#1,#2X{\def\temp{#2}\ifx\temp\empty\let\neext\eat
   \else\let\neext\br@ck@t\def\temp{,}\fi
   \def\teemp{#1}\ifx\teemp\empty\else\rcite{#1}\fi\temp\neext#2X}
\def\resetbr@cket{\gdef\[##1]{[\rtag{##1}]}}
\def\references{\resetbr@cket\newsection References\par}

\newtoks\symb@ls\newtoks\s@mb@ls\newtoks\p@gelist\n@wcount\ftn@mber
    \ftn@mber=1\newif\ifftn@mbers\ftn@mbersfalse\newif\ifbyp@ge\byp@gefalse
\def\defm@rk{\ifftn@mbers\n@mberm@rk\else\symb@lm@rk\fi}
\def\n@mberm@rk{\xdef\m@rk{{\the\ftn@mber}}%
    \global\advance\ftn@mber by 1 }
\def\rot@te#1{\let\temp=#1\global#1=\expandafter\r@t@te\the\temp,X}
\def\r@t@te#1,#2X{{#2#1}\xdef\m@rk{{#1}}}
\def\b@@st#1{{$^{#1}$}}\def\str@p#1{#1}
\def\symb@lm@rk{\ifbyp@ge\rot@te\p@gelist\ifnum\expandafter\str@p\m@rk=1 
    \s@mb@ls=\symb@ls\fi\write\f@nsout{\number\count0}\fi \rot@te\s@mb@ls}
\def\byp@ge{\byp@getrue\n@wwrite\f@nsin\openin\f@nsin=\jobname.fns 
    \n@wcount\currentp@ge\currentp@ge=0\p@gelist={0}
    \re@dfns\closein\f@nsin\rot@te\p@gelist
    \n@wread\f@nsout\openout\f@nsout=\jobname.fns }
\def\m@kelist#1X#2{{#1,#2}}
\def\re@dfns{\ifeof\f@nsin\let\next=\relax\else\read\f@nsin to \f@nline
    \ifx\f@nline\v@idline\else\let\t@mplist=\p@gelist
    \ifnum\currentp@ge=\f@nline
    \global\p@gelist=\expandafter\m@kelist\the\t@mplistX0
    \else\currentp@ge=\f@nline
    \global\p@gelist=\expandafter\m@kelist\the\t@mplistX1\fi\fi
    \let\next=\re@dfns\fi\next}
\def\symbols#1{\symb@ls={#1}\s@mb@ls=\symb@ls} 
\def\bigsymbol{\textstyle}
\symbols{\bigsymbol\ast,\dagger,\ddagger,\sharp,\flat,\natural,\star}
\def\ftnumbers{\ftn@mberstrue} \def\ftsymbols{\ftn@mbersfalse}
\def\paginal{\byp@ge} \def\resetftnumbers{\ftn@mber=1}
\def\ftnote#1{\defm@rk\expandafter\expandafter\expandafter\footnote
    \expandafter\b@@st\m@rk{#1}}

\long\def\jump#1\endjump{}
\def\ssum{\mathop{\lower .1em\hbox{$\textstyle\Sigma$}}\nolimits}

\def\qed{\nobreak\kern 1em \vrule height .5em width .5em depth 0em}
\def\newneq{\hbox{\rlap{\hbox to 1\wd9{\hss$=$\hss}}\raise .1em 
   \hbox to 1\wd9{\hss$\scriptscriptstyle/$\hss}}}
\def\subsetne{\setbox9 = \hbox{$\subset$}\mathrel{\hbox{\rlap
   {\lower .4em \newneq}\raise .13em \hbox{$\subset$}}}}
\def\supsetne{\setbox9 = \hbox{$\subset$}\mathrel{\hbox{\rlap
   {\lower .4em \newneq}\raise .13em \hbox{$\supset$}}}}

\def\vbar{\mathchoice{\vrule height6.3ptdepth-.5ptwidth.8pt\kern-.8pt}
   {\vrule height6.3ptdepth-.5ptwidth.8pt\kern-.8pt}
   {\vrule height4.1ptdepth-.35ptwidth.6pt\kern-.6pt}
   {\vrule height3.1ptdepth-.25ptwidth.5pt\kern-.5pt}}
\def\f@dge{\mathchoice{}{}{\mkern.5mu}{\mkern.8mu}}
\def\b@c#1#2{{\rm \mkern#2mu\vbar\mkern-#2mu#1}}
\def\b@b#1{{\rm I\mkern-3.5mu #1}}
\def\b@a#1#2{{\rm #1\mkern-#2mu\f@dge #1}}
\def\bb#1{{\count4=`#1 \advance\count4by-64 \ifcase\count4\or\b@a A{11.5}\or
   \b@b B\or\b@c C{5}\or\b@b D\or\b@b E\or\b@b F \or\b@c G{5}\or\b@b H\or
   \b@b I\or\b@c J{3}\or\b@b K\or\b@b L \or\b@b M\or\b@b N\or\b@c O{5} \or
   \b@b P\or\b@c Q{5}\or\b@b R\or\b@a S{8}\or\b@a T{10.5}\or\b@c U{5}\or
   \b@a V{12}\or\b@a W{16.5}\or\b@a X{11}\or\b@a Y{11.7}\or\b@a Z{7.5}\fi}}

\catcode`\X=11 \catcode`\@=12




\let\thischap\jobname

\def\partof#1{\csname returnthe#1part\endcsname}
\def\chapof#1{\csname returnthe#1chap\endcsname}

\def\setchapter#1,#2,#3;{%
  \expandafter\def\csname returnthe#1part\endcsname{#2}%
  \expandafter\def\csname returnthe#1chap\endcsname{#3}%
}

\setchapter 300a,A,II.A;
\setchapter 300b,A,II.B;
\setchapter 300c,A,II.C;
\setchapter 300d,A,II.D;
\setchapter 300e,A,II.E;
\setchapter 300f,A,II.F;
\setchapter 300g,A,II.G;
\setchapter  E53,B,N;
\setchapter  88r,B,I;
\setchapter  600,B,III;
\setchapter  705,B,IV;
\setchapter  734,B,V;

\def\cprefix#1{
\edef\theotherpart{\partof{#1}}\edef\theotherchap{\chapof{#1}}%
\ifx\theotherpart\thispart
   \ifx\theotherchap\thischap 
    \else 
     \theotherchap%
    \fi
   \else 
     \theotherchap\fi}

\def\sectioncite[#1]#2{%
     \cprefix{#2}#1}

\edef\thispart{\partof{\thischap}}
\edef\thischap{\chapof{\thischap}}

\def\lastpage of '#1' is #2.{\expandafter\def\csname lastpage#1\endcsname{#2}}


\def\spuriousreset{}


\expandafter\ifx\csname citeadd.tex\endcsname\relax
\expandafter\gdef\csname citeadd.tex\endcsname{}
\else \message{Hey!  Apparently you were trying to
\string\input{citeadd.tex} twice.   This does not make sense.} 
\errmessage{Please edit your file (probably \jobname.tex) and remove
any duplicate ``\string\input'' lines}\endinput\fi

\def\sciteu{\sciteerror{undefined}}

\def\sciteerror#1#2{{\mathortextbf{\scite{#2}}}\complainaboutcitation{#1}{#2}}
\def\mathortextbf#1{\hbox{\bf #1}}
\def\complainaboutcitation#1#2{%
\vadjust{\line{\llap{---$\!\!>$ }\qquad scite$\{$#2$\}$ #1\hfil}}}

\sectno=-1   
\localtags
\jjtags
\NoBlackBoxes
\define\mr{\medskip\roster}
\define\sn{\smallskip\noindent}
\define\mn{\medskip\noindent}
\define\bn{\bigskip\noindent}
\define\ub{\underbar}
\define\wilog{\text{without loss of generality}}
\define\ermn{\endroster\medskip\noindent}
\define\dbca{\dsize\bigcap}
\define\dbcu{\dsize\bigcup}
\define \nl{\newline}
\magnification=\magstep 1
\documentstyle{amsppt}

{    
\catcode`@11

\ifx\alicetwothousandloaded@\relax
  \endinput\else\global\let\alicetwothousandloaded@\relax\fi

\gdef\subjclass{\let\savedef@\subjclass
 \def\subjclass##1\endsubjclass{\let\subjclass\savedef@
   \toks@{\def\usualspace{{\rm\enspace}}\eightpoint}%
   \toks@@{##1\unskip.}%
   \edef\thesubjclass@{\the\toks@
     \frills@{{\noexpand\rm2000 {\noexpand\it Mathematics Subject
       Classification}.\noexpand\enspace}}%
     \the\toks@@}}%
  \nofrillscheck\subjclass}
} 


\expandafter\ifx\csname alice2jlem.tex\endcsname\relax
  \expandafter\xdef\csname alice2jlem.tex\endcsname{\the\catcode`@}
\else \message{Hey!  Apparently you were trying to
\string\input{alice2jlem.tex}  twice.   This does not make sense.}
\errmessage{Please edit your file (probably \jobname.tex) and remove
any duplicate ``\string\input'' lines}\endinput\fi

\expandafter\ifx\csname bib4plain.tex\endcsname\relax
  \expandafter\gdef\csname bib4plain.tex\endcsname{}
\else \message{Hey!  Apparently you were trying to \string\input
  bib4plain.tex twice.   This does not make sense.}
\errmessage{Please edit your file (probably \jobname.tex) and remove
any duplicate ``\string\input'' lines}\endinput\fi

\def\renewcommand{\newcommand}	       
\edef\cite{\the\catcode`@}%
\catcode`@ = 11
\let\@oldatcatcode = \cite
\chardef\@letter = 11
\chardef\@other = 12
%
%
%
%
\def\@innerdef#1#2{\edef#1{\expandafter\noexpand\csname #2\endcsname}}%
%
%
\@innerdef\@innernewcount{newcount}%
\@innerdef\@innernewdimen{newdimen}%
\@innerdef\@innernewif{newif}%
\@innerdef\@innernewwrite{newwrite}%
%
%
%
\def\@gobble#1{}%
%
%
%
\ifx\inputlineno\@undefined
   \let\@linenumber = \empty 
\else
   \def\@linenumber{\the\inputlineno:\space}%
\fi
%
%
%
\def\@futurenonspacelet#1{\def\cs{#1}%
   \afterassignment\@stepone\let\@nexttoken=
}%
\begingroup 
\def\\{\global\let\@stoken= }%
\\ 
\endgroup
\def\@stepone{\expandafter\futurelet\cs\@steptwo}%
\def\@steptwo{\expandafter\ifx\cs\@stoken\let\@@next=\@stepthree
   \else\let\@@next=\@nexttoken\fi \@@next}%
\def\@stepthree{\afterassignment\@stepone\let\@@next= }%
%
%
%
\def\@getoptionalarg#1{%
   \let\@optionaltemp = #1%
   \let\@optionalnext = \relax
   \@futurenonspacelet\@optionalnext\@bracketcheck
}%
%
%
\def\@bracketcheck{%
   \ifx [\@optionalnext
      \expandafter\@@getoptionalarg
   \else
      \let\@optionalarg = \empty
      \expandafter\@optionaltemp
   \fi
}%
\def\@@getoptionalarg[#1]{%
   \def\@optionalarg{#1}%
   \@optionaltemp
}%
%
%
%
\def\@nnil{\@nil}%
\def\@fornoop#1\@@#2#3{}%
\def\@for#1:=#2\do#3{%
   \edef\@fortmp{#2}%
   \ifx\@fortmp\empty \else
      \expandafter\@forloop#2,\@nil,\@nil\@@#1{#3}%
   \fi
}%
\def\@forloop#1,#2,#3\@@#4#5{\def#4{#1}\ifx #4\@nnil \else
       #5\def#4{#2}\ifx #4\@nnil \else#5\@iforloop #3\@@#4{#5}\fi\fi
}%
\def\@iforloop#1,#2\@@#3#4{\def#3{#1}\ifx #3\@nnil
       \let\@nextwhile=\@fornoop \else
      #4\relax\let\@nextwhile=\@iforloop\fi\@nextwhile#2\@@#3{#4}%
}%
%
%
%
\@innernewif\if@fileexists
\def\@testfileexistence{\@getoptionalarg\@finishtestfileexistence}%
\def\@finishtestfileexistence#1{%
   \begingroup
      \def\extension{#1}%
      \immediate\openin0 =
         \ifx\@optionalarg\empty\jobname\else\@optionalarg\fi
         \ifx\extension\empty \else .#1\fi
         \space
      \ifeof 0
         \global\@fileexistsfalse
      \else
         \global\@fileexiststrue
      \fi
      \immediate\closein0
   \endgroup
}%
%
%
%
%
\def\bibliographystyle#1{%
   \@readauxfile
   \@writeaux{\string\bibstyle{#1}}%
}%
\let\bibstyle = \@gobble
%
%
\let\bblfilebasename = \jobname
\def\bibliography#1{%
   \@readauxfile
   \@writeaux{\string\bibdata{#1}}%
   \@testfileexistence[\bblfilebasename]{bbl}%
   \if@fileexists
      \nobreak
      \@readbblfile
   \fi
}%
\let\bibdata = \@gobble
%
%
\def\nocite#1{%
   \@readauxfile
   \@writeaux{\string\citation{#1}}%
}%
\@innernewif\if@notfirstcitation
%
%
\def\cite{\@getoptionalarg\@cite}%
%
%
\def\@cite#1{%
   \let\@citenotetext = \@optionalarg
   \printcitestart
   \nocite{#1}%
   \@notfirstcitationfalse
   \@for \@citation :=#1\do
   {%
      \expandafter\@onecitation\@citation\@@
   }%
   \ifx\empty\@citenotetext\else
      \printcitenote{\@citenotetext}%
   \fi
   \printcitefinish
}%
\newif\ifweareinprivate
\weareinprivatetrue
\ifx\shlhetal\undefinedcontrolseq\weareinprivatefalse\fi
\ifx\shlhetal\relax\weareinprivatefalse\fi
\def\@onecitation#1\@@{%
   \if@notfirstcitation
      \printbetweencitations
   \fi
   \expandafter \ifx \csname\@citelabel{#1}\endcsname \relax
      \if@citewarning
         \message{\@linenumber Undefined citation `#1'.}%
      \fi
     \ifweareinprivate
      \expandafter\gdef\csname\@citelabel{#1}\endcsname{%
\strut 
\vadjust{\vskip-\dp\strutbox
\vbox to 0pt{\vss\parindent0cm \leftskip=\hsize 
\advance\leftskip3mm
\advance\hsize 4cm\strut\openup-4pt 
\rightskip 0cm plus 1cm minus 0.5cm ?  #1 ?\strut}}
         {\tt
            \escapechar = -1
            \nobreak\hskip0pt\pfeilsw
            \expandafter\string\csname#1\endcsname
             \pfeilso
            \nobreak\hskip0pt
         }%
      }%
     \else  
      \expandafter\gdef\csname\@citelabel{#1}\endcsname{%
            {\tt\expandafter\string\csname#1\endcsname}
      }%
     \fi  
   \fi
   \csname\@citelabel{#1}\endcsname
   \@notfirstcitationtrue
}%
%
%
\def\@citelabel#1{b@#1}%
%
%
\def\@citedef#1#2{\expandafter\gdef\csname\@citelabel{#1}\endcsname{#2}}%
%
%
%
\def\@readbblfile{%
   \ifx\@itemnum\@undefined
      \@innernewcount\@itemnum
   \fi
   \begingroup
      \def\begin##1##2{%
         \setbox0 = \hbox{\biblabelcontents{##2}}%
         \biblabelwidth = \wd0
      }%
      \def\end##1{}
      %
      %
      \@itemnum = 0
      \def\bibitem{\@getoptionalarg\@bibitem}%
      \def\@bibitem{%
         \ifx\@optionalarg\empty
            \expandafter\@numberedbibitem
         \else
            \expandafter\@alphabibitem
         \fi
      }%
      \def\@alphabibitem##1{%
         \expandafter \xdef\csname\@citelabel{##1}\endcsname {\@optionalarg}%
         \ifx\biblabelprecontents\@undefined
            \let\biblabelprecontents = \relax
         \fi
         \ifx\biblabelpostcontents\@undefined
            \let\biblabelpostcontents = \hss
         \fi
         \@finishbibitem{##1}%
      }%
      \def\@numberedbibitem##1{%
         \advance\@itemnum by 1
         \expandafter \xdef\csname\@citelabel{##1}\endcsname{\number\@itemnum}%
         \ifx\biblabelprecontents\@undefined
            \let\biblabelprecontents = \hss
         \fi
         \ifx\biblabelpostcontents\@undefined
            \let\biblabelpostcontents = \relax
         \fi
         \@finishbibitem{##1}%
      }%
      \def\@finishbibitem##1{%
         \biblabelprint{\csname\@citelabel{##1}\endcsname}%
         \@writeaux{\string\@citedef{##1}{\csname\@citelabel{##1}\endcsname}}%
         \ignorespaces
      }%
      %
      %
      \let\em = \bblem
      \let\newblock = \bblnewblock
      \let\sc = \bblsc
      \frenchspacing
      \clubpenalty = 4000 \widowpenalty = 4000
      \tolerance = 10000 \hfuzz = .5pt
      \everypar = {\hangindent = \biblabelwidth
                      \advance\hangindent by \biblabelextraspace}%
      \bblrm
      \parskip = 1.5ex plus .5ex minus .5ex
      \biblabelextraspace = .5em
      \bblhook
      \input \bblfilebasename.bbl
   \endgroup
}%
%
%
\@innernewdimen\biblabelwidth
\@innernewdimen\biblabelextraspace
%
%
%
\def\biblabelprint#1{%
   \noindent
   \hbox to \biblabelwidth{%
      \biblabelprecontents
      \biblabelcontents{#1}%
      \biblabelpostcontents
   }%
   \kern\biblabelextraspace
}%
%
%
%
\def\biblabelcontents#1{{\bblrm [#1]}}%
%
%
\def\bblrm{\rm}%
%
%
\def\bblem{\it}%
%
%
\def\bblsc{\ifx\@scfont\@undefined
              \font\@scfont = cmcsc10
           \fi
           \@scfont
}%
%
%
\def\bblnewblock{\hskip .11em plus .33em minus .07em }%
%
%
\let\bblhook = \empty
%
%
%
\def\printcitestart{[}
\def\printcitefinish{]}
\def\printbetweencitations{, }
\def\printcitenote#1{, #1}
%
%
%
\let\citation = \@gobble
%
%
%
\@innernewcount\@numparams
%
%
\def\newcommand#1{%
   \def\@commandname{#1}%
   \@getoptionalarg\@continuenewcommand
}%
%
%
\def\@continuenewcommand{%
   \@numparams = \ifx\@optionalarg\empty 0\else\@optionalarg \fi \relax
   \@newcommand
}%
%
%
\def\@newcommand#1{%
   \def\@startdef{\expandafter\edef\@commandname}%
   \ifnum\@numparams=0
      \let\@paramdef = \empty
   \else
      \ifnum\@numparams>9
         \errmessage{\the\@numparams\space is too many parameters}%
      \else
         \ifnum\@numparams<0
            \errmessage{\the\@numparams\space is too few parameters}%
         \else
            \edef\@paramdef{%
               \ifcase\@numparams
                  \empty  No arguments.
               \or ####1%
               \or ####1####2%
               \or ####1####2####3%
               \or ####1####2####3####4%
               \or ####1####2####3####4####5%
               \or ####1####2####3####4####5####6%
               \or ####1####2####3####4####5####6####7%
               \or ####1####2####3####4####5####6####7####8%
               \or ####1####2####3####4####5####6####7####8####9%
               \fi
            }%
         \fi
      \fi
   \fi
   \expandafter\@startdef\@paramdef{#1}%
}%
%
%
%
%
\def\@readauxfile{%
   \if@auxfiledone \else 
      \global\@auxfiledonetrue
      \@testfileexistence{aux}%
      \if@fileexists
         \begingroup
            \endlinechar = -1
            \catcode`@ = 11
            \input \jobname.aux
         \endgroup
      \else
         \message{\@undefinedmessage}%
         \global\@citewarningfalse
      \fi
      \immediate\openout\@auxfile = \jobname.aux
   \fi
}%
%
%
\newif\if@auxfiledone
\ifx\noauxfile\@undefined \else \@auxfiledonetrue\fi
%
%
%
%
\@innernewwrite\@auxfile
\def\@writeaux#1{\ifx\noauxfile\@undefined \write\@auxfile{#1}\fi}%
%
%
%
\ifx\@undefinedmessage\@undefined
   \def\@undefinedmessage{No .aux file; I won't give you warnings about
                          undefined citations.}%
\fi
%
%
\@innernewif\if@citewarning
\ifx\noauxfile\@undefined \@citewarningtrue\fi
%
%
%
\catcode`@ = \@oldatcatcode

\def\pfeilso{\leavevmode
            \vrule width 1pt height9pt depth 0pt\relax
           \vrule width 1pt height8.7pt depth 0pt\relax
           \vrule width 1pt height8.3pt depth 0pt\relax
           \vrule width 1pt height8.0pt depth 0pt\relax
           \vrule width 1pt height7.7pt depth 0pt\relax
            \vrule width 1pt height7.3pt depth 0pt\relax
            \vrule width 1pt height7.0pt depth 0pt\relax
            \vrule width 1pt height6.7pt depth 0pt\relax
            \vrule width 1pt height6.3pt depth 0pt\relax
            \vrule width 1pt height6.0pt depth 0pt\relax
            \vrule width 1pt height5.7pt depth 0pt\relax
            \vrule width 1pt height5.3pt depth 0pt\relax
            \vrule width 1pt height5.0pt depth 0pt\relax
            \vrule width 1pt height4.7pt depth 0pt\relax
            \vrule width 1pt height4.3pt depth 0pt\relax
            \vrule width 1pt height4.0pt depth 0pt\relax
            \vrule width 1pt height3.7pt depth 0pt\relax
            \vrule width 1pt height3.3pt depth 0pt\relax
            \vrule width 1pt height3.0pt depth 0pt\relax
            \vrule width 1pt height2.7pt depth 0pt\relax
            \vrule width 1pt height2.3pt depth 0pt\relax
            \vrule width 1pt height2.0pt depth 0pt\relax
            \vrule width 1pt height1.7pt depth 0pt\relax
            \vrule width 1pt height1.3pt depth 0pt\relax
            \vrule width 1pt height1.0pt depth 0pt\relax
            \vrule width 1pt height0.7pt depth 0pt\relax
            \vrule width 1pt height0.3pt depth 0pt\relax}

\def\pfeilsw{ \leavevmode 
            \vrule width 1pt height0.3pt depth 0pt\relax
            \vrule width 1pt height0.7pt depth 0pt\relax
            \vrule width 1pt height1.0pt depth 0pt\relax
            \vrule width 1pt height1.3pt depth 0pt\relax
            \vrule width 1pt height1.7pt depth 0pt\relax
            \vrule width 1pt height2.0pt depth 0pt\relax
            \vrule width 1pt height2.3pt depth 0pt\relax
            \vrule width 1pt height2.7pt depth 0pt\relax
            \vrule width 1pt height3.0pt depth 0pt\relax
            \vrule width 1pt height3.3pt depth 0pt\relax
            \vrule width 1pt height3.7pt depth 0pt\relax
            \vrule width 1pt height4.0pt depth 0pt\relax
            \vrule width 1pt height4.3pt depth 0pt\relax
            \vrule width 1pt height4.7pt depth 0pt\relax
            \vrule width 1pt height5.0pt depth 0pt\relax
            \vrule width 1pt height5.3pt depth 0pt\relax
            \vrule width 1pt height5.7pt depth 0pt\relax
            \vrule width 1pt height6.0pt depth 0pt\relax
            \vrule width 1pt height6.3pt depth 0pt\relax
            \vrule width 1pt height6.7pt depth 0pt\relax
            \vrule width 1pt height7.0pt depth 0pt\relax
            \vrule width 1pt height7.3pt depth 0pt\relax
            \vrule width 1pt height7.7pt depth 0pt\relax
            \vrule width 1pt height8.0pt depth 0pt\relax
            \vrule width 1pt height8.3pt depth 0pt\relax
            \vrule width 1pt height8.7pt depth 0pt\relax
            \vrule width 1pt height9pt depth 0pt\relax
      }


\def\widestnumber#1#2{}

\def\citewarning#1{\ifx\shlhetal\relax 
    \else
    \par{#1}\par
    \fi
}

\def\rm{\fam0 \tenrm}

\def\fakesubhead#1\endsubhead{\bigskip\noindent{\bf#1}\par}



%
%
%

%

\font\textrsfs=rsfs10
\font\scriptrsfs=rsfs7
\font\scriptscriptrsfs=rsfs5

\newfam\rsfsfam
\textfont\rsfsfam=\textrsfs
\scriptfont\rsfsfam=\scriptrsfs
\scriptscriptfont\rsfsfam=\scriptscriptrsfs

\edef\oldcatcodeofat{\the\catcode`\@}
\catcode`\@11

\def\Cal@@#1{\noaccents@ \fam \rsfsfam #1}

\catcode`\@\oldcatcodeofat


\expandafter\ifx \csname margininit\endcsname \relax\else\margininit\fi

\long\def\red#1\endred{}
\long\def\green#1\endgreen{}
\long\def\blue#1\endblue{}
\long\def\private#1\endprivate{}

\def\endred{ \unmatched endred! }
\def\endgreen{ \unmatched endgreen! }
\def\endblue{ \unmatched endblue! }
\def\endprivate{ \unmatched endprivate! }

\ifx\latexcolors\undefinedcs\def\latexcolors{}\fi

\def\emptycs{}
\def\evaluatelatexcolors{%
        \ifx\latexcolors\emptycs\else
        \expandafter\xxevaluate\latexcolors\xxfertig\evaluatelatexcolors\fi}
\def\xxevaluate#1,#2\xxfertig{\setupthiscolor{#1}%
        \def\latexcolors{#2}}


\font\smallfont=cmsl7
\def\rutgerscolor{\ifmmode\else\endgraf\fi\smallfont
\advance\leftskip0.5cm\relax}
\def\setupthiscolor#1{\edef\tmptmpcs{\noexpand\bgroup\noexpand\rutgerscolor
\noexpand\def\noexpand\currentcolor{#1}%
\noexpand}%
\expandafter\let\csname#1\endcsname\tmptmpcs
\def\tmptmpcs{\checkColorUnmatched{#1}\popthecolor}
\expandafter\let\csname end#1\endcsname\tmptmpcs}

\def\checkColorUnmatched#1{\def\expectcolor{#1}%
    \ifx\expectcolor\currentcolor   
    \else \edef\failhere{\noexpand\tryingToClose '\currentcolor' with end\expectcolor}\failhere\fi}

\def\currentcolor{???}

\def\popthecolor{\ifmmode\else\endgraf\fi\egroup}

\expandafter\def\csname#1\endcsname{}

\evaluatelatexcolors

 \let\outerhead\head
 \def\head{\innerhead}
 \let\innerhead\outerhead

 \let\outersubhead\subhead
 \def\subhead{\innersubhead}
 \let\innersubhead\outersubhead

 \let\outersubsubhead\subsubhead
 \def\subsubhead{\innersubsubhead}
 \let\innersubsubhead\outersubsubhead

 \let\outerproclaim\proclaim
 \def\proclaim{\innerproclaim}
 \let\innerproclaim\outerproclaim

 %
 %
 %
 %

\def\demo#1{\medskip\noindent{\it #1.\/}}
\def\enddemo{\smallskip}

\def\remark#1{\medskip\noindent{\it #1.\/}}
\def\endremark{\smallskip}

\pageheight{8.5truein}
\topmatter
\title{Polish Algebras, shy from freedom } \endtitle
\author {Saharon Shelah \thanks {\null\newline 
This research was partially supported by the Israel Science Foundation
founded by the Israel Academy of Sciences and
Humanities. Publication 771. \null\newline
I would like to thank Alice Leonhardt for the beautiful typing. 
} \endthanks} \endauthor 

\affil{The Hebrew University of Jerusalem \\
Einstein Institute of Mathematics \\
Edmond J. Safra Campus, Givat Ram \\
Jerusalem 91904, Israel
 \medskip
 Department of Mathematics \\
 Hill Center-Busch Campus \\
  Rutgers, The State University of New Jersey \\
 110 Frelinghuysen Road \\
 Piscataway, NJ 08854-8019 USA} \endaffil
\endtopmatter
\document

\newpage

\head {Annotated Content} \endhead
 \spuriousreset
\bn
\S0 $\quad$ Introduction
\bn
\S1 $\quad$ Metric groups and metric models
\mr
\item "{${{}}$}"  [We give basic definitions and some relations.]
\endroster
\bn
\S2 $\quad$ Semi-metric groups: on automorphism groups of uncountable
structures
\mr
\item "{${{}}$}"  [We prepare the ground to treating the automorphism
group of a structure of cardinality strong limit of countable
cofinality, e.g. $\beth_\omega$; this is the ``semi".  We also
consider replacing ``the automorphism group of ..." by other derived
structures.]
\endroster
\bn
\S3 $\quad$ Compactness of metric algebras
\mr
\item "{${{}}$}"  [The main lemma gives a sufficient condition for
solvability of a set of equations of some form.  We then deduce
sufficient conditions for non-freeness.]
\endroster
\bn
\S4 $\quad$  Conclusion
\mr
\item "{${{}}$}"  [We show that Polish groups are not free, also even
the semi-metric vesion suffice.  We then derive the conclusion on
automorphism groups.]
\endroster
\bn
\S5 $\quad$  Quite free but not free abelian groups
\mr
\item "{${{}}$}"  [We show that for every $n$ there is a very
explicit definition of an abelian group (Borel and even $F_\sigma$)
which is free if the continuum is at least $\aleph_{n+1}$.  This
applies to group and large family of varieties (families of
algebra defined by a set of equations).  We note that this works for
subgroups of $\Bbb Z^\omega$.]
\endroster
\bn
\S6 $\quad$ Beginning of stability theory
\mr
\item "{${{}}$}"  [We prove some basic results.]
\endroster
\newpage

\head {\S0 Introduction} \endhead  \resetall \sectno=0
 \spuriousreset
\bigskip

Our first motivation was the question:  can a countable structure have an
automorphism group, which a free uncountable group?  This is answered
negatively in \cite{Sh:744}.

This was a well known problem in group theory at least in England
(David Evans in a meeting in Durham 1987)
and we thank Simon Thomas for telling us about it.  
Independently in descriptive set theory, Howard and Kechris
\cite{BeKe96} ask if
there is an uncountable free Polish group, i.e. which is on a complete
separable metric space.  A related result (before \cite{Sh:744}) 
was gotten by Solecki [3] who proved that the group of
automorphisms of a countable structure cannot be an uncountable free
abelian group.  Having the problem arise 
independently supported the feeling that it is a natural problem.

The idea of the proof in \cite{Sh:744} was to prove that such a group has some
strong algebraic completeness or compactness, 
more specifically for any sequence
$\langle d_n:n < \omega \rangle$ of elements of the group converging
to the 
identity many countable sets of equations are solvable.  This is parallel
in some sense to Hensel lemma for the $p$-adics, and seem to me
interesting in its own right.

Lecturing in a conference in Rutgers, February 2001, I was asked
whether I am really speaking on Polish groups.
We can prove this using a more restrictive condition on the set of equations.
Parallel theorems, e.g. holds for semi groups and for metric algebras,
e.g. with non isolated unit ($e$ is a unit means $\{e\}$ is a subalgebra).
Here we do the general case.
\bn
More specifically we prove (see Conclusion \scite{c.2}, \scite{c.1}(1)).

\proclaim{\stag{0.1} Theorem}  1) There is no Polish group which as a
group is free and uncountable. \nl
2) Slightly more generally, assume
\mr
\item "{$(a)$}"  $G$ is a metric space
\sn
\item "{$(b)$}"  $G$ is a group with continuous $xy,x^{-1}$
\sn
\item "{$(c)$}"  $G$ is complete
\sn
\item "{$(d)$}"  the density of $G$ is $<|G|$.
\ermn
\ub{Then} $G$ is not free.
\endproclaim
\bn
\margintag{0.2}\ub{\stag{0.2} Thesis}:  If $G$ is a Polish algebra satisfying one of the
compactness conditions 
defined below, \ub{then} it is in fact large in the sense of
lots of sets of equations has a solution. 

A reader interested just in this theorem can read just \S3.
\bn
\margintag{0.3}\ub{\stag{0.3} Question}:  What are the restrictions on Aut$(\Bbb A)$ 
for uncountable structures $\Bbb A$?

We also prove that if $\Bbb A$ is a structure of cardinality $\mu$,
$\mu$ is strong limit of cofinality $\aleph_0$ (e.g. $\beth_\omega$)
and the automorphism
group of $\Bbb A$ is of cardinality $> \mu$, \ub{then} it is far from
being free; this does not follow directly from \scite{0.2}(2) as the
natural metric considered here does not satisfy all the conditions.

In \cite{Sh:744} this is proved in the special case where
$\Bbb A = \dbcu_n P^{\Bbb A}_n$ satisfying $|P_n^{\Bbb A}| < \mu$.
\bn
\centerline{ $* \qquad * \qquad *$}
\bn
\ub{Note}:  An arbitrary subgroup e.g. of the symmetric group of size
$\aleph_1$ can consistently be made into an automorphism group by
Just, Shelah, Thomas \cite{JShT:654}.  So $\beth_\omega$ is
more interesting from this point of view.
\bn
\ub{Ideas regarding Aut$(\Bbb A),\|\Bbb A\| = \beth_\omega$}:  
Let the set of elements of $\Bbb A$
be $\beth_\omega$.  For $f,g \in \text{ Aut}(\Bbb A)$ let

$$
d(f,g) = \sup\{2^{-n}:f \restriction \beth_n \ne g \restriction
\beth_n \text{ or } f^{-1} \restriction \beth_n \ne g^{-1}
\restriction \beth_n\}.
$$
\mn
Again a complete metric space and a Polish group.  But we need to
prove that there is a non-trivial convergent sequence.

This space is not necessarily of density $\le \beth_\omega$.  A way
out is to define

$$
\align
d_1(f,g) = \sup\{2^{-n}:&(\exists \alpha < \beth_n) \dsize \bigvee_m
(f(\alpha) < \beth_n \neq\equiv g(\alpha) < \beth_m \\
  &\text{ or } f(\alpha) \ne g(\alpha) \text{ both } < \beth_n \\
  &\text{ or similarly with } g^{-1},g^{-1}\}.
\endalign
$$
\mn
The purpose of introducing $d_1$ is to decrease the density.  Now
``$d_1(f,g) < 2^{-n}$" is an e.g. relation with $\le 2^{\beth_n}$
classes (rather than $(\beth_\omega)^{\beth_n}$ classes).  This is a
complete metric space with density $\le \beth_\omega$.  However, there are
problems with ``$(d(f_1,g_1),d(f_2,g_2) < 2^{-n} \Rightarrow d(f_1
\circ f_2,g_1 \circ g_2) < 2^{-n+1}$".

So we start the proof with $d_1$, find a non trivial converging
sequence $\langle f_n:n < \omega \rangle$ in $d_1$, prove it converges
pointwise, replace $(\beth_n:u < \omega)$ by $\langle A_n:n < \omega
\rangle$ and get a sequence converging in $d$.
\bn
\margintag{0.4}\ub{\stag{0.4} Question}:  Is there a model theory of Polish spaces?

Naturally we would like to develop a parallel to classification theory
(see \cite{Sh:c}).  A natural test problem is to 
generalize ``Morley theorem = Los conjecture". 
\ub{But} we only have one model so what is the meaning? \nl
Well, we may change the universe.  If we deal with abelian groups (or
any variety) it is probably more natural to ask when is such (Borel)
algebra free.
\bn
\margintag{0.5}\ub{\stag{0.5} Example}:  If 
$\Bbb P$ is adding $(2^{\aleph_0})^+$-Cohen subsets of $\omega$ then

$$
(\Bbb C)^{\bold V} \text{ and } (\Bbb C)^{\bold V[G]}
$$
\mn
are both algebraically closed fields of characteristic $0$ which are
not isomorphic (as they have different cardinalities).

So we restrict ourselves to forcing $\Bbb P_1 \lessdot \Bbb P_2$ such that

$$
(2^{\aleph_0})^{\bold V[\Bbb P_1]} = (2^{\aleph_0})^{\bold V[\Bbb P_1]}
$$
\mn
and compare the Polish models in $\bold V^{{\Bbb P}_1},
\bold V^{{\Bbb P}_2}$.    We may restrict our forcing notions to
c.c.c. or whatever...
\bn
\margintag{0.6}\ub{\stag{0.6} Example}:  Under any such interpretation
\mr
\item "{$(a)$}"  $\Bbb C =$ the field of complex numbers is categorical
\sn
\item "{$(b)$}"  $\Bbb R =$ the field of the reals is not \nl
(by adding $2^{\aleph_0}$ many Cohen reals).
\ermn
(Why?  Trivially: $\Bbb R^{\bold V[{\Bbb P}_2]}$ is complete
in $\bold V[\Bbb P_2]$ while $\Bbb R^{\bold V[\Bbb P_1]}$ in $\bold V^{[{\Bbb
P}_2]}$ is not complete but there are less trivial reasons).
\bn
\margintag{0.7}\ub{\stag{0.7} Conjecture}:  We have a dychotomy, i.e. either the model
is similar to categorical theories, or there are ``many complicated models".
\bn
\margintag{0.8}\ub{\stag{0.8} Thesis}:  Classification theory for such models resemble
more the case of $\Bbb L_{\omega_1,\omega}$ than the first order.
\bn
See \cite{Sh:h}; as support for this thesis we prove:
\proclaim{\stag{0.9} Theorem}  There is an $F_\sigma$ abelian group (i.e. a
$F_\sigma$-definition, in fact an explicit definition) such that
$\bold V \models ``G$ is a free abelian group" iff $\bold V \models
2^{\aleph_0} < \aleph_{736}$.
\endproclaim
\bn
\ub{Comments}:  In the context of the previous theorem we cannot do
better than $F_\sigma$, but we may hope for some other example which
is not a group or categoricity is not because of freeness.
\bn
\margintag{0.10}\ub{\stag{0.10} Conclusion}:  Freeness (of an $F_\sigma$-abelian
group) can stop at $\aleph_n$ (any $n$).
\bn
A connection with the model theories is that by 
Hart-Shelah \cite{HaSh:323} such things can also occur
in $L_{\omega_1,\omega}$ whereas (by \cite{Sh:87a}, \cite{Sh:87b} Theorem)
if $\dsize \bigwedge_n (2^{\aleph_n} < 2^{\aleph_{n+1}})$ and
$\psi \in \Bbb L_{\omega_1,\omega}$, categorical in every
$\aleph_n$, \ub{then} $\psi$ is categorical in every $\lambda$.  See
more in \cite{ShVi:648}.

The parallels here are still open.

This casts some light on the thesis that non-first order logics are
``more distant" from the ``so-called" mainstream mathematics.

Returning to stability theory per-se we have the modest:
\bigskip

\proclaim{\stag{0.11} Theorem}  For ``$\aleph_0$-stable Borel models" 
the theorem on the existence of indiscernibles can be generalized.
\endproclaim
\bn
We may consider another version of the interpretation of
``categoricity".   Of course, we can
use more liberal than $\bold L(A_2,r)$ or restrict the
$A_\ell$'s further (as in the forcing version).
\definition{\stag{0.12} Definition}  1) We say that $\Bbb A$ is
categorical in $\lambda \le 2^{\aleph_0}$ if: for some real $r$: for
every $A_1,A_2 \subseteq \lambda$ the models $\Bbb A^{\bold
L[A_1,r]},\Bbb A^{\bold L[A_2,r]}$ are isomorphic (in $\bold V$). \nl
2) For a class ${\frak K}$ of forcing notions and cardinal $\lambda$
such that for at least one $\Bbb P \in {\frak K},\Vdash_{\Bbb P}
``2^{\aleph_0} \ge \lambda"$, we have in $\bold V^{\Bbb P}$: the
structure $\Bbb A$ is categorical in $\lambda$ in the sense of part (1).
\sn

Comparing Definition \scite{0.12}(1) with 
the forcing version we lose when $\bold V = \bold L$, as
it says nothing, we gain as (when $2^{\aleph_0} > \aleph_1$) we do not
have to go outside the universe.  
Maybe best is categorical in $\lambda$ in $\bold
V^{\Bbb P}$ for every c.c.c. forcing notion $\Bbb P$ making
$2^{\aleph_0} \ge \lambda$.

Note also that it may be advisale to restrict ourselves to the case
$\lambda$ is regular as we certainly like to avoid the possibility
$(2^{\aleph_0})^{\bold L[A_1,r]} = \lambda < (2^{\aleph_0})^{\bold
L[A_2,r]}$ (see on this \cite[VII]{Sh:g}).
\enddefinition
\bigskip

Of course, any reasonable definition of unstability implies
non-categoricity: if we have many types we should have a perfect set
of them, hence adding Cohen subsets of $\omega$ adds more types
realized.  If we add $\langle \eta_i:i < 2^{\aleph_0} \rangle$ Cohen
reals for every $A \subseteq 2^{\aleph_0},\langle M^{[\eta_i:i \in
A]}:A \subseteq 2^{\aleph_0} \rangle$ are non-isomorphic over the
countable set of parameters, if we get $2^{2^{\aleph_0}}$
non-isomorphic models, we can forget the parameters and retain our
``richness in models".

Lately Blass asks on definable abelian subgroups of
$\Bbb Z^\omega$, answers are derived for this from \cite{Sh:402} and
\S5.   We may be more
humble than in \scite{0.4}.
\bn
\margintag{0.13}\ub{\stag{0.13} Question}:  Is there model theory for equational
theories, stressing free algebras? \nl
The material in \S1 - \S5 (except some generalizations) was presented
in a course in Rutgers, Sept. - Oct 2001 and I thank the audience for
their comments.  We shall continue elsewhere.
\bigskip

\demo{\stag{0.21} Notation}  1) Let $\omega$ denote the set of natural
numbers, and let $x < \omega$ mean ``$x$ is a natural number". \nl
2) Let $a,b,c,d$ denote members of $G$ (a group). \nl
3) Let $\bar d$ denote a finite sequence $\langle d_n:n < n^*
\rangle$,  and similarly in other cases. \nl
4) Let $k,\ell,m,n,i,j,r,s,t$ denote natural numbers.
\nl
5) Let $\varepsilon,\zeta,\xi$ denote reals $>0$.
\enddemo
\bigskip

\definition{\stag{0.22} Definition}  1) A group word is a sequence
$x^{r_1}_1x^{r_2}_2,\dotsc,x^{r_k}_k$ where the $x_\ell$ are variables or
elements of a group and $n_\ell \in \Bbb Z$ for $\ell
=1,\dotsc,k$. \nl
2) The word is reduced if $n_\ell \ne 0,x_\ell \ne x_{\ell +1}$. \nl
3) The length of a word $w = x^{n_1}_1 x^{n_2}_2 \ldots x^{n_k}_k$ is
$\dsize \sum^k_{\ell=1}|n_\ell|$. \nl
4) A group term $w(x_1,\dotsc,x_n)$ is a word of the form
$x^{r_1}_{\ell_1} x^{r_2}_{\ell_2} \ldots x^{r_k}_{i_k}$ with $i_\ell
\in \{1,\dotsc,k\},r_\ell \in \Bbb Z$ (actually $r_\ell \in
\{1,\dotsc,-1\}$ suffice).  For a group $G$ and $a_1,\dotsc,a_k \in
G$, the meaning of $b = w(a_1,\dotsc,a_k) \in G$ should be clear.
\enddefinition
\newpage

\head {\S1 metric groups and metric models} \endhead  \resetall \sectno=1
 \spuriousreset
\bigskip

We first define [semi]-metric [complete] group, and give a natural
major example: automorphism groups.  The natural example of
semi-metric group is the semi-group of endomorphism of a countable structure.
\definition{\stag{1.1} Definition}   $(G,{\frak d})$ is called a
metric group if:
\mr
\item "{$(a)$}"  $G$ is a group
\sn
\item "{$(b)$}"  $G$ is a metric space for the metric ${\frak d}$
\sn
\item "{$(c)$}"  the functions $xy,x^{-1}$ are continuous.
\ermn
2) $(G,{\frak d})$ is a metric semi-group when
\mr
\item "{$(a)$}"  $G$ is a semi-group
\sn
\item "{$(b)$}"  $G$ is a metric space for the metric ${\frak d}$
\sn
\item "{$(c)$}"  the function $xy$ is continuous (there is no $x^2$ as
$G$ is just a semi-group).
\ermn
3) Saying $(G,{\frak d})$ is complete, means complete as a metric space.
\enddefinition
\bn
\margintag{1.1A}\ub{\stag{1.1A} Notation}:  1) For a metric group $\bold M$ the 
metric is denoted by ${\frak d}_{\bold M}$ and the unit is denoted by
$e_{\bold M}$ and the group by $G_{\bold M}$.  When no confusion
arises ``$G$ is a metric group" means $(G,{\frak d}_G)$ is a metric
group.
\nl
2) Similarly for semi-groups.
\bn
Now we define cases closer to automorphism groups, in those cases the
proof is very similar to the one in \cite{Sh:744}.
\definition{\stag{1.2} Definition}  1) $G$ is a specially metric group \ub{if}: 
\mr
\item "{$(\alpha)$}"  $G$ is a metric group
\sn
\item "{$(\beta)$}"  for every $\zeta,\varepsilon \in \Bbb R^+$ there
is $\xi \in \Bbb R^+$ such that: if 
$x_1,x_2,y_1,y_2 \in \{x:{\frak d}_G(x,e_G) < \varepsilon\}$ 
\ub{then} ${\frak d}_G(x_1,x_2) < \xi \wedge {\frak
d}(y_1,y_2) < \xi$ implies ${\frak d}_G(x_1y_1,x_2y_2) < \zeta \wedge
{\frak d}(x^{-1}_1,x^{-1}_2) < \zeta$; this is a kind of uniform
continuity (inside the $\varepsilon$-neighborhood of 
$e_G$; this is harder if we
increase $\varepsilon$ and/or decrease $\zeta$
\sn
\item "{$(\gamma)$}"  for arbitrarily small $\zeta \in \Bbb R^+$ the
set $\{a \in G:{\frak d}(a,e_G) < \zeta\}$ is a subgroup of $G$.
\ermn
2) We say $\bar \zeta = \langle \zeta_n:n < \omega \rangle$ is
strongly O.K. for $G$ if:
\mr
\item "{$(a)$}"  $\zeta_n \in \Bbb R^+_n$
\sn
\item "{$(b)$}"  $\zeta_n$ satisfies clause $(\gamma)$ of part (2),
i.e. $\{a \in G:{\frak d}(a,e_G) < \zeta_n\}$ is a subgroup of $G$
\sn
\item "{$(c)$}"  $\zeta_{n+1} \le \zeta_n$ and $0 = 
\text{ inf}\{\zeta_n:n < \omega\}$
\sn
\item "{$(d)$}"  if $x_1,x_2,y_1,y_2 \in \{a \in G:{\frak d}_G(a,e_G) <
\zeta_0\}$ and ${\frak d}(x_1,x_2) < \zeta_{n+1},{\frak d}_G(y_1,y_2) <
\zeta_{n+1}$ and $r(1),r(2) \in \{1,-1\}$ \ub{then} ${\frak d}_G
(x^{r(1)}_1y^{r(2)}_1,x^{r(1)}_2y^{r(2)}_2) < \zeta_n$.
\ermn
3) We say $G$ is specially$^+$ metric group if in part (1) we have
$(\alpha), (\beta)$ and
\mr
\item "{$(\gamma)^+$}"  for every $\zeta \in \Bbb R^+$ the set $\{a
\in G:{\frak d}(a,e_G) < \zeta\}$ is a subgroup of $G$.
\ermn
4) We define similarly for semi groups omitting the operation $x^{-1}$
(this means omitting ``${\frak d}(x^{-1}_1,x^{-1}_2) < \zeta$" in
 clause $(\beta)$ of (1), and demanding $r(1),r(2)=1$ in clause (d) of
part (2).
\enddefinition
\bigskip

\demo{\stag{1.3} Observation}  1) If $G$ is a special metric group
\ub{then} there is a sequence $\bar \zeta$ which is strongly O.K. for $G$. \nl
2) We can in clause (d) of \scite{1.2}(2) above omit $r(1),r(2)$ and conclude
only ``${\frak d}(x_1y_1,x_2y_2) < \zeta_n$ and ${\frak
d}(x^{-1}_1,x^{-1}_2) < \zeta_n$.  This causes just slight changes
in the computations of length in the proof or replacing $\bar \zeta$
by a suitable subsequence. \nl
3) Every specially$^+$ metric group is a specially metric group. \nl
4) Parts (1), (3) holds for semi groups, too.
\enddemo
\bigskip

\demo{Proof}  Easy.
\enddemo
\bigskip

\definition{\stag{1.4} Definition}  1) Assume $\Bbb A$ is a countable
structure with automorphism group $G = \text{ Aut}(\Bbb A)$ and for notational
simplicity its set of elements is $\omega$, the set of natural numbers
(and, of course, it is infinite, otherwise trivial).

We define a metric ${\frak d} = {\frak d}_{\Bbb A} = {\frak
d}^{\text{aut}}_{\Bbb A}$ on $G$ by

$$
{\frak d}(f,g) = \text{ sup}\{2^{-n}:f(n) \ne g(n) \text{ or }
f^{-1}(n) \ne g^{-1}(n)\}.
$$
\mn
Let Aut$_{\Bbb A} = (\text{Aut}(\Bbb A),{\frak d}^{\text{aut}}_{\Bbb
A})$, but we may write $G_{\Bbb A}$ or $G^{\text{aut}}_{\Bbb A}$. \nl
2) Assume $\Bbb A$ is a countable 
(infinite) structure, let End$(\Bbb A)$ be the
semi-group of endomorphisms of $\Bbb A$; assume for simplicity that
its set of elements is $\omega$ and let ${\frak d}^{\text{end}}_{\Bbb A}$ be 
the following metric on End$(\Bbb A)$

$$
{\frak d}(f,g) = \text{ sup}\{2^{-n}:f(n) \ne g(n)\}.
$$
\mn
Let End$_{\Bbb A}$ be (End$(\Bbb A),{\frak d}^{\text{end}}_{\Bbb A})$,
we may write $G^{\text{end}}_{\Bbb A}$. \nl
3) If in part (2) we restrict ourselves to monomorphisms, we write
Mon$(\Bbb A)$, Mon$_{\Bbb A},{\frak d}^{\text{mon}}_{\Bbb A}$.  
\enddefinition
\bigskip

\proclaim{\stag{1.5} Claim}:  For $\Bbb A$ as above: \nl
1)  $(\text{\rm Aut}_{\Bbb A},{\frak d}^{\text{aut}}_{\Bbb A})$ is a complete 
separable specially$^+$ metric group. 
\nl
2) ({\rm End}$(\Bbb A),{\frak d}^{\text{end}}_{\Bbb A})$ and 
({\rm Mono}$(\Bbb A),{\frak d}^{\text{mono}}_{\Bbb A})$ are 
complete separable specially$^+$ metric semi-groups.
\endproclaim
\bigskip

\demo{Proof}  Easy.
\enddemo
\bn
We may think of a more general context.
\definition{\stag{1.6} Definition}  1) We say ${\frak a}$ is a metric
$\tau$-model ($\tau = \tau_{\frak a}$ is a vocabulary, that is a set
of function symbols and predicates; in the main case we say 
$\tau$-algebra when $\tau$ has no predicates only 
functions) \ub{if}
\mr
\item "{$(a)$}"  ${\frak a}$ 
is a metric space with metric ${\frak d}_{\frak a}$
\sn
\item "{$(b)$}"  $M_{\frak a} = M({\frak a})$ is a model or an
algebra with universe $|M_{\frak a}|$, 
(of course with a set of elements the same as the set of
points of the  metric space), with vocabulary $\tau = \tau_{\frak a}$
\sn
\item "{$(c)$}"  if $F \in \tau_{\frak a}$ is (an $n$-place) function symbol, \ub{then}
$F^{M({\frak a})}$ is (an $n$-place) continuous function from
$M_{\frak a}$ to $M_{\frak a}$ (for ${\frak d}_{\frak a}$, i.e. by the
topology which the metric ${\frak d}_G$ induces, of course)
\sn
\item "{$(d)$}"  if $R \in \tau_{\frak a}$ is an $n$-place predicate,
\ub{then} $R^{\frak a} = R^{M({\frak a})}$ is a closed subset of 
$M^n_{\frak a} = {}^n(|M_{\frak a}|)$.
\ermn
2) We say ${\frak a}$ is unitary if some $e \in \tau_{\frak a}$ is a unit
of ${\frak a}$ which means that $e$ an individual constant
and $\{e_{\frak a}\}$ is closed under $F^{M({\frak a})}$ for $F \in
\tau_G$ and $R \in \tau_G \Rightarrow
\langle e,\ldots \rangle \in R^{\frak a}$. \nl
3) We say ${\frak a}$ is complete if $(|M_{\frak a}|,{\frak d}_{\frak
a})$ is a complete metric space. \nl
4) We replace ``unitary" by ``specially$^\pm$ unitary" above if we add:
\mr
\item "{$(*)$}"  for every $\zeta \in \Bbb R^+$, the set $\{a \in
M_{\frak a}:{\frak d}_{\frak a}(a,e_{\frak a}) < \zeta\}$ is a
subalgebra of $\Bbb A$.  
\ermn
5) We replace unitary by ``specially-unitary" if some $\bar \zeta$
witness it which means:
\mr
\item "{$(a)$}"   $\bar \zeta$ is a decreasing sequence of positive
reals with limit zero
\sn
\item "{$(b)$}"  for every $F \in \tau_{\frak a}$ for some $n =
n(F,{\frak a})$ we have for every $m \in [n,\omega)$ the set $\{a \in
M_{\frak a}:{\frak d}_{\frak a}(a,e_{\frak a}) < \zeta_m\}$ is closed
under $F^{M({\frak a})}$.
\ermn
6) We add the adjective partial if we allow $F^{M({\frak a})}$ to be a
partial function, so in clause (c) of part (1) means now:
\mr
\item "{$(c)$}"   if $F \in \tau_{\frak a} = \tau(M_{\frak a})$ is an
$n$-place function symbol then the set 
$\{(F^{M({\frak a})}(a_1,\dotsc,a_n):
a_1,\dotsc,a_n \in M_{\frak a}$ and $F^{M({\frak a})}
(a_1,\dotsc,a_n)$ is 
well defined$\}$ is a closed subset of ${}^{n+1}(M_{\frak a})$.
\ermn
7) We say ${\frak a}$ is specially-unitary \ub{if} some $\bar \zeta$
witnesses it which means (a),(b) from part (5) and
\mr
\item "{$(d)$}"  for any $F \in \tau_{\frak a}$ a $k$-place function
for every $m \ge n(F,{\frak a})$ we also have:
\nl
if ${\frak d}_{\frak a}(x_\ell,y_\ell) < \zeta_{m+1}$ for $\ell =
1,\dotsc,k$ and ${\frak d}_{\frak a}(x_\ell,e_{\frak a}) <
\zeta_{m+1},{\frak d}_{\frak a}(y_\ell,e_{\frak a}) < \zeta_{m+1}$ then
${\frak d}_{\frak a}(F^{M({\frak a})}(x_1,\dotsc,x_k),F^{M({\frak
a})}(y_1,\dotsc,y_\ell))$ is $< \zeta_m$.
\endroster
\enddefinition
\newpage

\head {\S2 Semi metric groups: on automorphism groups of uncountable
structures} \endhead  \resetall 
 \spuriousreset
\bigskip

It seems natural investigating
the automorphism groups of a model $\Bbb A$ say of cardinality, e.g.
$\beth_\omega$, intending to put in a framework where we shall be able
to prove it is not a free group of
cardinality $> \beth_\omega$. Now choose $\bar A = \langle A_n:n <
\omega\rangle$ such that, e.g. the universe of $\Bbb A$ is
$\beth_\omega$ and $A_n = \beth_n$.  For such $\bar A$ there is a 
natural metric on Aut$(\Bbb A)$ under which it is a 
complete metric group, \ub{but} usually of
too big density; there is another natural metric on Aut$(\Bbb A)$
under which it is a complete metric space with density $\beth_\omega$
\ub{but} multiplication is not continuous and Cauchi sequences may
not converge to any point.  To get the desired results we use
``semi-metric" defined in \scite{1.7}, which combine the
two metrics; in other words we weaken the completeness demand; 
hopefully this will have other applications as well, e.g. also for
Borel groups.  We also look at some generalizations: replacing
Aut$(\Bbb A)$ by other derived structures.
\bn
The reader may skip this section if not interested in the results
concerning $\beth_\omega$.  One way to present what we are doing is
\definition{\stag{1.7} Definition}   1) We say $G$ is an indirectly complete
 metric group if:
\mr
\item "{$(a)$}"  $G$ is a group
\sn
\item "{$(b)$}"  $G$ is a metric space under ${\frak d}_G$
\sn
\item "{$(c)$}"  if $\bar c = \langle c_n:n < \omega \rangle$
satisfies ${\frak d}_G(c_n,c_{n+1}) < 1/2^n$,
\ub{then} letting $d_n = c^{-1}_{2n} c_{2n+1}$ and $\bar d =
\langle d_n:n < \omega \rangle$ we have
{\roster
\itemitem{ $(*)$ }  for some metric ${\frak d}' = {\frak d}_{G,\bar
d}$, under ${\frak d}',G$ is a metric group, and
$\bar d$ converges to some $c$ under ${\frak d}'$, see below.
\endroster}
\ermn
1A) We say $G$ is an indirectly complete metric group and is defined
similarly only in $(*)$ we replace ``metric" by ``complete metric".
\nl
1B) Similarly for add another adjective or several of them.
\nl
2) Similarly for semi-groups and for algebras (and models, see
Definition \scite{1.6}) with one change: in $(*)$
of clause (c) we let $d_n = c_n$.
\enddefinition
\bn
\margintag{1.7A}\ub{\stag{1.7A} Discussion}:  1) In \scite{1.7}(1), $(*)$, clearly $c=e_G$.
 So $G$ in Definition \scite{1.7}(1)
is not necessarily a metric group,
i.e. product and inverse are not necessarily continuous.  Also note
that a metric group which is a complete metric space is not
necessarily an indirectly metric group. 
\nl
2) For groups $G$ in \scite{1.7}(1), clause (c), we can 
replace ``${\frak d}(c_n,c_{n+1}) < 1/2^n$" by $\langle c^{-1}_{2n} 
c_{2n+1}:n < \omega \rangle$ converge to $e_G$ with no significant
difference in the results. 
\nl
3) There are some other variants which can serve as well:
we can just demand ``$\bar c$ is a Cauchi sequence" and/or add ``$\bar
c$ converge to $c'$", and/or in the conclusion say ``some $\omega$-subsequence
$\bar d'$ of $\bar d$ converge to some $d$".
\nl
4) Similarly in \scite{1.7}(2).
\bigskip

\definition{\stag{1.8} Definition}  Assume $\Bbb A$ is a structure. \nl
1) $\bar A$ is an $\omega$-representation of $\Bbb A$ if $\bar A = \langle
A_n:n < \omega \rangle,A_n \subseteq A_{n+1}$ for $n < \omega$ and
$\cup\{A_n:n < \omega\}$ is the universe of $\Bbb A$. \nl
2) For every $\omega$-representation $\bar A$ of $\Bbb A$ 
let Aut$_{\bar A}(\Bbb A) = \{f \in \text{ Aut}(\Bbb A)$: 
for every $n < \omega$ for some $m < \omega$ we have
$(\forall x \in A_n)(f(x) \in A_m \and f^{-1}(x) \in A_m)\}$. 
\nl
3) If $\bar A$ is an $\omega$-representation of $\Bbb A$ and $G = \text{
Aut}(\Bbb A)$ \ub{then} we define ${\frak d} = {\frak d}_{\Bbb A,\bar A} =
{\frak d}^{\text{aut}}_{{\frak A},\bar A}$, a metric \footnote{this is proved
in \scite{1.9}} on $G$ by

$$
\align
{\frak d}(f,g) = \text{ sup}\{2^{-n}:&\text{ there is } a 
\in A_n \text{ such that} \\
  &\,\text{for some } (f',g') \in
\{(f,g),(f^{-1},g^{-1})\} \\
  &\,\text{one of the following possibilities holds} \\
  &\,(a) \quad \text{ for some } m < \omega \text{ we have } f'(a)
\in A_m \Leftrightarrow g'(a) \notin A_m, \\
  &\,(b) \quad f'(a) \ne g'(a) \text{ are in } A_n\}.
\endalign
$$
\mn
4) If $\bar A$ is an $\omega$-representation of $\Bbb A$ and $G = \text{
Aut}(\Bbb A)$ then we define ${\frak d}' = {\frak d}'_{{\Bbb A},\bar A}$ a
metric on $G$ by ${\frak d}'(f,g) = \text{ sup}\{2^{-n}:f \restriction
A_n = g \restriction A_n$ and $f^{-1} \restriction A_n = g^{-1}
\restriction A_n\}$.
\enddefinition
\bigskip

\proclaim{\stag{1.9} Claim}  Assume $\bar A$ is an $\omega$-representation of an
infinite structure $\Bbb A$. 
\nl
1) $(\text{\rm Aut}(\Bbb A),{\frak d}_{{\Bbb A},\bar A})$ is an
 indirectly complete metric group with density $\le 2^{\aleph_0} +
\Sigma\{2^{|A_n|}:n < \omega\}$; 
in fact it is a complete metric space (but in general not a metric
group). 
\nl
2) $(\text{\rm Aut}(\Bbb A),{\frak d}'_{{\Bbb A},\bar A})$ is a complete
metric group of density $\le \Sigma\{\|\Bbb A\|^{|A_n|}:n <
\omega\}$; so if each $A_n$ is finite the density is $\le \aleph_0$. \nl
3) If the universe of $\Bbb A$ is $\omega$ and $A_n = n =
\{0,\dotsc,n-1\}$ \ub{then} $\bar A = \langle A_n:n < \omega \rangle$ is
an $\omega$-representation of $\Bbb A$ with each $A_n$ finite and 
${\frak d}'_{\Bbb A,\bar A} = 
{\frak d}^{\text{aut}}_{\Bbb A}$ from Definition \scite{1.4}(1) and
under it {\rm Aut}$(\Bbb A)$ is a complete separable specially metric group.
\endproclaim
\bigskip

\demo{Proof}  1) Let ${\frak d} = {\frak d}_{\Bbb A,\bar A}$.
First we show that
\mr
\item "{$(*)_1$}"  ${\frak d}$ is a metric (even ultrametric). \nl
[Why?  Clearly, for $f,g \in \text{ Aut}(\Bbb A)$ we have
${\frak d}(f,g)$ is a non-negative real and ${\frak
d}$ is symmetric (i.e. ${\frak d}(f,g) = {\frak d}(g,f))$ and ${\frak d}(f,g) =
0 \Leftrightarrow f=g$.  Mainly we should prove that ${\frak d}(f_1,f_3)
\le {\frak d}(f_1,f_2) + {\frak d}(f_2,f_3)$.  Now if $f_1 = f_2 \vee f_2 =
f_3$ this is obvious, and also if ${\frak d}(f_1,f_2) =1 \vee {\frak
d}(f_1,f_2) = 1$ this is obvious (as ${\frak d}(f_1,f_3) \le 1$)
so assume ${\frak d}(f_\ell,f_{\ell +1}) =
2^{-n(\ell)}$ and $n(\ell) > 0$ for $\ell =1,2$.  
So if $n < n(1)$ and $n < n(2)$  and $m <\omega$ then
we have, for every $a \in A_n$:
{\roster
\widestnumber\item{$(iii)$}
\itemitem{ $(i)$ }  $f_1(a) \in A_m \Leftrightarrow f_2(a) \in A_m$.
\nl
[Why?  As $n < n(1) \text{ and } {\frak d}(f_1,f_2) = 2^{-n(1)}$.]
\sn
\itemitem{ $(ii)$ }   $f_2(a) \in A_m \Leftrightarrow f_3(a) \in A_m$.
\nl
[Why?  As $n < n(2),{\frak d}(f_2,f_3) = 2^{-n(2)}$.]
\nl
hence together
\sn
\itemitem{ $(iii)$ }  $f_1(a) \in A_m \Leftrightarrow f_3(a)
\in A_m$.
\nl
Similarly for $f^{-1}_1,f^{-1}_2,f^{-1}_3$ so
\sn
\itemitem{ $(iv)$ }   $f^{-1}_1(a) \in A_m \Leftrightarrow
f^{-1}_3(a) \in A_m$
\sn
\itemitem{ $(v)$ }  if 
$f_1(a) \in A_n \text{ (equivalently } f_3(a) \in A_n$) then 
$f_1(a) = f_3(a)$. 
\nl
[Why?  First also $f_2(a) 
\in A_n$ by (i).  Second $f_1(a) = f_2(a)$ (as 
$n < n(1)$ and ${\frak d}(f_1,f_2) = 2^{-n(1)})$ 
and third $f_2(a) = f_3(a)$ (as $n < n(2)$ 
and ${\frak d}(f_2,f_3) = 2^{-n(2)})$ hence together
$f_1(a) = f_3(a)$.]
\sn
\itemitem{ $(vi)$ }  Similarly $f^{-1}_1(a) \in A_n \equiv
f^{-1}_3(a) \in A_n \Rightarrow f^{-1}_1(a) = f^{-1}_3(a)$. 
\endroster}
\ermn
Together (and there is $n$ such that $n<n(1),n<n(2))$ 
this gives ${\frak d}(f_1,f_3) \le 2^{\text{Min}\{n(1),n(2)\}}
= \text{ Max}\{{\frak d}(f_1,f_2),{\frak d}(f_2,f_3)\}$.  So $(*)_1$ holds.]
\mr
\item "{$(*)_2$}"  ${\frak d}(f^{-1},g^{-1}) = {\frak d}(f,g)$ \nl
[Why?  Read the definition of ${\frak d}$ in \scite{1.8}(3).]
\sn
\item "{$(*)_3$}"  $G$ as a metric space under ${\frak d}_G$ has
density $\le \Sigma\{2^{|A_n|}:n < \omega\}$ \nl
[Why?  Just look at the definition using easy cardinal arithmetic.]
\ermn
We may wonder whether (Aut$(\Bbb A),{\frak d})$ is complete, i.e. whether
every ${\frak d}$-Cauchi sequence $\langle f_n:n < \omega
\rangle$ in $G$, ${\frak d}$-converge to some $f \in G$.
\mn
Before answering we prove 
a weaker substitute (in fact, it is the one we shall use for
proving the indirect completeness below).

We say that $\langle f_n:n < \omega \rangle$ weakly converge to $f$ if
(they are $\in G$, or are just permutations of $\Bbb A$ and) for every
$\alpha \in \Bbb A$ the sequence $\langle f_n(a):n < \omega \rangle$ is
eventually constant and moreover is eventually equal to $f(a)$. 
\nl
Now we note that:
\mr
\item "{$(*)_4$}"  if $\langle f_n:n < \omega \rangle$ is a ${\frak
d}$-Cauchi sequence, \ub{then} it weakly converge, i.e. for every
$a \in \Bbb A,\langle f_n(a):n < \omega \rangle$ is eventually
constant and $\langle f^{-1}_n(a):n < \omega \rangle$ is
eventually constant so the limit $f$ is a well defined permutation of
$\Bbb A$, moreover belongs to $G = \text{ Aut}(\Bbb A)$.
\ermn
[Why?  Let $x \in \Bbb A$ so for some $n(1) < \omega$ we have $x \in
A_{n(1)}$.  As $\langle f_n:n < \omega \rangle$ is a ${\frak
d}$-Cauchi sequence for some $n(2) < \omega$ 
we have $n \ge n(2) \Rightarrow {\frak d}(f_n,f_{n(2)}) <
2^{-n(1)}$ hence $n \ge n(2) \Rightarrow \dsize \wedge_m (f_n(x) \in
A_m \equiv f_{n(2)}(x) \in A_m)$ by the definition of ${\frak d}$.  
Let $m(1)$ be such that $f_{n(2)}(x) \in
A_{m(1)}$, and let $n(3) \ge n(2)$ be such that $n \ge n(3)
\Rightarrow {\frak d}(f_n,f_{n(3)}) < 2^{-m(1)}$ so necessarily, again
by the definition of the metric, the sequence
$\langle(f_n(x):n = n(2),n(2)+1,\ldots \rangle$ is constant; 
and call this value
$f(x)$.  Similarly concerning $f^{-1}_n(x)$.]
\nl
Note that the above argument shows that 
``if $\bar f$ do ${\frak d}$-converge to $f$ then $\bar f$
weakly converges to some $f$", also it is interesting to note 
(though not explicitly used)
\mr
\item "{$(*)_5$}"  assume $w(x_1,\dotsc,x_n)$ is a group term and
$f^\ell_n \in G$ for $\ell \in \{1,\dotsc,n\},k < \omega$ and
$\langle f^\ell_k:k < \omega \rangle$ does ${\frak d}\text{-converge to }
g_\ell \text{ (or just weakly converge to } g_\ell)$ for 
$\ell \in \{1,\dotsc,k\}$ and we let
$f_k = w(f^\ell_1,\dotsc,f^\ell_n) \text{ for } k < \omega
\text{ \ub{then} } \langle f_k:k < \omega \rangle \text{ weakly
converge to } g =: w(g_1,\dotsc,g_n)$.
\ermn
Very nice, \ub{but} multiplication \footnote{E.g. let $\Bbb A$ be a
trivial structure (i.e. with the empty vocabulary so $G$ is the group
of permutations of $\Bbb A$) and $b_n \ne c_n \in A_{n+1} \backslash
A_n$ and $a_n \in A_0$ be pairwise distinct (for $n < \omega$).  Let
$h_k$ exchange $a_n,c_n$ if $n < k$ and is the identity otherwise.
Let $f$ interchange $a_n,b_n$ if $n < \omega$ and is the identity
otherwise.  Let $g_k$ interchange $b_n,c_n$ if $n < k$ and is the identity
otherwise.  Now $\langle g_k:k < \omega\rangle$ is a Cauchi sequence
and so has a limit $g$, but $h_k = f g_k f$ and $\langle h_k:k <
\omega\rangle$ is not a Cauchi sequence as ${\frak d}(h_{k_1},h_{k_1})
= 1$ if $k_1 < k_2$ as witnessed by $a = a_{k_1}$.}
 is not continuous in general for
this metric, ${\frak d}$.  We also have (though not actually used)
\mr
\item "{$(*)_6$}"  $(\text{Aut}(\Bbb A),{\frak d})$ is a complete
metric space. \nl
[Why?  Let $\bar f = \langle f_n:n < \omega \rangle$ be a ${\frak
d}$-Cauchi sequence, by $(*)_4$ there is $f \in \text{ Aut}(\Bbb A)$
to which $\bar f$ weakly converges.  Now let $n(*) < \omega$ be given
so for some $n(1) < \omega$ we have: $n \ge n(1) \Rightarrow {\frak
d}(f_n,f_{n(1)}) < 2^{-n(*)}$ and we shall prove $n \ge n(1) = {\frak
d}(f,f_n) < 2^{-n(*)}$.  So for each $c \in A_{n(*)}$ we shall check
clauses (a),(b) in the definition of ${\frak d}$ in Definition
\scite{1.8}(3).  First for every $m < \omega$ we have $n \ge n(1)
\Rightarrow f_n(c) \in A_m \equiv f_{n(1)}(c) \in A_m$, but $\langle
f_n(c):n < \omega \rangle$ is eventually constantly $f(c)$ hence
$f(c) \in A_m \equiv f_{n(1)} \in A_m$.  This takes care of clause
(a). \nl
As for clause (b), similarly $n \ge n(1) \wedge
[\{f_n(a),f_{n(1)}(a)\} \subseteq A_n \Rightarrow f_n(a) = f_{n(1)}(a)
\rangle$ hence $n \ge n(1) \wedge [\{f(a),f_{n(1)}(a)\} \subseteq
A_n \Rightarrow f(a) = f_{n(1)}(a) \rangle$.  The same holds for
$\langle f^{-1}_n:n < \omega \rangle,f^{-1}$ so we are done.]
\ermn
But we have to prove that (Aut$(\Bbb A),{\frak d})$ is an indirectly complete
 semi-metric group.  The only clause left is (c) of Definition
\scite{1.7}.  So assume $g_n \in G,{\frak d}(g_n,g_{n+1}) < 1/2^n$,
hence by $(*)_4 + (*)_6$ the sequence 
$\langle g_n:n < \omega \rangle$ converge to some $g \in G$ by
the metric and also weakly converge to $g$.  Let $f_n =
g^{-1}_{2n} g_{2n+1}$, easily $\langle f_n:n < \omega \rangle$ weakly
converge to $e_G = \text{ id}_{\Bbb A}$, let $f = e_G$; so it suffices
to find a metric ${\frak d}'$ such that $(G,{\frak d}')$ is a complete
specially metric group in which $\langle f_n:n < \omega \rangle$
converge to $f$.  We prove this assuming just
\mr
\item "{$\boxtimes$}"   $\bar f = \langle f_n:n < \omega \rangle $ is an
$\omega$-sequence of members of $G$ which weakly converge.
\ermn
Let $B_n = B^n_{\bar f} = 
\{a \in \Bbb A:a \in A_n$ and for every $m \in (n,\omega)$ we
have $f_m(a) = f_n(a) \and f^{-1}_m(a) =
f^{-1}_n(a)\}$.   Clearly $\bar B = \langle
B^n_{\bar f}:n < \omega \rangle$ is an
increasing $\omega$-sequence of subsets of
$\Bbb A$ with union (the universe of) $\Bbb A$.
Recall that ${\frak d}' =: {\frak d}'_{{\Bbb A},\bar B}$ was defined by

$$
{\frak d}'(f,g) = \text{ inf}\{2^{-n}:f \restriction B_n \ne g
\restriction B_n \text{ or } f^{-1} \restriction B_n \ne g^{-1}
\restriction B_n\}.
$$
\mn
Now by parts (2),(3) 
\mr
\item "{$(*)_6$}"  the group $G$ with ${\frak d}'$ is a complete
metric group 
\ermn
and obviously
\mr
\item "{$(*)_7$}"  ${\frak d}(f_n,f_m) < \text{ Max}\{2^{-n},2^{-m}\}$
and ${\frak d}(f_n,f) \le 2^{-n}$ 
\ermn
hence
\mr
\item "{$(*)_8$}"  $\langle f_n:n < \omega \rangle$ converge to $f$ by
${\frak d}'$.
\ermn
So by $(*)_6 + (*)_8$ we are done. \nl
(E.g. Why ``specially" (in part (1) which we are proving)?  As for each $n <
\omega,G_n = \{f \in G:{\frak d}'(f,e_G) < 2^{-n}\}$ is a subgroup of
$G$; the other requirements also are just like the proof of
\scite{1.5}.)
\nl
2),3) Left to the reader.  \hfill$\square_{\scite{1.9}}$
\enddemo
\bn
\margintag{1.11A}\ub{\stag{1.11A} Discussion}  We may consider cases like the
endomorphism semi group of a structure and the endomorphism ring of
an abelian group. The beautiful terms below are as in \cite{Sh:61}.
\nl
We may consider deriving from $\Bbb A$ other structures (in addition
to the automorphism group, endomorphism and monomorphism semi-group);
see also \scite{1.32}.
\bigskip

\definition{\stag{1.12} Definition}   Let $\Bbb A$ be a
structure. \nl
1) A term $\sigma(x_1,\dotsc,x_n)$ in the vocabulary of $M$ 
is called $\Bbb A$-beautiful \ub{if} for
every function symbol $F$ of $\tau_{\Bbb A}$, say $m$-place, the
equations

$$
\align
\sigma(&F(x^1_1,x^1_1,\dotsc,x^1_m),F(x^2_1,x^2_2,\dotsc,x^2_m),\dotsc,
F(x^n_1,x^n_2,\dotsc,x^n_m)) \\
 &=
F(\sigma(x^1_1,x^2_1,\dotsc,x^n_1),\sigma(x^1_2,x^2_2,\dotsc,x^n_2),\dotsc,
\sigma(x^1_m,x^1_m,\dotsc,x^1_m))
\endalign
$$
\mn
are satisfied by $\Bbb A$. 
\nl
3) If $\bar A$ is an $\omega$-representation of $\Bbb A$ \ub{then} we let
${\frak d}^{\text{end}}_{\Bbb A,\bar A},
{\frak d}^{\text{end}'}_{\Bbb A,\bar A}$,
both two-place functions from End$({\Bbb A})$ to $\Bbb R^{\ge 0}$ be
defined as follows: 
\nl
${\frak d}^{\text{end}}_{\Bbb A,\bar A}(f,g) = \text{ sup}\{2^{-n}:(\forall x
\in A_n)(f(x) \in A_n \vee g(x) \in A_n \rightarrow f(x) = g(x))$
and $(\forall x \in A_n)(\forall m < \omega)[f(x) \in
A_m \equiv g(x) \in A_m]\}$, \nl
${\frak d}^{\text{end}'}_{\Bbb A,\bar A}(f,g) =
\text{ sup}\{2^{-n}:f \restriction A_n = g \restriction A_n\}$. \nl
4) End$^+_{\Bbb A} = \text{ End}^+(\Bbb A)$ is the structure with: \nl
\ub{universe} End$(\Bbb A)$, \nl
\ub{functions}: \nl
composition (o), a two-place function 
\nl
$e$, the identity of $\Bbb A$, an individual constant 
\nl
$F^{\text{end}}_{\Bbb A,\sigma}$, for
every $\Bbb A$-beautiful term $\sigma(x_1,\dotsc,x_n)$, an $n$-place
function, where

$$
f = F^{\text{end}}_{A,\sigma}(f_1,\dotsc,f_n)
$$
\mn
is defined by

$$
f(x) = \sigma(f_1(x),\dotsc,f_n(x)).
$$
\mn
5) Aut$^+_{\Bbb A} = \text{ Aut}^+({\Bbb A})$ (or Mono$^+_{\Bbb A} =
\text{ Mono}^+({\Bbb A}))$ is defined similarly restricting ourselves
to beautiful terms $\sigma(x_1,\dotsc,x_n)$ which maps Aut$_{\Bbb A}$
to Aut$_{\Bbb A}$ (or Mono$_{\Bbb A}$ to Mono$_{\Bbb A}$).
\enddefinition
\bigskip

\proclaim{\stag{1.12A} Claim}  For any structure ${\frak A}$: \nl
1) For any ${\frak A}$-beautiful term $\sigma(x_1,\dotsc,x_n)$ the
function $F^{\text{end}}_{{\frak A},\sigma}$ is a full function from
{\rm End}$({\Bbb A})$ into {\rm End}$({\Bbb A})$. 
\nl
2) {\rm End}$^+_{\frak A}$ is a full structure (i.e. the functions are full
not (strictly) partial).
\endproclaim
\bigskip

\proclaim{\stag{1.13} Claim}  Assume
\mr
\item "{$(a)$}"  $\Bbb A$ is a structure
\sn
\item "{$(b)$}"  $\bar A$ is an $\omega$-representation of $\Bbb A$
\sn
\item "{$(c)$}"  ${\frak d} = {\frak d}^{\text{end}}_{\Bbb A,\bar A}$ and
${\frak d}' = {\frak d}^{\text{end}'}_{{\Bbb A},\bar A}$. 
\ermn
1) ({\rm End}$_{\Bbb A},{\frak d})$ is an indirectly complete metric
algebra (which is a semi group). 
\nl
2) ({\rm Mono}$_{\Bbb A},{\frak d} \restriction \text{ Mono}_{\Bbb
A})$ is an indirectly complete metric algebra which is a semi group. \nl
3) ({\rm End}$_{\Bbb A},{\frak d}')$ is an indirectly complete metric
algebra. 
\nl
4) $(\text{\rm Mono}_{\Bbb A},{\frak d}' \restriction 
\,\text{\rm Mono}_{{\frak A},\bar A})$ is an 
indirectly completed metric algebra. 
\nl
5) In parts (3) + (4) the density is $\le \dsize \sum_{n < \omega}
\|\Bbb A\|^{|A_n|}$ and in parts (1) + (2) the density is $\le \dsize
\sum_{n < \omega} 2^{|A_n|} + 2^{\aleph_0}$. 
\nl
6) Assume that for every beautiful $F \in \tau_{\frak A}$, for every $n$ large
enough the set $A_n$ is closed under $F$.  \ub{Then}
$(\text{\rm End}^+_{\frak A},{\frak d})$ is an indireclty complete metric
algebra and $(\text{\rm End}^+_{\Bbb A},{\frak d}')$ is an indirectly  
complete metric algebra.
\endproclaim
\bigskip

\demo{Proof}  Like \scite{1.9}.
\enddemo
\bn
We may wonder: \nl
\margintag{1.20}\ub{\stag{1.20} Question}:  Can we have an uncountable Polish algebra
$G$ (so $\tau_G$ is finite or just countable) which is free for some
variety?

Of course, if $\tau_G$ is empty this holds; obviously we may discard many.
\bigskip

\demo{\stag{1.21} Example}  The answer to \scite{1.20} is yes.
\enddemo
\bigskip

\demo{Proof}  Note that if $G$ is the
vector space over a countable field $F$ with basis $\langle x_\eta:\eta \in
{}^\omega 2 \rangle$, it is a metric space with countable
density, i.e. is Polish, as we can define
${\frak d}(x,y) = \|x-y\|$ where $\|0\| = 0$ and for $x =
\sum\{q_{\eta,x} x_\eta:\eta \in {}^\omega 2\} \ne 0$ (so $\{\eta \in
{}^\omega 2:q_{\eta,x} \ne 0\}$ is finite and, of course, $\eta \in
{}^\omega 2 \Rightarrow q_{\eta,x} \in F)$ we let $\|x\|$ be $2^{-n}$
where $n$ is minimal such that for some $\nu \in {}^n2$ we have
$0 \ne \Sigma\{q_{n(\eta)}:\nu \triangleleft \eta \in {}^\omega 2$ and
$q_{n(\eta)} \ne 0\}$.

Clearly $G$ is a metric space; it has density $\aleph_0 +
|F|$, and the addition and substruction are continuous, and it is separable
for $F$ countable.  So, for countable $F$, the 
completion $\hat G$ is an (additive) Polish group
and vector space.  Moreover, if $F = \Bbb Z/p
\Bbb Z$ for $p$ a prime, $\hat G$ is the free group of the appropriate
variety.  \hfill$\square_{\scite{1.20}}$
\enddemo
\bigskip

\remark{\stag{1.22} Remark}  Another metric on the same space is: 
for $x \in G$ let $n(x) = \text{ Min}\{n$: if
$q_{\eta,x} \ne 0 \wedge q_{\nu,x} \ne 0$ and $\eta \ne \nu$ then $n
\ge \ell g(\eta \cap \nu)\}$ and let ${\frak d}(x,y)$ be:
\endremark
\bn 
\ub{Case 1}:  $n(x) \ne n(y)$ or $(\exists \eta,\nu \in {}^\omega
2)[\eta \restriction n(x) = \nu \restriction n(y) \and q_{\eta,x} \ne
0 \and q_{\nu,y} \ne 0 \and q_{\eta,x} \ne q_{\eta,y}]$. \nl
Then ${\frak d}(x,y) = 2$.
\bn
\ub{Case 2}:  Otherwise

$$
\align
{\frak d}(x,y) = \text{ sup}\{2^{-n}:&n \ge n(x) \text{ and} \\
  &(\forall \eta,\nu \in {}^\omega 2)[\eta \restriction n(x) = \nu
\restriction n(y) \and q_{\eta,x} \ne 0 \and q_{\nu,y} \ne 0
\rightarrow \\
  &\ell g(\eta \cap \nu) \ge n]\}.
\endalign
$$
\mn
Now $G$ under ${\frak d}$ is a complete metric space, but it is not a
metric group.
\bn
\centerline {$* \qquad * \qquad *$}
\bn
Closely related to semi metric (see Definition \scite{1.7}, but not
enough for our theorems) is:
\definition{\stag{1.30} Definition}  1) We say ${\frak a} = (M,{\frak d},\bold U)$ is a
metric-topological algebra \ub{if}:
\mr
\item "{$(a)$}"  $M$ is an algebra
\sn 
\item "{$(b)$}"  ${\frak d}$ is a metric (on the universe of $M$)
\sn
\item "{$(c)$}"  $\bold U$ is a Hausdorff topology (i.e. the family of open
sets) on the universe of $M$ \nl
such that
\sn
\item "{$(d)$}"  the operations of $M$ are continuous by
${\frak d}$ and by $\bold U$
\sn
\item "{$(e)$}"  every open ${\frak d}$-ball i.e. set of the form
$\{a:{\frak d}(a,a_0) < \zeta\}$, is open also by the
topology $\bold U$.
\ermn
2) We say ${\frak a} = (M,{\frak d},\bold U)$ is complete \ub{if} every
${\frak d}$-Cauchi sequence converge to some point of $M$ by the
topology $\bold U$ though not necessarily by ${\frak d}$.
\enddefinition

\proclaim{\stag{1.31} Claim}  Assume
\mr
\item "{$(a)$}"  $\Bbb A$ is a structure
\sn
\item "{$(b)$}"   $\bar A$ is an $\omega$-representation of $\Bbb A$
\sn
\item "{$(c)$}"  $G = Aut(\Bbb A)$
\sn
\item "{$(d)$}"   ${\frak d} = {\frak d}_{\Bbb A,\bar A}$, see
Definition \scite{1.8}
\sn
\item "{$(e)$}"  $\bold U$ is the topology on $G$ such that a
neighborhood basis for $f \in G$ is $\{U_{f,X}:X \subseteq \Bbb A$
finite$\}$ where $u_{f,X} = \{g \in G:f \restriction X = g
\restriction x$ and $f^{-1} \restriction X = g^{-1} \restriction X\}$.
\ermn
\ub{Then} $(G,{\frak d},\bold U)$ is a complete metric-topological
algebra.
\endproclaim
\bigskip

\demo{Proof}  Included in the proof of Claim \scite{1.9}.
\hfill$\square_{\scite{1.31}}$
\enddemo
\bn
\margintag{1.32}\ub{\stag{1.32} Discussion}:  1) We can generally use topology
instead of metric.  What is the gain? \nl
2) Instead of automorphisms we can consider a universal Horn Theory
$T$ in a vocabulary $\tau^+ = \tau^*_T$ extending $\tau_{\Bbb A}$, 
e.g. $\tau^+ = \tau_{\frak A} \cup\{F^*\},F^*$ a function symbol with
arity $n^*$.  So

$$
\text{Exp}_T(\Bbb A) =: \{\Bbb A^+:\Bbb A^+ \text{ a } \tau^+
\text{-expansion of } {\frak A} \text{ and is a model of } T\},
$$ 
\mn
if e.g. $\tau^+ \backslash \tau_{\frak A} = \{F\}$, we may replace
Exp$_T(A)$ by $\{F^{{\Bbb A}^+}:\Bbb A^+ \in \text{ Exp}_T(\Bbb A)\}$.  \nl
We may define

$$
\align
\text{beautiful}(T,\Bbb A) = \bigl\{ \bar \sigma(\bar x):&\sigma
\text{ a } \tau_{\Bbb A} \text{-term}, \bar x = (x_1,\dotsc,x_k) \\
  &\text{and if } \Bbb A^+_\ell \in \text{ Exp}_T(\Bbb A) \text{ for }
\ell = 1,\dotsc,k \text{ and} \\
  &F \in \tau^+ \backslash \tau_{\Bbb A} \text{ we define the
arity}(F) \text{-place function } F^*_\sigma \\
  &\text{from } \Bbb A \text{ to } \Bbb A \text{ by } 
F^*_\sigma(\bar a) = \sigma(F^{{\Bbb A}_1}(\bar a),\dotsc,
F^{{\Bbb A}_k}(\bar a)) \\
  &\text{and let } \Bbb A^+ = (\Bbb A,F^*_\sigma)_{F \in \tau^+ \backslash
\tau_{\Bbb A}} \text{ then } \Bbb A^+ \text{ is a model of } T
\bigr\},
\endalign
$$
\mn 
we can consider
more complicated operations. So $(\text{Exp}_T(\Bbb A),H_{\bar
\sigma(\bar x)})_{\bar \sigma(\bar x) \in \text{ beautiful}(\Bbb A,T)}$ is a
generalization of Aut$(\Bbb A)$ where $H_{\bar \sigma(\bar x)}(\Bbb
A^+_1,\dotsc,\Bbb A^+_k) = \Bbb A^+$ is defined as above.
\newpage

\head {\S3 Compactness of metric algebras} \endhead  \resetall \sectno=3
 \spuriousreset
\bn
Note that below if $u_n = \{t_n\} = \{n\}$ we may write $x_n$ instead
of $\bar x_n,n$ instead $t \in u_n$ and $d_n$ instead of $\bar d_n$.
\proclaim{\stag{s.1} The completeness Lemma}  Assume 
${\frak a}$ is a Polish algebra $M
= M_{\frak a}$ (so with countable vocabulary) such that
\mr
\item "{$\circledast(a)$}"  $\langle u_n:n < \omega \rangle$ is a
sequence of pairwise disjoint non-empty finite sets 
\sn
\item "{${{}}(b)$}"  $\bar x_n = \langle x_t:t \in u_n \rangle$
\sn
\item "{${{}}(c)$}"  $\bar \sigma_n(\bar x_{n+1}) = \langle
\sigma_{n,t}(\bar x_{n+1}):t \in u_n \rangle$ is a sequence of
$\tau_M$-terms, possibly with parameters (from $M_{\frak a}$) so
$\bar \sigma_n(\bar d) = \langle \sigma_{n,t}(\bar d):t \in u_n
\rangle$ for any $\bar d = \langle d_s:s \in u_{n+1} \rangle,d_s \in
M$; if ${\frak a}$ is a Polish group, the $\sigma_n$ are so called words
\sn
\item "{${{}}(d)$}"  $\zeta = \langle \zeta_n:n < \omega \rangle$ is a
sequence of positive reals converging to 0
\sn
\item "{${{}}(e)$}"  $\bar d_{n+1} = \langle d_{n+1,t}:t \in u_{n+1}
\rangle$ with each $d_{n+1,t}$ an element of $M$ such that if $\bar
d'_{n+1} = \langle d'_{n+1,t}:t \in u_{n+1} \rangle$ is of distance $<
\zeta_{n+1}$ from $\bar d_{n+1}$, (that is $d'_{n+1,t} \in 
\,\text{\rm Ball}_G(d_{n+1,t},\zeta_{n+1})$ 
for each $t \in u_{n+1}$), \ub{then}
$\bar \sigma_n(\bar d'_{n+1}) \in \,\text{\rm Ball}(\bar d_n,\zeta_n)$
which means: $t \in u_n \Rightarrow \bar \sigma_{n,t}
(\ldots,d'_{n+1,s},\ldots)_{s \in u_{n+1}} \in \,\text{\rm Ball}
(d_{n,t},\zeta_n)$
\sn
\item "{${{}}(f)$}"  for every $n < \omega$ and a position real
$\varepsilon$ there is $m>n$ such that
{\roster
\itemitem{ $(*)_\varepsilon$ }  if $d'_{m,t} \in \,\text{\rm Ball}
(d_{m,t},\zeta_n)$ for every 
$t \in u_m$ \ub{then} the distance between \nl
$\bar \sigma_n(\bar \sigma_{n-1}(\ldots,\bar \sigma_{m-1}(\bar
d'_m),\ldots),\bar \sigma_n(\bar \sigma_{n+1}
(\ldots,\bar \sigma_{m-1}(\bar d_m))\ldots)$ is $< \varepsilon$.
\endroster}
\ermn
\ub{Then} there are $d^*_{n,t} \in M$ for $n < \omega,t \in u_n$ which solves
the set of equations

$$
d^*_{n,t} = \sigma_{n+1}(d^*_{n+1,s})_{s \in u_{n+1}}
$$
and satisfies

$$
d^*_{n,t} \in \,\text{\rm Ball}_G(d_{n+1},\zeta_n).
$$
\endproclaim
\bigskip

\remark{\stag{s.2} Remark}  1) In ``special" versions we have $\bar d_n
= \bar \sigma_n(\bar d_{n+1})$ (and in \cite{Sh:744} we have $d_n =
\sigma_n(d_{n+1}))$ but here there is no ``the true solution which we
perturb".
\nl
2) Condition (e) in Lemma \scite{s.1} says 
that if in large $n$ we perturb $\bar d_n$ with
error $< \zeta_n$ and compute down by the $\bar \sigma$'s we still get
a reasonable $\bar d_k$ for every $k < n$ but not necessarily a very
good one.
\endremark
\bigskip

\demo{Proof}  For every $k$ we shall define $\langle c^k_{n,t}:
t \in u_n,n < \omega \rangle$, a sequence of elements of the algebra.

First, if $n \ge k$ let $c^k_{n,t} = d_{n,t}$.  Second, we define
$\bar c^k_n = \langle c^k_{n,t}:t \in u_n \rangle$ by downward
induction on $n \le k$.
\sn
\ub{$n=k$} by the first case.
\sn
\ub{$n<k$} let $\bar c^k_n = \bar \sigma_n(\bar c^k_{n+1})$.
\sn
Next we show
\mr
\item "{$(*)_1$}"  $c^k_{n,t} \in \text{ Ball}(d_{n,t},\zeta_n)$. \nl
[Why?  If $n \ge k$ this is trivial as $c^k_{n,t} = d_{n,t}$.  If $n \le k$ by
downward induction on $n$, using condition (e) of the assumptions.]
\sn
\item "{$(*)_2$}"  for every positive real $\varepsilon > 0$ and $n <
\omega$, there is $m > n$ such that if $k \ge m$ then $t \in u_n
\Rightarrow c^k_{n,t} \in \text{ Ball}(c^m_{n,t},\varepsilon)$. \nl
[Why?  Given $n < \omega$ and $\varepsilon > 0$ choose $m$ as in clause (f)
of the assumption.  Let $k \ge m$.  By $(*)_1$, $t \in  u_m \Rightarrow
c^k_{m,t} \in \text{ Ball}(d^k_{n,t},\zeta_m)$.  By the way the
$c^k_{i,t}$ were defined for $t \in u_n$ we have $\bar c^k_n = \bar
\sigma_n(\bar \sigma_{n+1}(\ldots,\bar \sigma_{n-1}(\bar
c^k_m)\ldots))$ and similarly $\bar c^m_n = \bar \sigma_n(\bar \sigma_{n+1}
\ldots \bar \sigma_{m-1}(\bar d_m)\ldots)$. 
\nl
Condition (f) i.e. the choice of $m$ tells us that the desired results holds.]
\ermn
So, for each $n$ and $t \in u_n$ the sequence 
$\langle c^k_{n,t}:k < \omega \rangle$ is a Cauchy sequence by
$(*)_2$; 
hence it converges to some $c_{n,t} \in M$.  Now
\mr
\item "{$\boxtimes$}"  the sequence $\langle c_{n,t}:n < \omega,t \in u_n
\rangle$ forms a solution: for every $n < \omega$ and $t \in u_n$ the
equation $c^k_{n,t} = \sigma_{n,t}(\ldots,c^k_{n+1,s})$ is satisfied
whenever $n > k$ hence in the limit 
$c_{n,t} = \sigma_{n,t}(\ldots,c_{n+1,s},\ldots)_{s \in u_{n+1}}$.
\hfill$\square_{\scite{s.1}}$ 
\endroster
\enddemo
\bn
Recall about groups \nl
\margintag{s.3}\ub{\stag{s.3} Fact}:  A free group is torsion free and the group is
not divisible, in fact, every element $c$ has at most one $n$-th root for
each $n = 1,2,\ldots$ and has no root for every large enough $n$
except when $c$ is the unit.
\bn
\margintag{s.4}\ub{\stag{s.4} Fact}:  Every countable subgroup at a free group $G$
is contained in a countable subgroup which is a retract of $G$.
\bn
We now give a criterion to show non-freeness.  We could use $\bar x$
instead of $x$, of course.
\proclaim{\stag{s.5} Claim}  1)
\mr
\item "{$(a)$}"  ${\frak a}$ is a complete metric algebra, $M = M_a$
with unit $e_M$
\sn
\item "{$(b)$}"  $B \subseteq M$ is countable with $e_M$ belonging to
the closure of $B \backslash \{e_u\}$
\sn
\item "{$(c)$}"  $\Xi$ is a set of terms of the form $\sigma(x,\bar
y)$
\sn
\item "{$(d)$}"  if $\sigma(x,\bar y) \in \Xi$ and $\bar b \in B$ and
$c \in G$ then $\{x \in M:c = \sigma(x,\bar b)\}$ is finite (or at
most a singleton)
\sn
\item "{$(e)$}"  for every finite $A \subseteq M$ (or $A \subseteq M$
a singleton) and $\zeta$ a positive real there are a sequence $\bar b$
from $B$ and term $\sigma(x,\bar y) \in \Xi$ such that
{\roster
\itemitem{ $(\alpha)$ }  $\sigma(e_M,\bar b) \in \,\text{\rm Ball}(e_M,\zeta)$
\sn
\itemitem{ $(\beta)$ }  $\sigma(c,\bar b) \notin A$ for every $c \in M$.
\endroster}
\ermn
\ub{Then} no countable subalgebra of $M$ containing $B$ is a retract
(in the algebraic sense) of $M$.  Hence $M$ is not free for any variety.
\nl 
2) We can omit $e_{\frak a}$, i.e. omit clause (b) and the last phrase
of clause (a) and change clause (e) to
\mr
\item "{$(e)'$}"  for any finite $A \subseteq M_{\frak a}$ [or $A
\subseteq M_{\frak a}$ a singleton] and real $\zeta > 0$ and $d \in M_{\frak
a}$ \ub{there is} a term $\sigma(x,\bar y) \in \Xi$ and sequence $\bar
b \in {}^{\ell g(\bar y)} B$ and element $d' \in M_{\frak a}$ such
that
{\roster
\itemitem{ $(\alpha)$ }  $\sigma^{M_{\frak a}}(d',\bar b) \in 
\,\text{\rm Ball}_{\frak a}(d,\zeta)$
\sn
\itemitem{ $(\beta)$ }  for no $c \in M_{\frak a}$ do we have
$\sigma^{M_{\frak a}}(c,\bar b) \in A$.
\endroster}
\endroster
\endproclaim
\bigskip

\remark{\stag{s.5A} Remark}  We can similarly phrase sufficient
conditions for ``$M_{\frak a}$ is unstable in $\aleph_0$"
[for quantifier free formulas, see \S6].
\endremark
\bigskip

\demo{Proof}  1) Like the proof of part (2) below except that we add
to $(*)$:
\mr
\item "{$(\eta)$}"  $e_n = e_M$.
\ermn
2) We rely on \scite{s.1}.

Assume toward contradiction that $M$ is a countable reduct of $M_{\frak a}$
which includes $B$, so we can choose
$h^*$, a homomorphism
from $M_{\frak a}$ onto $M$ which extends id$_M$.  Let $\langle a_n:n
< \omega \rangle$ list $M$.  Let $u_n = \{n\}$.  We choose $\bar b_n$
and $\sigma_n(x,\bar y_n)$ and $\zeta_n$ by induction on $n$ such that
\mr
\item "{$(*)(\alpha)$}"  $\sigma_n(x,\bar y_n) \in \Xi$
\sn
\item "{$(\beta)$}"  $\bar b_n$ a sequence from $M$ of length 
$\ell g(\bar y_n)$
\sn
\item "{$(\gamma)$}"  $\zeta_n$ a positive real, $\zeta_{n+1} <
\zeta_n/2$ and $\zeta_{n+1,\ell}$ is a positive real $< \zeta_n$
\sn
\item "{$(\delta)$}"  $e_n \in M_{\frak a}$  
\sn
\item "{$(\varepsilon)$}"  $\sigma_n(e_{n+1},\bar b_n) \in \text{
Ball}_G(e_n,\zeta_n)$, moreover if $c' \in \text{
Ball}_G(e_{n+1},\zeta_{n+1})$ \ub{then} $\sigma_n(c,\bar b_n) \in \text{ Ball}
(e_n,\zeta_n)$  
\sn
\item "{$(\zeta)$}"  if $k < n$ and $c',c'' \in \text{
Ball}_G(e_{n+1},\zeta_{n+1})$ and we define the terms
$\sigma_{n+1,\ell}(x)$ for $\ell \le n+1$, with parameters by downward
induction on $\ell$ as follows
$\sigma_{n+1,n+1}(x) = x,\sigma_{n+1,\ell}(x) =
\sigma_\ell(\sigma_{n+1,\ell +1}(x),\bar b_\ell)$ 
\ub{then} the ${\frak d}$-distance between
$\sigma_{n+1,k}(c')$ and $\sigma_{n+1,k}(c'')$ is $< \zeta_n$.
\ermn
Let us carry the induction, in stage $n$ we choose $e_n,\zeta_n$ and
$\sigma_{n-1}(x,\bar b_{n-1})$ if $n > 0$.
\enddemo
\bn
\ub{Case 1}:  $n=0$.

This is straightforward.
\bn
\ub{Case 2}:  $n=k+1$.

Let $D$ be the set of $\bar c = \langle c_m:m \le k \rangle$ 
which satisfies
\mr
\item "{$\boxtimes_k(i)$}"  $m < k \Rightarrow c_m =
\sigma_m(c_{m+1},h^*(\bar b_m))$
\sn
\item "{$(ii)$}"   $c_0 = a_k$.  
\ermn
We can prove by induction on $m \le k$ that $\{c_m:\bar c \in D\}$
is finite, and let $A = \{c_k:\bar c \in D\}$.  By clause $(e)'$ of the
assumptions (see \scite{s.5}(2)), there are $ r < \omega,\sigma = \sigma(x,y_0,
\dotsc,y_{r-1}) \in \Xi$ and $\bar b \in {}^r(M_{\frak a})$
and $d'$ as there.  We let $\bar d_k = \bar b,\sigma_k = \sigma,\bar
y_k = \langle y_0,\dotsc,y_{r-1} \rangle,e_n = d'$.

Lastly, we should choose $\zeta_n \in \Bbb R^+$.  There are several
demands but each holds for every small enough $\zeta > 0$, more exactly
one for clause $(\varepsilon)$ and for each $m<n$, one for clause $(\zeta)$.

Having carried the induction, clearly \scite{s.1} apply
hence there is a solution $\langle d^*_n:n < \omega \rangle$,
that is $M_{\frak a} \models d^*_m = \sigma_m(d^*_{m+1},\bar b_m)$ for
$m < \omega$.  
But $h^*$ is a homomorphism from $M_{\frak a}$ into $M$
so $\langle h^*(d^*_n):n < \omega \rangle$ satisfies all the equations
in $\boxtimes_k$ hence by our choice in stage $n = k+1,h^*(d^*_0) \ne
a_k$. As this holds for every $k$ and $\{a_k:k < \omega\}$ list
the elements of $M$ we are done.  \hfill$\square_{\scite{s.5}}$
\bigskip

\remark{\stag{3.5} Remark}  1) If we phrase algebraic compactness, it is
preserved by taking reducts. \nl
2) In a reasonalbe variant we 
can replace ``$M$ countable" by $\|M\| < \text{ cov}$(meagre);
we'll return to this elsewhere. \nl
3) We can change the demand on $\Xi$: at most one solution in clause
(e), $A$ a singleton in clause (f). \nl
4) This suffices for groups.
\endremark
\newpage

\head {\S4 Conclusions} \endhead  \resetall \sectno=4
 \spuriousreset
\bn
\margintag{c.1}\ub{\stag{c.1} Conclusion}  1) If $(G,{\frak d})$ is a complete metric group of
density $<|G|$, \ub{then}:
\mr
\item "{$(a)$}"  $G$ is not free,
\sn
\item "{$(b)$}"  if $G$ is $\aleph_1$-free then for some
countable $A \subseteq G$, there is no countable reduct $B,M$ of
$M_{\frak a}$ including $A$.
\ermn
2) It suffices that $G$ is an indirectly complete metric group and as
a metric space it is of density $< |G|$. 
\nl
3) Instead ``density" $<|G|$" it is enough to assume that the 
topology induces by the metric is not discrete.
\bigskip

\demo{Proof}  1), 2) Easy.  
Let $\mu$ = density$(G)$ and $y_i \in G$ be pairwise
distinct for $i < \mu^+$.  Without loss of generality $y_i \notin
\langle \{y_j:j < i\} \rangle_G$.  So for some increasing sequence
$\langle i_n:n < \omega \rangle$ the sequence $\langle y_{i,n}:n < \omega
\rangle $ is a Cauchi sequence.  

For part (1), by 
completeness it converges say to $y^*$, the convergence is 
for ${\frak d}_G$.
Now $\langle b_n:n < \omega \rangle =: 
\langle (y^*)^{-1} y_{i_{2n+1}}:n < \omega \rangle$ converges to
$e_G$, the members are pairwise distinct so \wilog \, $\ne e$.

However for part (2) we know that some $(G,{\frak d}')$ equal to
$(G,{\frak d})$ as a group but with a different metric; is a complete
metric group with an $\omega$-sequence of members of $G \backslash
\{e_G\}$ converging to $e_G$.

Let $\Xi = \{x^my_1:m < \omega\}$ and $B = \{b_n:n < \omega\}$.  
Now we shall apply \scite{s.5}. 
In the assumptions, clauses (a)-(c) are obvious.  As for clause (d) we
are using:  equations of the form $x^m a' = a''$ has at most one solution in
$G$, see \scite{s.3}.  
We are left with clause (e), so we are given a real $\zeta > 0$ and a
finite set $A \subseteq G$ (in fact, a singleton is enough).  We can
choose $b \in B \backslash A \backslash \{e_G\}$ of distance $< \zeta$
from $e_G$.  Let $\sigma(x,y) = xy$ and $\bar b = \langle b \rangle$
where $n < \omega$ is the minimal $n > 1$ such that $[a \in A
\Rightarrow ab^{-1}$ has no $n$-th root].  This is possible, see Fact
\scite{s.3} so 
$\sigma(e_M,\bar b) = b \in \text{ Ball}_G(e_M,\zeta)$ as required
in subclause $(\alpha)$ of clause (e) and $\sigma(e_M,\bar b) = b
\notin A$ as required in subclause $(\beta)$ of clause (e) so by
\scite{s.5} we are done.
\nl
3) Choose $\langle y_n:n < \omega \rangle$ converging to some $y^*$
such that $\langle y_n:n < \omega \rangle \char 94 \langle y^*
\rangle$ is with no repetitions, possible on $(G,{\frak d})$ is not
discrete.  Now continue as above.  \hfill$\square_{\scite{c.1}}$ 
\enddemo
\bn
In particular \nl
\margintag{c.2}\ub{\stag{c.2} Conclusion}:  There is no free uncountable Polish group.
\bigskip

\proclaim{\stag{c.3} Claim}  (1) In the proof of Proposition
\scite{c.1}(b) we do not use all the strength of ``$G$ is free".
E.g. if $(G,{\frak d})$ is a complete metric group then $(a)
\Rightarrow \neg(b)$ where: 
\mr
\item "{$(a)$}"  for some group words 
$w_n(x_1,\dotsc,x_{r_n})$ for $n < \omega$ possibly with 
parameters in $G$ we have
{\roster
\itemitem{ $(\alpha)$ }  there are $y_i \in G$ for $i < \mu^+$ (where
$\mu =$ {\rm density}$(G)$) such that $i < j \Rightarrow y_i \ne y_j$
\sn
\itemitem{ $(\beta)$ }  for some $k$, for every $b,a_2,\dotsc,a_{2_n}
\in G$ the set $\{a_1 \in G:G \models w_n(a_1,\dotsc,a_{r_n}) = b\}$
has at most $k$ members
\sn
\itemitem{ $(\gamma)$ }  for every real $\zeta > 0$, finite $A
\subseteq G \backslash \{e_G\}$ and an infinite set $B \subseteq G$
such that $e_G$ belongs to its closure, there are $n < \omega$ and
$b_2,\dotsc,b_{r_n} \in B$ such that $w_r(e_\mu,b_2),\ldots) \in
\,\text{\rm Ball}(e_G,\zeta)$ and $A$ is disjoint to $\{w_r(c,b_r,\ldots):c
\in G\}$
\endroster}
\item "{$(b)$}"  if $X$ is a countable subset of $G$, \ub{then} there is a
countable subgroup $H$ of $G$ which includes $X$ and is a reduct of
$G$, that is there is a projection from $G$ onto $H$.
\ermn
2) The uncountable free abelian group falls under this criterion, in
fact, any uncountable strongly $\aleph_1$-free abelian group 
also satisfies this criterion. 
\nl
3) In part (1) we can weaken (b) to
\mr
\item "{$(b)^-$}"  $G$ is strongly $\aleph_1$-free
\ermn
or just:
\mr
\item "{$(b)^{- -}$}"   for 
every countable $X \subseteq G$ there is a countable
subgroup $H$ of $G,X \subseteq H$ such that: if $H \subseteq H'
\subseteq G,H'$ countable then $H$ is a reduct of $H'$, i.e. there is
a projection from $H'$ onto $H$.
\endroster
\endproclaim
\bigskip

\demo{Proof}  The same as the proof of \scite{c.1}.
\enddemo
\bigskip

\remark{\stag{c.3A} Remark}  1) We 
may consider for a metric space a group rank: the
objects being finite approximation to the system of elements we
actually use $\langle \sigma_n(\bar d_{n+1},\bar b_n):n < \omega
\rangle$ in \scite{s.1} (or in \scite{s.5}).
\nl
2) The results above confirms the thesis that the compactness
conditions say that $G$ is ``large", ``rich".  \nl
3) Note that we can expand $M_{\frak a}$ by
individual constants, equivalently consider terms with parameters. \nl
4) In the applications of \scite{c.3}, we do not actually use $r>2$. \nl
5) Concerning semi groups we intend to say it in a continuation. \nl
6) We may consider assumption ``some $h:M_A \rightarrow N$ is a homomorphisms
onto $N,N$ countable", and rephrase the criterion in \scite{s.5}. \nl
7) We may consider just $\|N\| < 2^{\aleph_0}$, so have to split into
two so we get ${}^\omega 2$ cases among which at least one ``succeeds".
\endremark
\bn
\margintag{c.4}\ub{\stag{c.4} Conclusion}:   1) Assume $\Bbb A$ is a countable
structure.  Then Aut$(\Bbb A)$, the group of automorphisms of $\Bbb
A$, is not a free uncountable group, in fact it satisfies the
conclusions of \scite{c.1}, \scite{c.3}. \nl
2) Assume $\Bbb A$ is a structure of cardinality $\lambda$ and
$\lambda = \mu = \beth_\omega$ or more generally $\lambda =
\Sigma\{\lambda_n:n < \omega\},2^{\lambda_n} < 2^{\lambda_{n+1}},\mu =
\Sigma\{2^{\lambda_n}:n < \omega\} < 2^\lambda$.  \ub{Then} Aut$(\Bbb A)$
cannot be free of cardinality $> \mu$, in fact, it satisfies the
conclusions of \scite{c.1}, \scite{c.3}.
\bigskip

\demo{Proof}  1) By \scite{1.5}, Aut$(\Bbb A)$ is a Polish group and
apply \scite{c.2}. \nl
2) Without loss of generality the universe of $\Bbb A$ is $\lambda$,
using $\bar A = \langle \lambda_n:n < \omega \rangle$ we know by
\scite{1.9}(1) that
$(\text{Aut}(\Bbb A),{\frak d}_{\Bbb A,\bar A})$ is a complete 
semi-metric group and apply \scite{c.1}(2).  \hfill$\square_{\scite{c.4}}$
\enddemo
\bigskip

\proclaim{\stag{c.6} Claim}  For complete specially metric groups the
proof of \cite{Sh:744} works, similarly for algebras.
\endproclaim
\newpage

\head {\S5 Quite free but not free abelian groups} \endhead  \resetall \sectno=5
 \spuriousreset
\bigskip

If uncountable Polish groups are not free, we may look at wider
classes: $F_\sigma$, Borel analytic, projective $\bold L[\Bbb R]$. 
\bn
\margintag{6.0}\ub{\stag{6.0} Question}:   1) Is the freeness of a reasonably
definable abelian group absolute? 
\nl
2) For which cardinals $\lambda$ does $\lambda$-freeness imply freeness (or
$\lambda^+$-freeness) for nicely definable abelian groups, in
particular for $\lambda = \aleph_\omega$?
\nl
3) Similarly for other varieties (or any case when ``free" is definable
like universal Horn theory). 

This is connected also to
\cite{Sh:402} whose original aim was a question of
Marker ``are there non-free Whitehead Borel Abelian groups".  But
already in \cite{Sh:402} it seems to me the basic question is to 
clarify freeness in such groups; that is, question \scite{6.0} above.
\nl
Blass asked about definable subgroups of $\Bbb Z^\omega$ (see question
\scite{6.10}): by \cite{Sh:402} and the construction here we quite
resolve this.

Recall that \cite{Sh:402} analyze
$\aleph_1$-free abelian groups which are $\Sigma^1_1$ or so.  A natural
dividing line was suggested; the complicated half was proved to be not
Whitehead, and at least for me is an analog to not $\aleph_0$-stable.  The low
half is $\aleph_2$-free.  So under CH we were done, but what if
$2^{\aleph_0} > \aleph_1$?  Are they also free?  This was left open by
\cite{Sh:402}.

We shed some light by giving an example (an $F_\sigma$ one) 
showing that the non-CH case in \cite{Sh:402}
is a real problem.  This resolves the original problem: it is
consistent that there are non-free Whitehead groups, this is derived in
\scite{6.16}.  But what about the further question,
e.g. \scite{6.0}(2)?  The examples seem to indicate (at least to me) 
that the picture in
\cite{Sh:87a}, \cite{Sh:87b} is the right one here, connecting theories
of $\psi \in \Bbb L_{\omega_1,\omega}$ with $\Sigma^1_1$-models.  Also
related are \cite{EM2}, \cite{MkSh:366} on almost freeness for
varieties, and see \cite{EM} on abelian groups.  In particular we
conjecture  ``every $\aleph_\omega$-free Borel group is free".

We shall use freely the well known theorem saying
\mr
\item "{$\boxtimes$}"  a subgroup of a free abelian group is a free
abelian group.
\endroster
\bigskip

\definition{\stag{6.1} Definition}  For $k(*) < \omega$ we define
an abelian group $G = G_{k(*)}$ and is generated
by $\{x_{m,\bar \eta,\nu}:m \le k(*) \text{ and } \nu \in
{}^{\omega >} 2 \text{ and } 
\bar \eta = \langle \eta_\ell:\ell \le k(*),\ell \ne m \rangle$ where
$\eta_\ell \in {}^\omega 2\} \cup \{y_{\bar \eta,n}:n < \omega$ and
$\bar \eta = \langle \eta_\ell:\ell \le k(*) \rangle$ where $\eta_\ell \in
{}^\omega 2\}$ freely except the equations:
\mr
\item "{$\boxtimes_{\bar \eta,n}$}"  $(n!)y_{\bar \eta,n+1} = y_{\bar \eta,n}
+$ \nl
$\sum\{x_{m,\bar \eta_m,\nu}:m \le k(*) \text{ and } \bar \eta_m =
\bar \eta \restriction\{m' \le k(*):m' \ne m\}$ and \nl

$\qquad \qquad \quad \nu = \eta_m \restriction n\}$.
\ermn
(Note that if $m_1 < m_2 \le k(*)$ then $\bar\eta_{m_1} \ne \bar
\eta_{m_2}$ having different index sets).
\enddefinition
\bigskip

\demo{Explanation}  A canonical example of a non-free group is
$(\Bbb Q,+)$.  Other examples are related to it after we divide by
something.  The $y$'s here play that role of providing (hidden) copies
of $\Bbb Q$.  What about $x$'s?   For each $\bar \eta \in \Lambda$ we use
$m \le k(*)$ to give $\langle y_{\bar \eta,n}:n < \omega \rangle,k(*)$
``chances", ``opportunities" to avoid having $(\Bbb Q,+)$ as 
a quotient, one for each cardinal
$\le \aleph_{k(*)}$.  More specifically, if $H \subseteq G$ is the
subgroup which is generated by $X = \{x_{m,\bar \eta,\nu}:m \ne
m(*)$ and $\eta$ is a function from $\{\ell \le k(*):\ell \ne m\}$ to
$\omega$ and $\nu \in {}^{\omega >} 2\}$, still in $G/H$ the 
$\{y_{\bar \eta,n}:n < \omega\}$ does not
generate a copy of $\Bbb Q$, as witnessed by 
$\{x_{m(*),\bar \eta_{m(*)},\eta_{m(*)} \restriction n}:n < \omega\}$.
\enddemo
\bigskip

\proclaim{\stag{6.1.7} Claim}  The abelian group $G_{k(*)}$ is a
Borel group,  even an $F_\sigma$-one that is the set of 
elements and the graphs of $+$ and the
function $x \mapsto -x$ (i.e. $\{(x,y,z):G_{k(*)} \models ``x + y =
z"\}$ hence also $\{(x,-x):x \in G_{k(*)}\}$) are $F_\sigma$-sets; hence
Borel.
\endproclaim
\bigskip

\demo{Proof}  Let cd be a one-to-one function from the set of
finite sequences of natural numbers onto the set of natural numbers
and we define:
\mr
\item "{$\odot_1$}"  $(a) \quad$ code$(x_{m,\bar\eta,\nu}) = 
\langle \text{cd}(\langle m,\text{cd}(\nu),\dotsc,\eta_\ell(i),
\ldots\rangle_{\ell \le k(*),\ell\ne m}:i < \omega\rangle$ 
so it
\nl

\hskip25pt  $\in {}^\omega \omega$ and
let ${\Cal X} = \{\text{code}(x_{n,\bar \eta,\nu}):(n,\bar \eta,\nu)$
as in Definition \scite{6.1}$\}$
\sn
\item "{${{}}$}"  $(b) \quad$ code$(y_{\bar \eta,n}) = \langle
\text{cd}(n,\dotsc,\eta_\ell(i)\ldots)_{\ell \le n(*)}:i <
\omega\rangle$ and 
\nl

\hskip25pt ${\Cal Y} = \{\text{code}(y_{\bar \eta,n}):(\bar
\eta,n)$ as in Definition \scite{6.1}$\}$
\sn
\item "{${{}}$}"  $(c) \quad$ for a 
sequence $\bar a = \langle a_\ell:\ell < n\rangle$ of
integers let $\rho_{\bar a} = \langle\text{sign}(a_\ell)$:
\nl

\hskip25pt $\ell <
n\rangle$, sign$(a_\ell)$ is 0,1,2 if $a_\ell$ is negative, zero,
positive respectively.
\ermn
We say $\nu$ represents $x \in G_{k(*)}$ as witnessed by $\langle
(z_\ell,a_\ell,m):\ell < n\rangle$ \ub{when}:
\mr
\item "{$\odot_2$}" $(a) \quad G \models x = 
\dsize \sum_{\ell < n} a_\ell z_\ell$,
\sn
\item "{${{}}$}"  $(b) \quad z_\ell \in \{x_{n,\bar \eta,\nu}:(\eta,\bar
\eta,\nu)$ as in Definition \scite{6.1}$\} \cup \{y_{\bar \eta,m}:
(\bar \eta,m)$ as
\nl

\hskip25pt  in Definition \scite{6.1}$\}$
\sn
\item "{${{}}$}"  $(c) \quad \langle z_\ell:\ell < n\rangle$ is without
repetitions
\sn
\item "{${{}}$}"  $(d) \quad \langle \text{cd}(z_\ell) \restriction
m:\ell < n\rangle$ are pairwise distinct
\sn
\item "{${{}}$}"  $(e) \quad$ if $\langle(z'_\ell,a'_\ell,m'):\ell <
n'\rangle$ satisfies clauses (a)-(d), then $m \le m'$
\sn
\item "{${{}}$}"   $(f) \quad$ if $n=0$ then $m=0$
\sn
\item "{${{}}$}"  $(g) \quad \text{ cd}(z_0) <_{\text{lex}} \text{
cd}(z_1) <_{\text{lex}} \ldots$
\sn
\item "{${{}}$}"  $(h) \quad \nu = \langle \text{cd}(\langle n\rangle \char 94
\text{ sign}(\bar a) \char 94 \langle |a_\ell|:\ell < n\rangle \char
94 \langle \text{cd}(z_\ell)(i):\ell \le k(*)\rangle):i < \omega\rangle$ .
\ermn
Now for $n <\omega,\bar a = \langle a_\ell:\ell < n\rangle \in {}^n
\Bbb Z,i < \omega$ and $\bar\varrho = \langle \varrho_\ell:\ell <
n\rangle \in {}^n({}^i Z)$ is $<_{\text{lex}}$-increasing hence
without repetitions (and if $n=0$ we let $i=0$) we let

$$
\align
Z_{\bar a,\bar\varrho} = \{\nu:&\nu \text{ represent some } x \in G_{k(*)}
\text{ as witnessed by} \\
  &\langle(z_\ell,a_\ell):\ell < n\rangle \text{ and } 
\text{ cd}(z_\ell) = \varrho_\ell \text{ for } \ell < n\}.
\endalign
$$
\mn
Let ${\Cal Y}$ be the set of such pairs $(\bar a,\bar\rho)$
\mr
\item "{$(*)_1$}"   $\langle Z_{\bar a,\bar\varrho}:(\bar
a,\bar\varrho) \in {\Cal Y}\rangle$ is a sequence of pairwise disjoint
 closed subsets of ${}^\omega \omega$ 
\sn
\item "{$(*)_2$}"   every member of $G$ is represented by one and only
one member of $Z := \cup\{Z_{\bar a,\bar\varrho}:(\bar a,\bar\varrho)
\in {\Cal Y}\}$.
\ermn
[Why?  For any $i < n$ clearly $\{x_{m,\bar\eta,\nu}:(m,\bar\eta,\nu)$
as in Definition \scite{6.1}$\} \cup\{y_{\bar\eta,i}:\bar\eta \in
\Lambda\}$ generates freely a subgroup $G'_{k(*),i}$ of $G_{k(*)}$
such that the quotient $G_{k(*)}/G'_{k(*),i}$ is torsion.  The rest
should be clear, too.]
\mr
\item "{$(*)_3$}"  ${\Cal U} = \{(\nu_1,\nu_2,\nu_3):\nu_\ell$
represent $x_\ell \in G_{k(*)}$ for $\ell=1,2,3$ and $G_{k(*)} \models
``x_1 + x_2 = x_3"\}$ is the graph of a two-place function
\sn
\item "{$(*)_4$}"  for any $(\bar a_\ell,\bar\varrho_\ell) \in {\Cal
Y}$ for $\ell=1,2,3$ the set $\{(\nu_1,\nu_2,\nu_3)
 \in {\Cal U}:\nu_\ell \in Z_{\bar a_\ell,\bar\varrho_\ell}$ for $\ell=1,2,3\}$
 is a closed set.
\ermn
Clearly we are done.  \hfill$\square_{\scite{6.1.7}}$
\enddemo
\bn
As a warm up we note:
\proclaim{\stag{6.2} Claim}  $G_{k(*)}$ is an $\aleph_1$-free abelian
group.
\endproclaim
\bigskip

\demo{Proof}  Let $U \subseteq {}^\omega 2$ be countable (and
infinite) and define $G'_U$ like $G$
restricting ourselves to $\eta_\ell \in U$; by the L\"owenheim-Skolem
argument it suffices to
prove that $G'_U$ is a free abelian group.  List ${}^{k(*)+1}U$ without
repetitions as $\langle \bar \eta_t:t < \omega \rangle$, and choose $s_t <
\omega$ such that $[r < t \and \bar \eta_r \restriction k(*) = \bar \eta_t
\restriction k(*) \Rightarrow \emptyset = \{\eta_{t,k(*)} 
\restriction \ell:\ell \in [s_t,\omega)\} \cap 
\{\eta_{r,k(*)} \restriction \ell:\ell \in [s_r,\omega)\}]$.
\sn
Let
\medskip

$\qquad \quad Y_1 = \{x_{m,\bar \eta,\nu}:m < k(*),\bar \eta \in
{}^{k(*)+1 \backslash \{m\}}U
\text{ and } \nu \in {}^{\omega >}2\}$

$$
\align
Y_2 = \biggl\{ x_{m,\bar \eta,\nu}:&m = k(*),\bar \eta \in {}^{k(*)}U
\text{ and for no } t < \omega \text{ do we have} \\
  &\bar \eta = \bar \eta_t \restriction k(*) \and 
\nu \in \{\eta_{t,k(*)} \restriction \ell:s_t \le \ell < \omega\} \biggr\}
\endalign
$$

$\qquad \quad Y_3 = \{y_{\bar \eta_t,n}:t < \omega \text{ and } n \in
[s_t,\omega)\}$.
\mn
Now
\mr
\item "{$(*)_1$}"  $Y_1 \cup Y_2 \cup Y_3$ generates $G'_U$.
\ermn
[Why?  Let $G'$ be the subgroup of $G'_U$ which $Y_1 \cup Y_2 \cup Y_3$
generates.  First we prove by induction on $n < \omega$ that for $\bar \eta
\in {}^{k(*)}U$ and $\nu \in {}^n 2$ we have $x_{k(*),\bar
\eta,\nu} \in G'$.  If $x_{k(*),\bar \eta,\nu} \in Y_2$ this is clear;
otherwise, by the definition of $Y_2$ for some $\ell < \omega$ and $t
< \omega$ such that $\ell \ge s_t$ we have $\bar \eta = \bar \eta_t
\restriction k(*),\nu = \eta_{t,k(*)} \restriction \ell$.

Now
\mr
\item "{$(a)$}"   $y_{\bar \eta_{t,\ell +1}},
y_{{\bar \eta}_t,\ell}$  are in $Y_3 \subseteq G'$.
\ermn
Hence by the equation $\boxtimes_{\bar \eta,n}$ in Definition
\scite{6.1}, clearly $x_{k(*),\bar \eta,\nu} \in G'$.  So as $Y_1
\subseteq G' \subseteq G'_U$, all the generators of the 
form $x_{m,\bar \eta,\nu}$
with each $\eta_\ell \in U$ are in $G'$.  Also we have
\mr
\item "{$(b)$}"  $x_{m,\bar \eta_t \restriction \{i \le k(*),i \ne
m\},\nu}$ belong to $Y_1 \subseteq G'$ if $m < k(*)$.
\ermn
Now for each $t < \omega$ we prove that all the generators $y_{{\bar
\eta}_t,n}$ are in $G'$.
If $n \ge s_t$ then clearly $y_{\bar \eta_t,n} \in Y_3 \subseteq G'$.
So it suffices to prove this for
$n \le s_t$ by downward induction on $n$; for $n = s_t$ by an earlier
sentence, for $n < s_t$ by $\boxtimes_{\bar \eta,n}$.
The other generators are in this subgroup so we
are done.]
\mr
\item "{$(*)_2$}"  $Y_1 \cup Y_2 \cup Y_3$ generates $G'_U$ freely. 
\ermn
[Why?  Translate the equations.

Alternatively, let $\langle z_\alpha:\alpha < \alpha(*)\rangle$ list
the set of generators of $G'_U$ without repetition such that for some
increasing continuous $\langle \alpha_i:i \le \omega + \omega\rangle$
we have $\alpha_0 = 0,\alpha_{\omega + \omega} = \alpha(*)$ and
\mr
\item "{$(a)$}"   $\{z_\alpha:\alpha < \alpha_1\} = Y_1 \cup \cup Y_2
\cup Y_3$
\sn
\item "{$(b)$}"  $\{z_\alpha:\alpha \in [\alpha_{1+n},\alpha_{1+n+1})\}
= \{x_{k(*),\bar \eta,\nu} \in G'_U:x_{k(*),\bar\eta,\nu} \notin Y_2$
\nl

\quad and $\ell g(\nu) = n\}$
\sn
\item "{$(c)$}"  $\{z_\alpha:\alpha \in 
[\alpha_{\omega +r},\alpha_{\omega +r+1})\} =
\{y_{\bar \eta_t,n}:t < \omega,n < s_t$ and $r= s_t -n\}$.
\ermn
Now the proof above shows that:
\mr
\item "{$\circledast$}"  there is a one-to-one function from the set
$\Xi$ of equations defining $G'_U$ onto $[\alpha_1,\alpha_\omega)$
such that:
\nl
if the equation $\varphi$ is mapped to the ordinal $\alpha$ then: if
$z_\beta$ appears in the equation then $\beta \le \alpha$ and
$z_\alpha$ appears in the equation and its coefficient is $1$ or $-1$.
\ermn
This clearly suffices.]   \hfill$\square_{\scite{6.2}}$
\enddemo
\bn
Now systematically
\definition{\stag{6.3} Definition}  1) For $U \subseteq {}^\omega 2$
let $G_U$ be the subgroup of $G$ generated by \nl
$Y_U = \{y_{\bar
\eta,n}:\bar \eta \in {}^{k(*)+1}(U)$ and $n < \omega\} \cup
\{x_{m,\bar \eta,\nu}:m \le k(*)$ and $\bar \eta \in
{}^{(k(*)+1)\backslash\{m\}}(U)$ and $\nu \in{}^{\omega >} 2\}$.  Let
$G^+_U$ be the divisible hull of $G_U$ and $G^+ = G^+_{({}^\omega 2)}$.
\nl
2) For $U \subseteq {}^\omega 2$ and
finite $u \subseteq {}^\omega 2$
let $G_{U,u}$ be the subgroup \footnote{note that if $u=\{\eta\}$ then
$G_{U,u} = G_U$}  of 
$G$ generated by $\cup\{G_{U \cup (u \backslash \{\eta\})}:
\eta \in u\}$; and for $\bar \eta \in {}^{k(*)
\ge} U$ let $G_{U,\bar \eta}$ be the
subgroup of $G$ generated by $\cup \{G_{U \cup \{\eta_k:k < \ell
g(\bar \eta) \text{ and } k\ne \ell\}}:\ell < \ell g(\bar \eta)\}$.
\nl
3) For $U \subseteq {}^\omega 2$ let $\Xi_U = \{\text{the equation }
\boxtimes_{\bar \eta,n}:\bar \eta \in {}^{k(*)+1} U$ and $n <
\omega\}$.  Let $\Xi_{U,u} = \cup\{\Xi_{U \cup u \backslash
\{\beta\}}:\beta \in u\}$.
\enddefinition
\bigskip

\proclaim{\stag{6.4} Claim}  0) If $U_1 \subseteq U_2 \subseteq
{}^\omega 2$ \ub{then} $G^+_{U_1} \subseteq G^+_{U_2}
\subseteq G^+$.
\nl
1) For any $n(*) < \omega$, the abelian group
$G^+_U$ (which is a vector space over $\Bbb Q$), has the basis 
$Y^{n(*)}_{U_i} := \{y_{\bar \eta,n(*)}:\bar \eta \in {}^{k(*)+1}(U)\} 
\cup \{x_{m,\bar \eta,\nu}:m \le k(*),\bar \eta \in
{}^{(k(*)+1)\backslash\{m\}}(U)$ and $\nu \in {}^{\omega >}2\}$. 
\nl
2) For $U \subseteq {}^\omega 2$ the abelian group
$G_U$ is generated by $Y_U$ freely (as an abelian group) except the
set $\Xi_U$ of equations.
\nl
3) If $U_m \subseteq {}^\omega 2$ for $m < m(*)$ \ub{then} the
subgroup $G_{U_0} + \ldots + G_{U_{m(*)-1}}$ of $G$ is 
generated by $Y_{U_0} \cup Y_{U_1} \cup
\ldots \cup Y_{U_{m(*)-1}}$ freely (as an abelian group) except the
equations in $\Xi_{U_0} \cup \Xi_{U_1} \cup \ldots \cup \Xi_{U_{m(*)-1}}$
provided that
\mr
\item "{$\circledast$}"  if $\eta_0,\dotsc,\eta_{k(*)} \in
\cup\{U_m:m < m(*)\}$ such that 
\nl

$(\forall \ell \le k(*))
(\exists m < m(*))[\{\eta_0,\dotsc,\eta_1\}
\backslash \{\eta_\ell\} \subseteq U_m)$ 
\nl

\ub{then} for some $m < m(*)$
we have $\{\eta_0,\dotsc,\eta_{k(*)}\} \subseteq U_m$.
\ermn
4) If $U_\ell = U \backslash U'_\ell$ for $\ell < m(*) \le k(*)$ and
$\langle U'_\ell:\ell < m(*)\rangle$ are pairwise disjoint then
$\circledast$ holds.
\nl
5) $G_{U,u} \subseteq G_{U \cup u}$ if $U \subseteq {}^\omega 2$ and
$u \subseteq {}^\omega 2 \backslash U$; moreover $G_{U,u}
\subseteq_{\text{pr}} G_{U \cup u} \subseteq_{\text{pr}} G$.
\nl
6) If $\langle U_\alpha:\alpha < \alpha(*)\rangle$ is
$\subseteq$-increasing continuous \ub{then} also $\langle
G_{U_\alpha}:\alpha < \alpha(*)\rangle$ is $\subseteq$-increasing
   continuous.
\nl
7) If $U_1 \subseteq U_2 \subseteq U \subseteq {}^\omega 2$ and $u \subseteq
{}^\omega 2 \backslash U$ is finite, $|u| < k(*)$ and $U_2
\backslash U_1 = \{\eta\}$ and $v=u \cup\{\eta\}$ \ub{then}
$(G_{U,u} + G_{U_2 \cup u})/(G_{U,u} + G_{U_1 \cup u})$ is isomorphic to
$G_{U_1 \cup v}/G_{U_1,v}$.
\nl
8) If $U \subseteq {}^\omega 2$ and $u \subseteq {}^\omega 2
\backslash U$ has $\le k(*)$ members \ub{then} $(G_{U,u} + G_u)/G_{U,u}$
is isomorphic to $G_u/G_{\emptyset,u}$.
\endproclaim
\bigskip

\demo{Proof}  0), 1) Obvious. 
\nl
2),3),4)   Follows.
\nl
5) First, $G_{U,u} \subseteq G_{U \cup u}$ follows by the definition.
Second, we deal with proving $G_{U,u} \subseteq_{\text{pr}} G_{U \cup u}$.
So let $|u| = m(*)+1$ and $\langle \eta_\ell:\ell \le m(*)\rangle$
list $u$, necessarily with no repetitions and let $U_\ell =
U \cup (u \backslash \{\eta_\ell\})$ (so $G_{U,u} =
G_{U_0} + \ldots + G_{U_{m(*)}})$ and assume $z \in G_{U \cup u},a \in
\Bbb Z \backslash \{0\}$ and az belongs to $G_{U_0} + \ldots +
 G_{U_{m(*)}}$ so it has the form 
$\Sigma\{b_i x_{m_i,\bar \eta_i,\nu_i}:i < i(*)\} + 
\Sigma\{c_j y_{\bar\rho_j,n_j}:j < j(*)\}$ with $b_i,c_j 
\in \Bbb Z$ and $\bar \eta_i,\bar\rho_j$ are (finite) sequences of 
members of $U_{\ell(i)},U_{k(j)}$ respectively and are as required in
Definition \scite{6.1} where $\ell(i),k(j) < m(*)$. 

Now similarly
as $z \in G_{U \cup u}$, we can find $z = \Sigma\{b'_i x_{m'_i,\bar
\eta'_i,\nu'_i}:i < i'(*)\} + \Sigma\{c'_j y_{\bar\rho'_j,n'_j}:j<j'(*)\}$.

By the equations in Definition \scite{6.1} \wilog \, for some $n(*)$
we have: $i < i(*) \Rightarrow n_i = n(*)$ and $i < i'(*) \Rightarrow
n'_i = n(*)$.  Also \wilog \, in each of the sequences
$\langle(m_i,\bar\eta_i,\nu_i):i <
i(*)\rangle,\langle\bar\rho_j:j<j(*)\rangle$ is with no repetitions,
and also in
$\langle(m'_i,\eta'_i,\nu'_i):i<i(*)\rangle,\langle\rho'_j:j<j(*)\rangle$
there is no repetition (for $\langle \rho_j:j<j(*)\rangle$ and
$\langle \rho'_j:j<j'(*)\rangle$ we use $n_i = n(*),n'_i = n(*))$.  Together
\mr
\item "{$\circledast$}"  $\Sigma\{b_i
x_{m_i,\bar\eta_i,\nu_i}:i<i(*)\} + \Sigma\{c_j
y_{\bar\rho_j,n(*)}:j<j(*)\} = \Sigma\{a b'_i
x_{m'_i,\bar\eta'_i,\nu_i}:i<i(*)\} + \Sigma\{a c'_j
y_{\bar\rho'_j,n(*)}:j<j'(*)\}$. 
\ermn
Now this equation holds in $G_{U \cup u}$ hence is $G$ and even in
$G^+$.  By part (1) and the ``no repetitions", after possible
permuting we get $i(*) = i'(*),j(*) = j'(*),(m_i,\bar\eta_i,\nu_i) =
(m'_i,\bar\eta'_i,\nu'_i),b_i = ab'_i$ for $i<i(*),\bar\rho_j =
\rho'_j$ for $j<j(*),c_j = ac'_j$ for $j<j(*)$.  But this proves that
$\{x_{m'_i,\bar\eta'_i,\nu_i}:i<i'(*)\} \cup
\{y_{\bar\rho'_j,n(*)}:j<j'(*)\} \subseteq G_{U,u}$ hence $z \in
G_{U,u}$ as required.

Third, the proof of $G_{U \cup u} \subseteq_{\text{pr}} G$ is similar. 
\nl
6) Easy.
\nl
7)  Clearly $U_1 \cup v = U_2 \cup u$ hence 
$G_{U_1 \cup u} \subseteq G_{U_1 \cup v} = G_{U_2 \cup u}$ hence
$G_{U,u} + G_{U_1 \cup u}$ is a subgroup of $G_{U,u} + G_{U_2 \cup u}$,
so the first quotient makes sense.

Hence by the isomorphism theorem 
$(G_{U,u} + G_{U_2 \cup u})/(G_{U,u} + G_{U_1 \cup u})$ is isomorphic
to $G_{U_2 \cup u}/(G_{U_2 \cup u} \cap (G_{U,u} + G_{U_1 \cup u}))$. 
Now $G_{U_1,v} \subseteq G_{U_1 \cup v} = G_{U_2 \cup u}$ 
and $G_{U_1,v} = \Sigma\{G_{U_1 \cup(v \backslash \{\nu\})}:\nu \in v\}
= \Sigma\{G_{U_1 \cup v \backslash \{\nu\}}:\nu \in u\} + G_{U_1,(v
\backslash \{\eta\})} \subseteq G_{U,u} + G_{U_1,u}$. 
Together $G_{U_1,v}$ is included in their intersection, i.e.
$G_{U_2 \cup u} \cap (G_{U,u} + G_{U_1,u})$ include $G_{U_1,v}$ and
using part (1) both has the same divisible hull inside $G^+$.  But
$G_{U_1,v}$ is a pure subgroup of $G$ by part (5)
hence of $G_{U_1 \cup v}$.  Hence
necessarily $G_{U_1 \cup u} \cap (G_{U,u} + G_{U_1,u}) = G_{U_1,v}$, so
as $G_{U_2 \cup u} = G_{U_1 \cup v}$ we are done.
\nl
8) The proof is similar to the proof of part (7).  Note that $G_{U,u}
\subseteq G_{U,u} + G_u$ hence the first quotient makes sense.  So by
an isomorphism theorem $(G_{U,u} + G_u)/G_{U,u}$ is isomorphic to
$G_u/(G_{U,u} \cap G_u)$.  Now $G_{U,u} \cap G_u$ includes
$G_{\emptyset,u}$ and using part (1) both has the same divisible hull
inside $G^+$.  But $G_{\emptyset,u}$ is a pure subgroup of $G_u$ by
part (5).  So necessarily $G_{U,u} \cap G_u = G_{\emptyset,u}$, so
$G_u/(G_{U,u} \cap G_u) = G_u/G_{\emptyset,u}$, so we are done.
 \hfill$\square_{\scite{6.4}}$
\enddemo
\bn
\ub{Discussion}:  For the reader's convenience 
we write what the group $G_{k(*)}$ is
for the case $k(*)=0$.  So, omitting constant indexes and replacing
sequences of length one by the unique entry we get that it is generated by
$y_{\eta,n}$ (for $\eta \in{}^\omega 2,n < \omega$) and $x_\nu$ (for $\nu \in
{}^{\omega >} 2$) freely as an abelian group except the equations
$(n!)y_{\eta,n+1} = y_{\eta,n} + x_{\eta \restriction n}$. \nl
Note that if $K$ is the countable subgroup generated by $\{x_\nu:\nu
\in {}^{\omega >} 2\}$ then $G/K$ is a divisible group of cardinality
continuum hence $G$ is not free.  So $G$ is $\aleph_1$-free but not free.
\bn
Now we have the main proof
\proclaim{\stag{6.5} Main Claim}  1) The abelian group $G_{U \cup
u}/G_{U,u}$ is free \ub{if} $U \subseteq {}^\omega 2,u \subseteq {}^\omega
2 \backslash U$ and $1 \le |u| \le k \le k(*)$ and $|U| \le
\aleph_{k(*)-k}$. \nl
2) If $U \subseteq {}^\omega 2$ and $|U| \le \aleph_{k(*)}$, 
\ub{then} $G_U$ is free.
\endproclaim
\bigskip

\demo{Proof}  1) We prove this by induction on $|U|$; \wilog \, $|u|=k$
as also $k' = |u|$ satisfies the requirements.
\enddemo
\bn
\ub{Case 1}:  $U$ is countable.

So let $\{\nu^*_\ell:\ell < k\}$ list $u$ be with no repetitions, now if
$k=0$, i.e. $u = \emptyset$ then $G_{U \cup u} = G_U = G_{U,u}$ so the
conclusion is trivial.  Hence we assume $u \ne \emptyset$, and let $u_\ell :=
u \backslash \{\nu^*_\ell\}$ for $\ell < k$.

Let $\langle \bar \eta_t:t < t^* \le \omega \rangle$ list with no
repetitions the set $\Lambda_{U,u} := \{\bar \eta \in 
{}^{k(*)+1}(U \cup u)$: for no
$\ell < k$ does $\bar \eta \in {}^{k(*)+1}(U \cup u_\ell)\}$.  Now
comes a crucial point: let $t < t^*$, 
for each $\ell < k$ for some $r_{t,\ell} \le
k(*)$ we have $\eta_{t,r_{t,\ell}} = \nu^*_\ell$ by the definition
of $\Lambda_{U,u}$, so 
$|\{r_{t,\ell}:\ell < k\}| = k < k(*)+1$ hence for some $m_t \le
k(*)$ we have $\ell < k \Rightarrow r_{t,\ell} \ne m_t$ so for each
$\ell < k$ the sequence $\bar \eta_t \restriction (k(*)+1 \backslash 
\{m_t\})$ is not
from $\{\langle \rho_s:s \le k(*)$ and $s \ne m_t \rangle:
\rho_s \in {}^\omega(U \cup u_\ell)$ for every $s \le k(*)$
such that $s \ne m_t\}$.  

For each $t < t^*$ we define $S(t) = \{m \le k(*):\{\eta_{t,s}:s
\le k(*) \and s \ne m\}$ is included in $U \cup u_\ell$ for no
$\ell \le k\}$.  So $m_t \in S(t) \subseteq \{0,\dotsc,k(*)\}$ and $m
\in S(t) \Rightarrow \bar \eta_t \restriction \{j \le k(*):j \ne m\}
\notin (U \cup u_\ell)$ for every $\ell \le k$.  For 
$m \le k(*)$ let $\bar \eta'_{t,m} := \bar \eta_t
\restriction \{j \le k(*):j \ne m\}$ and $\bar \eta'_t := \bar \eta'_{t,m_t}$.
Now we can choose $s_t <
\omega$ by induction on $t$ such that
\mr
\item "{$(*)$}"  if $t_1 < t,m \le k(*)$ and $\bar \eta'_{t_1,m} = 
\bar \eta'_{t,m}$, then \nl
$\eta_{t,m} \restriction s_t \notin \{\eta_{t_1,m}
\restriction \ell:\ell < \omega\}$.
\ermn
Let $Y^* = \{x_{m,\bar \eta,\nu} \in G_{U \cup u}:x_{m,\bar \eta,\nu}
\notin G_{U \cup u_\ell}$ for $\ell < k\} \cup \{ y_{\bar \eta,n} \in
G_{U \cup u}:y_{\bar \eta,n} \notin G_{U \cup u_\ell}$ for $\ell <
k\}$. 
\nl
Let
\medskip

$\quad \qquad Y_1 = \{x_{m,\bar \eta,\nu} \in Y^*$: for no $t < t^*$ do we have
$m = m_t \and \bar \eta = \bar \eta'_t\}$.

$$
\align
Y_2 = \{x_{m,\bar \eta,\nu} \in Y^*:&\,x_{m,\bar \eta,\nu} \notin
Y_1
\text{ but for no} \\
 &\,t < t^* \text{ do we have } m = m_t \and \bar \eta = \bar \eta'_t
\and \\
  &\,\eta_{t,m_t} \restriction s_t \trianglelefteq \nu \triangleleft
\eta_{t,m_t}\}
\endalign
$$

$\qquad \quad Y_3 = \{y_{\bar \eta,n}:y_{\bar \eta,n} \in Y^*
\text{ and } n \in [s_t,\omega) \text{ for the } t < t^* \text{
such that } \bar \eta = \bar \eta_t\}$. 
\sn
Now the desired conclusion follows from
\mr
\item "{$(*)_1$}"  $\{y + G_{U,u}:y \in Y_1 \cup Y_2 \cup Y_3\}$
generates $G_{U \cup u} /G_{U,u}$
\sn
\item "{$(*)_2$}"  $\{y + G_{U,u}:y \in Y_1 \cup Y_2 \cup Y_3\}$
generates $G_{U \cup u} /G_{U,u}$ freely.
\endroster
\bn
\demo{Proof of $(*)_1$}  It suffices to check that all the generators
of $G_{U \cup u}$ belong to $G'_{U \cup u} =: \langle Y_1 \cup Y_2
\cup Y_3 \cup G_{U,u} \rangle_G$.

First consider $x = x_{m,\bar \eta,\nu}$ where $\eta \in {}^{k(*)+1}(U
\cup u),m < k(*)$ and $\nu \in {}^n 2$ for some $n < \omega$.  If $x \notin
Y^*$ then $x \in G_{U,u_\ell}$ for some $\ell < k$ but $G_{U \cup
u_\ell} \subseteq G_{U,u} \subseteq G'_{U \cup u}$ so we are done, hence
assume $x \in Y^*$.  If $x \in Y_1 \cup Y_2 \cup Y_3$ we are done so
assume $x \notin Y_1 \cup Y_2 \cup Y_3$.  As $x \notin Y_1$ for some
$t < t^*$ we have $m = m_t \and \bar \eta = \eta'_t$.  As $x \notin
Y_2$, clearly for some $t$ as above we have 
$\eta_{t,m_t} \restriction s_t \trianglelefteq
\nu \triangleleft \eta_{t,m_t}$.  Hence by Definition \scite{6.1} the
equation $\boxtimes_{{\bar \eta}_t,n}$ from Definition \scite{6.1}
holds, now $y_{{\bar \eta}_t,n},y_{{\bar \eta}_t,n+1} \in 
G'_{U \cup u}$.  So in order to deduce from the equation that 
$x = x_{m,\bar \eta,\nu} = 
x_{m_t,\bar \eta'_t,\eta_{t,m_t} \restriction n}$ belongs to $G_{U \cup u}$, it
suffices to show that $x_{j,\bar \eta'_{t,j},\eta_{t,j} \restriction n} 
\in G'_{U \cup u}$ for each $j \le k(*),j \ne m_t$.  But each such
$x_{j,\bar \eta'_{t,j},\eta_{t,j} \restriction n}$ belong to $G'_{U \cup u}$ 
as it belongs to $Y_1 \cup Y_2$. \nl
[Why?  Otherwise necessarily for some $r < t^*$ we have
$j = m_r,\bar \eta'_{t,j} = \bar \eta'_{r,m_r}$ 
and $\eta_{r,m_r} \restriction s_r \trianglelefteq \eta_t
\restriction n \triangleleft \eta_{r,m_r}$ so $n \ge s_r$ and as said
above $n \ge s_t$.  Clearly $r \ne t$ as $m_r
= j \ne m_t$, now as $\bar \eta'_{t,m_r} = \bar
\eta'_{r,m_r}$ and $\bar \eta_t \ne \bar \eta_r$ (as $t \ne r$) 
clearly $\eta_{t,m_r}
\ne \eta_{r,m_r}$.  Also $\neg(r < t)$ by $(*)$ above applied with
$r,t$ here standing for $t_1,t$ there  as $\eta_{r,m_r}
\restriction s_r \trianglelefteq \eta_{t,j} \restriction n \triangleleft
\eta_{r,m_r}$.  Lastly for if $t < r$, again $(*)$ applied with $t,r$
here standing for $t_1,t$ there as $n \ge m_t$ gives
contradiction.]
\nl
So indeed $x \in G'_{U \cup u}$.

Second consider $y = y_{\bar \eta,n} \in G_{U \cup u}$, if $y \notin
Y^*$ then $y \in G_{U,u} \subseteq G'_{U \cup u}$, so assume $y \in
Y^*$.  If $y \in Y_3$ we are done, so assume $y \notin Y_3$, so for
some $t,\bar \eta = \bar \eta_t$ and $n < s_t$.  We prove by downward
induction on $s \le s_t$ that $y_{\bar \eta,s} \in G'_{U \cup u}$, this
clearly suffices.  For $s=s_t$ we have $y_{\bar \eta,s} \in Y_3
\subseteq G'_{U \cup u}$; and if $y_{\bar \eta,s+1} \in G'_{U \cup u}$
use the equation $\boxtimes_{\bar \eta_t,s}$ from \scite{6.1}, in the 
equation
$y_{\bar \eta,s+1} \in G'_{U \cup u}$ and the $x$'s appearing in the
equation belong to $G'_{U \cup u}$ by the earlier part of the proof
(of $(*)_1$) so
necessarily $y_{\bar \eta,s} \in G'_{U \cup u}$, so we are done.
\enddemo
\bigskip

\demo{Proof of $(*)_2$}  We rewrite the equations in the new variables
recalling that $G_{U \cup u}$ is generated by the relevant variables
freely except the equations of $\boxtimes_{\bar \eta,n}$ from
Definition \scite{6.1}.  After rewriting, all the equations disappear.
\enddemo
\bn
\ub{Case 2}:  $U$ is uncountable.

As $\aleph_1 \le |U| \le \aleph_{k(*)-k}$, necessarily $k < k(*)$.

Let $U = \{\rho_\alpha:\alpha < \mu\}$ where $\mu = |U|$, 
list $U$ with no repetitions.
Now for each $\alpha \le |U|$ let $U_\alpha := 
\{\rho_\beta:\beta < \alpha\},u_\alpha = u \cup
\{\rho_\alpha\}$.  Now  
\mr
\item "{$\odot_1$}"  $\langle(G_{U,u} + G_{U_\alpha \cup u})/G_{U,u}:
\alpha < |U| \rangle$ is an 
increasing continuous sequence of subgroups of $G/G_{U,u}$
\nl
[Why?  By \scite{6.4}(6).]
\sn
\item "{$\odot_2$}"  $G_{U,u} + G_{U_0 \cup u}/G_{U,u}$ is free.
\nl
[Why?  This is $(G_{U,u} + G_{\emptyset \cup u})/G_{U,u} = (G_{U,u} +
G_u)/G_{U,u}$ which by \scite{6.4}(8) is isomorphic to
$G_u/G_{\emptyset,u}$ whichis free by Case 1.]
\ermn
Hence it suffices to prove that for each
$\alpha < |U|$ the group $(G_{U,u} + G_{U_{\alpha +1} \cup
u})/(G_{U,u} + G_{U_\alpha \cup u})$ is free.  But easily
\mr
\item "{$\odot_3$}"    this group is isomorphic to $G_{U_\alpha \cup
u_\alpha}/G_{U_\alpha,u_\alpha}$.
\nl
[Why?  By \scite{6.4}(7) with $U_\alpha,U_{\alpha +1},U,\rho_\alpha,u$
here standing for $U_1,U_2,U,\eta,u$ there.]
\sn
\item "{$\odot_4$}"  $G_{U_\alpha \cup u_\alpha}/G_{U_\alpha,u_\alpha}$
is free. 
\nl
[Why?  By the induction hypothesis, as $\aleph_0 + |U_\alpha| < |U|
\le \aleph_{k(*)-(k+1)}$ and $|u_\alpha| = k+1 \le k(*)$.] 
\ermn
2) If $k(*) = 0$ just use \scite{6.2}, so assume $k(*) \ge 1$.  Now
the  proof is similar to (but easier than) the proof of case (2)
inside the proof of part (1) above.
\nl
${{}}$  \hfill$\square_{\scite{6.5}}$
\bigskip

\proclaim{\stag{6.6} Claim}  If $U \subseteq {}^\omega 2$ and 
$|U| \ge \aleph_{k(*)+1}$ \ub{then} $G_U$ is not free.
\endproclaim
\bigskip

\demo{Proof}  Assume toward contradiction that $G_U$ is free and let
$\chi$ be large enough; for notational simplicity assume $|U| =
\aleph_{k(*)+1}$.  O.K. as a subgroup of a free abelian group is a
free abelian group.  We
choose $N_\ell$ by downward induction on $\ell \le k(*)$ such that
\mr
\item "{$(a)$}"  $N_\ell$ is an elementary submodel 
\footnote{${\Cal H}(\chi)$ is $\{x$: the transitive closure of $x$ has
cardinality $< \chi\}$ and $<^*_\chi$ is a well ordering of ${\Cal H}(\chi)$}
of $({\Cal H}(\chi),\in,<^*_\chi)$
\sn
\item "{$(b)$}"  $\|N_\ell\| = |N_\ell \cap \aleph_{k(*)}| =
\aleph_\ell$ and $\aleph_\ell +1 \subseteq N_\ell$
\sn
\item "{$(c)$}"   $G_U \in N_\ell$ and $N_{\ell +1},\dotsc,N_{k(*)}
\in N_\ell$.
\ermn
Let $G_\ell = G_U \cap N_\ell$, a subgroup of $G_U$.
Now
\mr
\item "{$(*)_0$}"  $G_U/(\Sigma\{G_\ell:\ell \le k(*)\})$ 
is a free (abelian) group 
\nl
[easy or see \cite{Sh:52}, that is: \nl
as $G_U$ is free we can prove by induction on $k \le k(*) + 1$ then
$G/(\Sigma\{G_\ell:\ell < k\})$ is free, for $k = 0$ this is the
assumption toward contradiction, the induction step is by Ax VIII in
\cite{Sh:52} for abelian groups and for $k=k(*)+1$ we get the desired 
conclusion.] 
\ermn
Now
\mr
\item "{$(*)_1$}"  letting $U^1_\ell$ be $U$ for $\ell = k(*)+1$ and
$\dbca^{k(*)}_{m=\ell} (N_m \cap U)$ for $\ell \le k(*)$; we have: $U^1_\ell$
has cardinality $\aleph_\ell$ for $\ell \le k(*) + 1$ \nl
[Why?  By downward induction on $\ell$.
For $\ell = k(*)+1$ this holds by an assumption.  For 
$\ell = k(*)$ this holds by clause (b).  For $\ell < k(*)$ this
holds by the choice of $N_\ell$ as the
set $\dbca^{k(*)}_{m=\ell +1} (N_m \cap U)$ has cardinality
$\aleph_{\ell +1} \ge \aleph_\ell$ and belong to $N_\ell$ and clause
(b) above.]
\sn
\item "{$(*)_2$}"  $U^2_\ell =: U^1_{\ell +1} \backslash (N_\ell \cap U)$
has cardinality $\aleph_{\ell +1}$ for $\ell \le k(*)$ \nl
[Why?  As $|U^1_{\ell +1}| = \aleph_{\ell +1} > \aleph_\ell =
\|N_\ell\| \ge |N_\ell \cap U|$.]
\sn
\item "{$(*)_3$}"  for $m < \ell \le k(*)$ the set 
$U^3_{m,\ell} =: U^2_\ell \cap
\dbca^{\ell-1}_{r = m} N_r$ has cardinality $\aleph_m$   \nl
[Why?  By downward induction on $m$.  For $m=\ell -1$ as $U^2_\ell \in
N_m$ and $|U^2_\ell| = \aleph_{\ell +1}$ and clause (b).  
For $m < \ell$ similarly.]
\ermn
Now for $\ell=0$ choose $\eta^*_\ell \in U^2_\ell$, possible by
$(*)_2$ above.  Then for $\ell >0,\ell \le k(*)$ choose 
$\eta^*_\ell \in U^3_{0,\ell}$.  This is possible by $(*)_3$.  So
clearly
\mr
\item "{$(*)_4$}"  $\eta^*_\ell \in U$ and $\eta^*_\ell \in N_m \cap U
\Leftrightarrow \ell \ne m$ for $\ell,m \le k(*)$. \nl
[Why?  If $\ell=0$, then by its choice, $\eta^*_\ell \in U^2_\ell$,
hence by the definition of $U^2_\ell$ in $(*)_2$ we
have $\eta^*_\ell \notin N_\ell$, and 
$\eta^*_\ell \in U^1_{\ell +1}$ hence $\eta^*_\ell \in N_{\ell +1}
\cap \ldots \cap N_{k(*)}$ by $(*)_1$ so $(*)_4$ holds for $\ell = 0$.
If $\ell > 0$ then by its choice, $\eta^*_\ell \in U^3_{0,\ell}$ but
$U^3_{m,\ell} \subseteq U^2_\ell$ by $(*)_3$ so 
$\eta^*_\ell \in U^2_\ell$ hence as
before $\eta^*_\ell \in N_{\ell +1} \cap \ldots \cap N_{k(*)}$ and
$\eta^*_\ell \notin N_\ell$.  Also by $(*)_3$ we have
$\eta^*_\ell \in \dbca^{\ell -1}_{r=0}
N_\ell$ so $(*)_4$ really holds.]
\ermn
Let $\bar \eta^* = \langle \eta^*_\ell:\ell \le k(*) \rangle$ and let $G'$
be the subgroup of $G_U$ generated by $\{x_{m,\bar \eta,\nu}:m \le
k(*)$ and $\bar \eta \in {}^{k(*)+1 \backslash \{m\}}U$
and $\nu \in {}^{\omega >} 2\} 
\cup \{y_{\bar \eta,n}:\bar \eta \in {}^{k(*)+1} U 
\text{ but } \bar \eta \ne \bar \eta^*$ and $n < \omega\}$.
Easily $G_\ell \subseteq G'$ recalling $G_\ell = N_\ell \cap G_U$
hence $\Sigma\{G_\ell:\ell \le k(*)\}
\subseteq G'$, but $y_{\bar \eta^*,0} \notin G'$ hence
\mr
\item "{$(*)_5$}"  $y_{\bar \eta^*,0} \notin \sum\{G_\ell:\ell \le
k(*)\}$. 
\ermn
But for every $n$
\mr
\item "{$(*)_6$}"  $\bar n!y_{\bar \eta^*,n+1} - y_{\bar \eta^*,n} = 
\Sigma\{x_{m,\bar \eta^* \restriction (k(*))+1 \backslash \{m\}),\eta^*_m
\restriction n}:m \le k(*)\} \in \Sigma\{G_\ell:\ell \le k(*)\}$. \nl
[Why?  $x_{m,\bar \eta^* \restriction (k(*))+1 \backslash \{m\}),\eta^*_m
\restriction n} \in G_m$ as
$\bar \eta^* \restriction (k(*))+1 \backslash \{m\}) \in N_m$ by
$(*)_4$.]
\ermn
We can conclude that in $G_U/\sum\{G_\ell:\ell \le k(*)\}$, the
element $y_{\bar \eta^*,0} + \sum\{G_\ell:\ell \le k(*)\}$ is not zero
(by $(*)_5$) but is divisible by every natural number by $(*)_6$. \nl
This contradicts $(*)_0$ so we are done. \hfill$\square_{\scite{6.6}}$
\enddemo
\bigskip

\demo{\stag{6.6G} Conclusion}  $G_{k(*)}$ is a Borel and even $F_\sigma$  
abelian group which is $\aleph_{k(*)+1}$-free \ub{but} if $2^{\aleph_0} \ge
\aleph_{k(*)+1}$ is not free and even not $\aleph_{k(*)+2}$-free. 
\enddemo
\bigskip

\demo{Proof}  $G_{k(*)}$ is Borel and $F_\sigma$ by \scite{6.1.7}, it is
$\aleph_{k(*)+1}$-free by \scite{6.5} and if $2^{\aleph_0} \ge
\aleph_{k(*)+1}$ it is not $\aleph_{k(*)+2}$-free by \scite{6.6}.
\hfill$\square_{\scite{6.6G}}$ 
\enddemo
\bn
Blass asks 
\nl
\margintag{6.10}\ub{\stag{6.10} Question}:  Suppose (a) + (b) below, 
does it follow that forcing with $\Bbb Q$ add reals?
\mr
\item "{$(a)$}"   $G$ is a Borel definition of an
uncountable abelian subgroup of ${}^\omega \Bbb Z$ (the Specker group)
which is not free
\sn
\item "{$(b)$}"   the forcing $\Bbb Q$ satisfies 
$\Vdash_{\Bbb Q} ``G^{\bold V}$ is free".  
\endroster
\bn
Now 
\nl
\margintag{6.11}\ub{\stag{6.11} Fact}:  For just Borel abelian group $G$: if CH, then
the answer to \scite{6.10} is yes, if not CH then the answer is not
for $\Bbb Q = \text{ Levy}(\aleph_1,2^{\aleph_0})$.
\bigskip

\demo{Proof}  First, assume CH holds and $G$ is as in (a) of \scite{6.10};
(or just defined absolutely enough such that $G^{\bold V}$ is a
subgroup of $G^{\bold V^{\Bbb Q}}$ for any forcing notion $\Bbb Q$ and
is still not free).  
\ub{Then} by \cite{Sh:402} the group $G^{\bold V}$ is non-free
in some strong way such that no forcing not collapsing $2^{\aleph_0}$ to
$\aleph_0$ can make it free (that is, for some countable $G_0 \subseteq
G^{\bold V},G^{\bold V}_0/G$ contains the direct sum of $2^{\aleph_0}$
finite rank non-free abelian groups).
\sn
This is a strong yes answer. 
\nl
On the other hand, if $2^{\aleph_0} > \aleph_1$ we
can find such group: for $k(*) \ge 1$, our 
$G_{k(*)}$ if $\aleph_1 < \aleph_{k(*)+1} \le
2^{\aleph_0}$, see below, is a strong negative answer.  So together this gives
answers to a question of Blass.  \hfill$\square_{\scite{6.11}}$
\enddemo
\bigskip

\proclaim{\stag{6.12} Corollary}  1) The group $G_{k(*)}$ is embeddable
into ${}^\omega \Bbb Z$, even purely. \nl
2) Hence forcing which does not add bounded subsets to $\aleph_{k(*)}$
can make it free (i.e. {\rm Levy}$(\aleph_\ell,2^{\aleph_0})$ if $\ell \le
k(*)$  while if our universe satisfies $2^{\aleph_0} >
\aleph_{k(*)}$ it is not free there).
\endproclaim
\bigskip

\demo{Proof}  1) For every $n < \omega$ we define a function $f_n$
from $Y$ to $G_{k(*)}$ where $Y$ is the set of generators of
$G_{k(*)}$, i.e.

$$
\align
Y = \{y_{\bar \eta,n+1}:&n < \omega,\bar \eta \in {}^{k(*)+1}({}^\omega
2)\} \cup \{x_{m,\bar \eta_m,\nu}:m \le k(*), \\
  &\bar \eta \in {}^{\{\ell:\ell \le k(*),\ell \ne n\}}({}^\omega 2) 
  \text{ and } \nu \in {}^{\omega >} 2\}.
\endalign
$$
\mn
First define a function $h_n$: for 
$\eta \in {}^{\omega \ge} 2,g_n(\eta)$ is a sequence
of length $\ell g(\eta)$ and \nl
\sn
$(h_n(\eta))(\ell) = \cases \eta(\ell)
&\text{ if } \ell < n \and \ell < \ell g(\eta) \\ 
0 &\text{ if } \ell \ge n \and \ell < \ell g(\eta) \endcases$ . \nl
\sn
For $\bar \eta = \langle \eta_\ell:\ell \in u \rangle \in
{}^u({}^{\omega \ge} 2)$ we let $h_n(\bar \eta) = \langle f_n(\eta_\ell):\ell
\in u \rangle$. \nl
Lastly, let

$$
f_n(y_{{\bar \eta},n+1}) = y_{h_n({\bar \eta}),n+1}
$$

$$
f_n(x_{m,{\bar \eta}_m,\nu}) = x_{m,h_n({\bar \eta}_m),h_n(\nu)}.
$$
\mn
Does $f_n$ induce a homomorphism from $G_{k(*)}$ into $G_{k(*)}$?  For
this it is enough to check that for every one of the relations from
Definition \scite{6.1}, its $f_n$-image is satisfied in $G_{k(*)}$,
but this is obvious as it is mapped to another one of the equations
in the definition of $G_{k(*)}$: the equation in $\boxtimes_{\bar
\eta,m}$ is mapped to the equation in $\boxtimes_{g_n(\bar \eta),m}$.

So $f_n$ extends to an endomorphism $\hat f_n$ of $G_{k(*)}$.  Easily
\mr
\item "{$\circledast$}"  if $L \subseteq G_{k(*)}$ is a finite rank
subgroup (so free) then for $n$ large enough $\hat f \restriction L$ is one
to one.
\ermn
Now the range of $\hat f_n$ is clearly countable hence free, say is
$\underset{\ell < \omega} {}\to \bigoplus \Bbb Z z_{n,\ell}$.  Hence for
some homomorphisms $g_{n,\ell}$ from Range$(f_n)$ to $\Bbb Z$ for
$\ell < \omega$  we have

$$
\align
z \in \text{ Rang}(\hat f_n) \Rightarrow &z =
\Sigma\{g_{n,\ell}(z)z_{n,\ell}:\ell < \omega\} \\
  &\text{ where } g_{n,\ell}(z)=0 \text{ for every } \ell \text{ large
enough}
\endalign
$$
\mn
Let $f_{n,\ell} = g_{n,\ell} \circ \hat f_n \in \text{
Hom}(G_{k(*)},\Bbb Z)$.  Those homomorphisms give, by renaming the
$f_{n,\ell}$'s, an embedding of $G_{k(*)}$ into ${}^\omega \Bbb Z$.
Looking at the construction, it is a pure one. \nl
2) By \scite{6.5}.    \hfill$\square_{\scite{6.12}}$
\enddemo
\bigskip

\proclaim{\stag{6.16} Claim}  Assume {\rm MA} $+ 2^{\aleph_0} >
\aleph_2$.

If $k(*) > 2$ \ub{then} $G = G_{k(*)}$ is a Whitehead Borel (abelian) group. 
\endproclaim
\bigskip

\demo{Proof}  By \scite{6.1.7} we know that $G_{k(*)}$ is a Borel
group.  Let $\langle \eta_\alpha:\alpha < 2^{\aleph_0}\rangle$ list
${}^\omega 2$ with no repetitions and ${\Cal U}_\alpha =
\{\eat_\beta:\beta < \alpha\}$.

So $\langle {\Cal U}_\alpha:\alpha < 2^{\aleph_0}\rangle$ be
$\subseteq$-increasing continuous with union ${}^\omega 2$ such that
${\Cal U}_0 = \emptyset,|{\Cal U}_\alpha| \le |\alpha|$; and let
$H_\alpha := G_{{\Cal U}_\alpha}$, see Definition \scite{6.3}(1).  So
$\langle H_\alpha:\alpha < 2^{\aleph_0}\rangle$ is a
$\subseteq$-increasing continuous sequence of subgroups of $G$ with
union $G$.  For $\alpha < 2^{\aleph_0}$, letting $u_\alpha =
\{u_\alpha\}$ recalling Definition \scite{6.3} we hae $G_{{\Cal
U}_\alpha \cup u_\alpha} = G_{{\Cal U}_{\alpha +1}} = H_{\alpha +1}$
and $G_{U_\alpha,u_\alpha} = G_{{\Cal U}_\alpha} = H_\alpha$, hence
$H_{\alpha +1}/H_\alpha = G_{{\Cal U}_\alpha \cup u_\alpha}/G_{{\Cal
U}_\alpha,u_\alpha}$ and by \scite{6.5}(1) the latter group is
$\aleph_2$-free so $H_{\alpha +1}/H_\alpha$ is $\aleph_2$-free.  As MA
holds and $|H_{\alpha +1}/H_\alpha| < 2^{\aleph_0}$ and $H_{\alpha
+1}/H_\alpha$ is $\aleph_2$-free we know that it is a Whitehead group.

As $H_\alpha$ is $\subseteq$-increasing continuous, $H_0 =\{0\}$ and
each $H_{\alpha +1}/H_\alpha$ is a Whitehead group, it follows that
$\cup\{H_\alpha:\alpha < 2^{\aleph_0}\}$ is a Whitehead group, which
means that $G$ is as required.  \hfill$\square_{\scite{6.16}}$
\enddemo
\newpage

\head{\S6 Beginning of stability theory} \endhead  \resetall \sectno=6
 \spuriousreset
\bigskip

We may consider the dividing line for abelian groups 
from \cite{Sh:402} and try to
generalize it for any simply defined (e.g. $\Sigma^1_1$ or Borel) model.  We
deal with having two possibilities, in the high, complicated
side we get a parallel of non $\aleph_0$-stability; in the low side we
have a rank.  But even for minimal formulas, the example in \S5 shows
that we are far from being done, still we may be able to say something on the
structure.

We may consider also ranks parallel to the ones for superstable
theories.  Note that there are two kinds of definability we are
considering: the model theoretic one and the set theoretic one.  See
more in \cite{Sh:F562}.
\bigskip

\demo{\stag{7.0} Convention}  If not said otherwise, ${\frak A}$ will
be a structure with countable vocabulary and its set of elements is a
set of reals.
\enddemo
\bigskip

\definition{\stag{6a.deL} Definition}  1)  For a structure ${\frak A}$,
an ${\frak A}$-formula $\varphi$ is a formula in the vocabulary of
${\frak A}$ with finitely many free variables, writing $\varphi =
\varphi(\bar x)$ means that $\bar x$ is a finite sequence of variables with no
repetitions including the free variables of $\varphi$. 
We did not specify the logic; we may assume it is 
$\subseteq \Bbb L_{\omega_1,\omega}$ or even 
$\Bbb L_{\omega_1,\omega}(\bold Q)$ where $\bold Q$ is the quantifier ``there
are uncountably many".
\nl
2) $\Delta$ denotes a set of such formulas and $\bar \varphi$ a pair
$(\varphi_0(\bar x),\varphi_1(\bar x))$ of formulas so $\bar \varphi$
is a $\Delta$-pair if $\varphi_0,\varphi_1 \in \Delta$. \nl
3) We say $\varphi$ (or $\Delta$ or $\bar \varphi$) is $\Sigma^1_1$ (or
$\Sigma^1_2$ or $\Delta^1_0$ (= Borel)) \ub{if} 
they are so as set theoretic formulas.
\enddefinition
\bigskip

\definition{\stag{7.1} Definition}  1) We say $({\frak A},\Delta)$ is a
$\Sigma^1_1$-candidate \ub{when}:
\mr
\item "{$(a)$}"  ${\frak A}$ is a $\Sigma^1_1$-model
\sn
\item "{$(b)$}"  $\Delta$ is a countable set of ${\frak A}$-formulas
which, are in the set theory sense, $\Sigma^1_1$ (we 
identify $\varphi$ and $\neg \neg \varphi$).
\ermn
We can replace being $\Sigma^1_1$ by $\Sigma^1_2$, etc., (naturally we
need enough absoluteness); if we
replace it by $\Gamma$ we write $\Gamma$-candidate.  If $\Gamma$ does
not appear we mean it is $\Sigma^1_1$ or understood from the content
normal.
\nl
2) If $({\frak A},\Delta)$ is a candidate we say ${\frak A}$ is locally
$(\aleph_0,\Delta)$-stable (or $({\frak A},\Delta)$ is
$\aleph_0$-stable), but we may omit ``locally"; \ub{when} $\Delta$ is
a countable set of ${\frak A}$-formulas and for 
$\chi$ large enough and $x \in {\Cal H}(\chi)$, 
for every countable $N \prec ({\Cal H}(\chi),\in,<^*_\chi)$ to which
$x$ belongs and $\bar a \in {}^m {\frak A}$ where $m < \omega$ 
the following weak definability condition on tp$_\Delta(\bar a,N
\cap {\frak A},{\frak A})$ holds:
\mr
\item "{$(*)$}"   letting $\Phi^m_{{\frak A},\Delta} =
\Phi^n_{({\frak A},\Delta)} = 
\{\bar \varphi(\bar x,\bar b):\bar \varphi(\bar x,\bar b) =
(\varphi_0(\bar x,\bar b),
\varphi_1(\bar x,\bar b))$ and $\bar x = \langle x_\ell:\ell < m
\rangle,\bar b \in {}^{\omega >} {\frak A},\varphi_0,\varphi_1 \in
\Delta$ and ${\frak A} \models \neg(\exists \bar x)
(\varphi_0(\bar x,\bar b) \and
\varphi_1(\bar x,\bar b))\}$, \nl
for some function $\bold c \in N$ with domain
$\Phi^m_{({\frak A},\Delta)}$ to $\{0,1\}$ we have:
\sn
\item "{$(**)$}"  if $\bar \varphi = (\varphi_0(\bar x,\bar
b),\varphi_1(\bar x,\bar b)) \in \Phi^m_{({\frak A},\Delta)} \cap N$
and $\ell < 2$ and ${\frak A} \models \varphi_\ell[\bar a,\bar b]$
then $\ell = \bold c(\bar \varphi)$.
\ermn
3) We say that $({\frak A},\Delta)$ is $\aleph_0$-unstable (or ${\frak
A}$ is $(\aleph_0,\Delta)$-unstable) \ub{if}:
there are $\bar a_\eta \in {}^m {\frak A}$ for $\eta \in {}^\omega 2$ 
and $\varphi_{\nu,0}(\bar x,\bar y_\nu) \in \Delta$ and 
$\varphi_{\nu,1}(\bar x,\bar y_\nu) \in 
\Delta$ and $\bar b_\nu \in {}^{\ell g(\bar y)} {\frak A}$ for
$\nu \in {}^{\omega >} 2$ such that:
\mr
\item "{$(a)$}"  ${\frak A} \models \neg(\exists \bar x)
(\varphi_{\nu,0}(\bar x,\bar b_\nu) \and \varphi_{\nu,1}(\bar x,\bar b_\nu))$
\sn
\item "{$(b)$}"   if 
$\nu \triangleleft \eta_0,\nu \triangleleft \eta_1,n = \ell
g(\nu)$ and $\eta_0(n) = 0,\eta_1(n)=1 \text{ \ub{then} } 
{\frak A} \models \varphi_{\nu,0}(\bar a_{\eta_0},\bar b_\nu) 
\wedge \varphi_{\nu,1}(\bar a_{\eta_1},\bar b_\nu)$.
\ermn
There are obvious absoluteness results (for $\bar \varphi \in
\Phi^m_{({\frak A},\Delta)},({\frak A},\Delta)$ is $\aleph_0$-unstable
and stable).
\enddefinition
\bigskip

\remark{\stag{7.2} Observation}  1) If $\Delta$ is 
closed under negation \ub{then} in Definition \scite{7.1}(2) we have
\mr
\item "{$(*)'$}"  for some $\bold c \in N$ we have: $\varphi(\bar x,\bar
y) \in \Delta \and \bar b \in {}^{\ell g(\bar y)}{\frak A} \and \bar b
\in N$ implies 
\sn
\item "{$(**)'$}"  ${\frak A} \models \varphi(\bar a,\bar b)$ iff
$\bold c(\varphi(\bar x,\bar b))=1$.
\ermn
2) In Definition \scite{7.1}(2) we can fix $x = ({\frak A},\Delta)$ and
omit $<^*_\chi$, at the expense of larger $\chi$.
\endremark
\bigskip

\demo{Proof}  Straight.
\enddemo
\bigskip

\proclaim{\stag{7.4} The End-Extention Indiscernibility existence
lemma}   Assume $({\frak A},\Delta)$ is an
$\aleph_0$-stable candidate. \nl
1) In Definition \scite{7.1}(2), the demand ``$N$ is countable" can be
omitted. \nl
2) Assume $\Delta$ is closed under negation and permuting the
variables, $m < \omega,\bar a_\alpha \in {}^m {\frak A}$ for $\alpha <
\lambda$ and $\aleph_0 < \lambda = \text{\rm cf}(\lambda)$ and $S
\subseteq \lambda$ is stationary and $A \subseteq {\frak A}$ has
cardinality $< \lambda$.  \ub{Then} for some stationary $S' \subseteq
S$ the sequence $\langle \bar a_\alpha:\alpha \in S' \rangle$ is a
$\Delta$-end extension indiscernible sequence over $A$ in ${\frak A}$
(see Definition \scite{7.4A}(4),(5) below). \nl
3) Moreover for any pregiven $n < \omega$ we can find stationary $S'
\subseteq S$ such that $\langle \bar a_\alpha:\alpha \in S' \rangle$
is $(\Delta,n)$-end extension indiscernible over $A$ in ${\frak
A}$. \nl
4) We can find a club $E$ of $\lambda$ and regressive function $f_n$
on $S \cap E$ for $n < \omega$ such that:
\mr
\widestnumber\item{$(ii)^+$}
\item "{$(i)$}"   if $\alpha,\beta \in S \cap E$ \ub{then}
$f_{n+1}(\alpha) = f_{n+1}(\beta) \Rightarrow f_n(\alpha) =
f_n(\beta)$
\sn
\item "{$(ii)$}"   if $n < \omega$ and $\gamma < \lambda$, \ub{then}
the sequence $\langle \bar a_\alpha:\alpha
\in S \cap E,f_n(\alpha) = \gamma \rangle$ is $(\Delta,n)$-end
extension indiscernible over $A$
\sn
\item "{$(ii)^+$}"  moreover, if $n < \omega$ and $\beta,\gamma < \lambda$
\ub{then} $\langle \bar a_\alpha:\alpha  \in S \cap E \backslash \beta \and
f_n(\alpha) = \gamma \rangle$ is $(\Delta,n)$-end extension
indiscernible over $A \cup \{\bar a_\gamma:\gamma < \beta\}$.
\endroster
\endproclaim
\bigskip

\remark{Remark}  Similar to \cite[III,4.23,pg.120-1]{Sh:c}, but before
proving we define:
\endremark
\smallskip
\definition{\stag{7.4A} Definition}  1) Let $({\frak A},\Delta)$ be a
candidate.  We say ``${\frak A}$ has $(\lambda,\Delta)$-order" when:
\mr
\item "{$(*)_\lambda$}"  for some $m(*) < \omega$ and 
$\bar \varphi(\bar x,\bar y) \in \Phi^{m(*)}_{{\frak A},\Delta}$ with 
$\ell g(\bar x) = \ell g(\bar y)$ linear orders  some 
${\Cal I} \subseteq {}^{m(*)} {\frak A}$ of cardinality $\lambda$, see
part (2) for definition. 
\ermn
2) We say $\bar \varphi(\bar x,\bar y)$ linear orders $\bold I
\subseteq {}^{m(\lambda)} {\frak A}$ \ub{if} for some $\langle \bar a_t:t
\in I \rangle$ we have:
\mr
\item "{$(a)$}"  $\bold I = \langle \bar a_t:t \in I \rangle$
\sn
\item "{$(b)$}"  $I$ is a linear order
\sn
\item "{$(c)$}"  $\bar \varphi = (\varphi_0(\bar x,\bar
y),\varphi_1(\bar x,\bar y))$ and contradictory in ${\frak A}$
\sn
\item "{$(d)$}"  if $s <_I t$ then ${\frak A} 
\models \varphi_0(\bar a_s,\bar a_t) \wedge \varphi_1[\bar a_t,\bar a_s]$.
\ermn
3) For a linear order $I$ (e.g. a set of ordinals), we say $\langle
\bar a_t:t \in J \rangle$ is a $\Delta$-end-extension indiscernible
(sequence over $A$) \ub{if} for 
any $n < \omega$ and $t_0 <_J < \ldots <_J t_{n-1}
<_J t$, the sequences $\bar a_{t_0} \char 94 \ldots \bar
a_{t_{n-2}} \char 94 \bar a_{t_{n-1}}$ and $\bar a_{t_0} \char 94
\ldots \char 94 \bar a_{t_{n-2}} \char 94 \bar a_{t_n}$ realizes the
same $\Delta$-type (over $A$) in ${\frak A}$. \nl
4) We say that $\langle \bar a_t:t \in J \rangle$ is
$(\Delta,n_0,n_1)$-end-extension indiscernible over $A$ in ${\frak A}$
when:
\mr
\item "{$(a)$}"  $J$ a linear order for some $m,\bar a_t \in
{}^m{\frak A},A \subseteq {\frak A}$
\sn
\item "{$(b)$}"  if $\langle r_\ell:\ell < n_0 \rangle,\langle
s_\ell:\ell < n_1 \rangle,\langle t_\ell:\ell < n_1 \rangle$ are
$<_J$-increasing sequences, $r_{n_0-1} <_J s_0,r_{n_0-1} <_J
t_0$ \ub{then} $\bar a_{r_0} \char 94 \ldots \char 94 \bar a_{r_{n_0-1}}
\char 94 \bar a_{s_0} \char 94 \ldots \char 94 \bar a_{s_{n_1-1}}$ and
\nl
$\bar a_{r_0} \char 94 \ldots \char 94 \bar a_{r_{n_0-1}} \char 94
\bar a_{t_0} \char 94 \ldots \char 94 \bar a_{t_{n_1-1}}$ realizes the
same $\Delta$-type over $A$ in ${\frak A}$
\sn
\item "{$(c)$}"  if $J$ has a last element we allow to decrease $n_0$
and/or $n_1$.
\ermn
5) If we omit $n_0$ this means for every $n_0$, (so ``$\Delta$-end
extension..." means $(\Delta,1)$-end extension.
\enddefinition
\bigskip

\demo{Proof of \scite{7.4}}  1) Let $N^* \prec 
({\Cal H}(\chi),\in,<^*_\chi)$ be such
that ${\frak A},\Delta \in N^*$.  Now for every countable $N \prec N^*$ to
which $({\frak A},\Delta)$ belongs there is 
$\bold c_N \in N$ as mentioned in the definition 
\scite{7.2}(2).  Hence by normality of the club filter on
$[N^*]^{\aleph_0}$, the family of countable subsets of $N^*$,
for some $\bold c^*$ the set $\bold N =
\{N:N \prec N^*$ is countable and $\bold c_N = \bold c^*\}$ is a stationary
subset of $[N^*]^{\aleph_0}$, so $\bold c^*$ can serve for $N$. \nl
2) Let $\langle N_\alpha:\alpha < \lambda \rangle$ be an increasing
continuous sequence of elementary submodels of $({\Cal
H}(\chi),\in,<^*_\chi)$ to which ${\frak A}$ belongs, such that
$\|N_\alpha\| < \lambda,N_\alpha \cap \lambda \in \lambda$ and $\alpha
\subseteq N_\alpha$ and $\langle \bar a_\alpha:\alpha < \lambda
\rangle \in N_0$ (hence $\bar a_\alpha \in N_{\alpha +1}$).  
For each $\alpha \in S$, applying \scite{7.4}(1)
to $N_\alpha,\bar a_\alpha$ we get $\bold c_\alpha \in N_\alpha$ as in
Definition \scite{7.1}(2).  
So for some $\bold c^*$ and some stationary subsets of $S' \subseteq S$ of
$\lambda$ we have $\alpha \in S' \Rightarrow \bold c_\alpha =
\bold c^*$.  Now $\Delta$-end extension indiscernibility follows.
\nl
3) We prove this by induction on $n$:
\mr
\item "{$\boxtimes^n_\lambda$}"  for all $m < \omega$
a stationary $S \subseteq \lambda,\bar a_\alpha \in {}^m{\frak A}$ for
$\alpha < \lambda$ \nl
there is a stationary $S' \subseteq S$ such that: 
\nl
if $\beta < \lambda,\alpha'_\ell \in S',\alpha''_\ell \in S'$ for
$\ell < n$ and $\beta \le \alpha'_0 < \alpha'_1 < \ldots$ and $\beta
\le \alpha''_0 < \alpha''_1 < \ldots$ then $\bar a_{\alpha'_0} \char
94 \ldots \bar a_{\alpha'_{n-1}},\bar a_{\alpha''_0} \char 94 \ldots
\char 94 \bar a_{\alpha''_{n-1}}$ realizes the same type over $A \cup
\{\bar a_\gamma:\gamma < \beta\}$.
\ermn
For $n=0$ the demand is empty so $S'=S$ is as required.  For $n=1$
apply part (2).  For $n+1 > 1$ by the induction hypothesis we can find
stationary $S_1 \subseteq S$ as required in $\boxtimes^n_\lambda$.  For
each $\alpha < \lambda$ we can choose $\beta_{\alpha,\ell} =
\beta(\alpha,\ell)$ for $\ell \le n$ such that $\alpha =
\beta_{\alpha,0} < \beta_{\alpha,1} < \ldots < \beta_{\alpha,n}$ and
$0 < \ell \le n \Rightarrow \beta_{\alpha,\ell} \in S_1$.  Let $\bar
a^*_\alpha = \bar a_{\beta_{\alpha,0}} \char 94 \ldots \char 94 \bar
a_{\beta_{\alpha,n}}$ so $\bar a^*_\alpha \in {}^{m(n+1)}{\frak A}$
and apply the induction hypothesis to $m \times (n+1),S_1,\langle \bar
a^*_\alpha:\alpha < \lambda \rangle$ getting a stationary $S_2
\subseteq S_1$ as required in $\boxtimes^n_\lambda$. \nl
We claim that $S_2$ is as required.  So assume $\beta \le \alpha'_0 <
\ldots < \alpha'_n < \lambda,\beta$ and $\beta \le \alpha''_0 < \ldots
< \alpha''_n < \lambda$ and $\alpha'_\ell,\alpha''_\ell \in S_2$.  Now
\mr
\widestnumber\item{$(iii)$}
\item "{$(i)$}"   $a_{\alpha'_0} \char 94 \bar a_{\alpha'_1} \char 94
\ldots \char 94 \bar a_{\alpha'_n}$ and $\bar a_{\alpha'_0} \char 94
\bar a_{\beta(\alpha'_0,1)} \char 94 \ldots \char 94 \bar
a_{\beta(\alpha'_0,n)}$ realizes the same $\Delta$-type over $A \cup
\{\bar a_\gamma:\gamma < \beta\}$ in ${\frak A}$ \nl
[why?  as $\beta(\alpha'_1,\ell) \in S_1,\alpha'_\ell \in S_2
\subseteq S_1$ and the choice of $S_1$)]
\sn
\item "{$(ii)$}"   
$\bar a_{\alpha'_0} \char 94 \bar
a_{\beta(\alpha'_0,1)} \char 94 \ldots \char 94 \bar
a_{\beta(\alpha'_0,n)}$ is equal to $\bar a^*_{\alpha'_0}$ \nl
[why?  by the choice of $\bar a^*_{\alpha'_0}$]
\sn
\item "{$(iii)$}"   $\bar a^*_{\alpha'_0},\bar a^*_{\alpha''_0}$
realizes the same $\Delta$-type over $A \cup \{\bar a^*_\gamma:\gamma
< \beta\}$ hence over $A \cup \{\bar a_\gamma:\gamma < \beta\}$
\nl
[why?  by the choice of $S_2$]. \nl
Similarly
\sn
\item "{$(iv)$}"   $\bar a^*_{\alpha''_0}$ is equal to
$\bar a_{\alpha''_0} \char 94 \bar
a_{\beta(\alpha''_0,1)} \char 94 \ldots \char 94 \bar
a_{\beta(\alpha''_0,n)}$ 
\sn
\item "{$(v)$}"   $\bar a_{\alpha''_0}
\char 94 \bar a_{\beta(\alpha''_0,1)} \char 94 \ldots \char 94 \bar
a_{\beta(\alpha''_0,n)}$
and
$\bar a_{\alpha''_0} \char 94 \bar a_{\alpha''_1}
\char 94 \ldots \char 94 \bar a_{\alpha''_n}$ 
realizes the same $\Delta$-type over $A \cup
\{\bar a_\gamma:\gamma < \beta\}$.
\ermn
By (i)-(v) the set $S_2$ is as required in
$\boxtimes^{n+1}_\lambda$. \nl
4) The proofs of parts (2), (3) actually give this.
\hfill$\square_{\scite{7.4}}$
\enddemo
\bigskip

\proclaim{\stag{7.7} The order/unstability lemma}  Assume that
\mr
\item "{$\boxtimes_1(a)$}"  $({\frak A},\Delta)$ is a candidate
\sn
\item "{$(b)$}"  $\varphi_0(\bar x,\bar y),\varphi_1(\bar x,\bar y)
\in \Delta$ are contradictory in ${\frak A}$
\sn
\item "{$(c)$}"  $J$ is a linear order of cardinality $\lambda$
\sn
\item "{$(d)_\lambda$}"   we have $\bar a_t \in {}^m {\frak A}$
for $t \in J$ satisfies ${\frak A} \models 
\varphi_0[\bar a_s,\bar a_t] \and \varphi_1[\bar a_t,\bar a_s]$
whenever $s <_J t$
\sn
\item "{$\boxtimes_2$}"  $\lambda \ge \aleph_{\omega_1}$ \ub{or} 
$J$ is uncountable with density $\mu < |J|$.
\ermn
\ub{Then} $({\frak A},\Delta)$ is $\aleph_0$-unstable; 
even more specifically the demand in Definition \scite{7.1}(3)
holds with $\varphi_{\nu,0} = \varphi_0,\varphi_{\nu,1}
= \varphi_1$.
\endproclaim
\bn
\margintag{7.7Q}\ub{\stag{7.7Q} Question}:  What can $\{\lambda:{\frak A}$ has a
$(\Delta,\lambda)$-order$\}$ be?

We first prove a claim from which we can derive the lemma.
\proclaim{\stag{7.7A} Claim}  Assume
\mr
\item "{$(a)$}"  $({\frak A},\Delta)$ is a
$\Sigma^1_{\ell(*)}$-candidate, $\ell(*) \in \{1,2\}$ and $m <
\omega,\Phi = \Phi^m_{({\frak A},\Delta)}$ or just $\Phi \subseteq
\Phi^m_{({\frak A},\Delta)}$
\sn
\item "{$(b)$}"  $\bar{\Cal P} = \langle {\Cal P}_\alpha:\alpha <
\omega_{\ell(*)}\rangle$
\sn
\item "{$(c)$}"  ${\Cal P}_\alpha$ is a non-empty family of subsets of ${}^m
{\frak A}$
\sn
\item "{$(d)$}"  if $\alpha < \beta < \omega_{\ell(*)}$ and $B \in {\Cal
P}_\beta$ \ub{then} for some $B_0,B_1 \in {\Cal P}_\alpha$ and pair
$(\varphi_0(\bar x,\bar b),\varphi_1(\bar x,\bar b)) \in \Phi$ we have
$\ell < 2 \and \bar a \in B_\ell \Rightarrow {\frak A} \models
\varphi_\ell(\bar a,\bar b)$
\sn
\item "{$(e)$}"  if $B \in {\Cal P}_\beta$ and $\alpha < \beta <
\omega_1$ and $F$ is a function with domain $B$ and countable range,
\ub{then} there is $B' \in {\Cal P}_\alpha$ such that $B' \subseteq B$
and $F \restriction B'$ is constant
\sn
\item "{$(f)$}"  if $\ell(*) = 2$ we then in clause (e), on {\rm Rang}$(F)$
we demand just $|\text{\rm Rang}(F)| \le \aleph_1$.
\ermn
\ub{Then} $({\frak A},\Delta)$ is $\aleph_0$-unstable.
\endproclaim
\bigskip

\demo{Proof of \scite{7.7} from \scite{7.7A}}

Let $\Phi = \{(\varphi_0(\bar x,\bar y),\varphi_1(\bar x,\bar y)\}$ 
and for $\alpha < \omega_{\ell(*)}$ let

$$
{\Cal P}_\alpha = \{\bold I:\bold I \subseteq {}^m {\frak A} \text{ is
linearly ordered by } \bar \varphi \text{ and has cardinality }
\ge \aleph_\alpha\}.
$$
\mn
This should be clear.  \hfill$\square_{\scite{7.7}}$
\enddemo
\bigskip

\demo{Proof of \scite{7.7A}}  For each $\varphi(\bar x) \in \Delta$
as $\{\bar a \in {}^{\ell g(\bar x)} {\frak A}:{\frak A} \models
\varphi[\bar a]\}$ is a $\Sigma^1_{\ell(*)}$-set and let
$\{\bar a \in {}^{\ell g(\bar x)}{\frak A}:{\frak A} \models
\varphi([\bar a])\} = \{\bar a$:for some 
$\alpha < \omega_{\ell(*)}$ and $\nu \in
{}^\omega \omega,(\bar a,\nu) \in \bold C_{\varphi,\alpha}\}$ where
for each $\alpha < \omega_{\ell(*)-1}$ we have 
$\bold C_{\varphi,\alpha}$ closed subset of 
${}^{(\ell g(\bar x)+1)}({}^\omega \omega)$.  We can find $F_0,F_1$
such that if $\varphi(\bar x) \in \Delta$ and ${\frak A} \models
\varphi(\bar a)$ then $F^0_\varphi(\bar a) < \omega_{\ell(*)-1}$ and
$F^1_\beta(\bar a) \in {}^\omega \omega$ witnessing this.
For notational simplicity and \wilog \, $m=1$. Let $W =
\{w:w \subseteq {}^{\omega >} 2$ is a front hence finite$\}$.

For $w \in W$ and $n < \omega$ let 
$Q_{n,w}$ be the family of objects ${\frak x} = (n,\bar a,\bar u,\bar
\nu,\bar \varphi) = (n^{\frak x},\Gamma^{\frak x},\ldots)$ such that:
\mr
\item "{$(*)_{n,w}$}"   for unboundedly many $\alpha <
\omega_{\ell(*)}$  we can find witness (or $\alpha$-witness) ${\frak y} = 
(\langle \bar a_\ell:\ell < n \rangle,\langle B_\rho:\rho \in w
\rangle)$ which means:
{\roster
\itemitem{ $(a)$ }  $\bar u = \langle (u^0_\rho,u^1_\rho):
\rho \in w \rangle$ and $\rho \in w \Rightarrow 
u^0_\rho,u^1_\rho \subseteq n$ and $\bar \varphi = \langle \bar
\varphi^\ell:\ell < n \rangle,\bar \varphi^\ell =
(\varphi^\ell_0(x,\bar y_\ell),\varphi^\ell_1(x,\bar y_\ell)) \in \Phi$
\sn
\itemitem{ $(b)$ }  $\bar a_\ell \in {}^{\ell g(\bar y_\ell)}{\frak A}$
\sn
\itemitem{ $(c)$ }  $B_\rho \in {\Cal P}_\alpha$
\sn
\itemitem{ $(d)$ }   if $\rho \in w,b \in B_\rho$ and $\ell < n$ then
$(\varphi^\ell_0(\bar x,\bar y),\varphi^\ell_1(\bar x,\bar y)) \in
\Phi,\ell g(\bar x) = m,\ell g(\bar y)$ arbitrary (but finite) and

$$
\ell \in u^0_\rho \Rightarrow {\frak A} \models \varphi^\ell_0[b,\bar a_\ell]
$$

$$
\ell \in u^1_\ell \Rightarrow {\frak A} \models \varphi^\ell_1[b,\bar a_\ell]
$$
\sn
\itemitem{ $(e)$ }  if $\nu \ne \rho$ are from $w$ then 
$(u^0_\rho \cap u^1_\nu \ne \emptyset) \vee (u^0_\nu \cap
u^1_\rho \ne \emptyset)$
\sn
\itemitem{ $(f)$ }   $\bar \nu = \langle \nu^\ell_{\rho,\ell}:\rho \in
w,i \in \{0,1\}$ and $\ell \in u^i_\rho \rangle$
\sn
\itemitem{ $(g)$ }  if $b \in B_\rho,i \in \{0,1\},\ell \in u^i_\rho$
then $F^0_{\varphi^\ell_i}(b,\bar a_\ell) =
\alpha^i_{\rho,\ell},(F^1_{\varphi^\ell_i}(b,\bar a_\ell))
\restriction n = \nu^i_{\rho,\ell}$.
\endroster}
\ermn
Clearly
\mr
\item "{$(*)_1$}"   $Q_{0,\{<>\}} \ne \emptyset$. \nl
[Why?  Let ${\frak x} = (0,<>,<>,<>,<>)$ and if $\alpha <
\omega_{\ell(*)}$ choose $\bold I \in {\Cal P}_\alpha$ we let 
$B_{<>} = \bold I$
\sn
\item "{$(*)_2$}"   if ${\frak x} \in Q_{n,w}$ and for $\rho \in
w,F_\rho$ is an $(n+1)$-place function with domain ${\frak A}$ and
range $\subseteq \omega_{\ell(*)-1}$ or just countable range,
\ub{then} there is $\langle {\frak y}^1_\alpha:\alpha <
\omega_{\ell(*)} \rangle$ such that ${\frak y}^1_\alpha$ is an
$\alpha$-witness for ${\frak x} \in Q_{n,w}$
and $\langle F_\rho(\bar a^{\frak y},\bar b_\rho):
\rho \in w \rangle$ is the same for all $b_\rho \in B^{{\frak
y}_\alpha}_\rho,\alpha < \omega_{\ell(*)}$ where $\ell g(\bar a) = n$ \nl
[why?  as ${\frak x}_1 \in Q_{n,\omega}$ we know that for some unbounded
$Y \subseteq \omega_{\ell(*)}$ for 
each $\alpha_1 \in Y$ there is an $\alpha$-witness $\langle
a^\alpha_\ell:\ell < n \rangle \char 94 \langle B^\alpha_\rho:\rho \in
w \rangle$ as required in $(*)_{n,\omega}$.  
Let $\alpha < \omega_1$ and $\beta(\alpha) = \text{
Min}(Y \backslash (\alpha +1))$.  Now for each $\rho \in w$ as
$B^{\beta(\alpha)}_\rho \in {\Cal P}_{\beta(\alpha)}$ and the set
$\{F_\rho(a^{\beta(\alpha)}_0,\dotsc,a^{\beta(\alpha)}_{n-1},\bar
b):\bar b \in B^\beta_\rho\}$ is countable, by clause (d) of the
assumption we can find 
$c^\alpha_\rho$ such that $B^{\alpha,*}_\rho \subseteq \{b \in
B^{\beta(\alpha)}_\rho:
F_\rho(a^{\beta(\alpha)}_0,\dotsc,a^{\beta(\alpha)}_{n-1},b) = 
c^\alpha_\rho\}$ belong to ${\Cal P}_\alpha$.  As the set $\{\langle
c^\alpha_\rho:\rho \in w \rangle:\alpha < \omega_{\ell(*)}\}$ is countable,
there is a sequence $\langle c^*_\rho:\rho \in w \rangle$ such that
the set $Y' = \{\alpha < \omega_1:c^\alpha_\rho = c^*_\rho$ for $\rho
\in \omega\}$ is uncountable for $\alpha <
\omega_{\ell(*)},\gamma_\alpha = \text{ Min}(Y' \backslash \alpha)$
and ${\frak y}_\alpha = (\langle a^{\beta(\gamma_\alpha)}_\ell:\ell < n
\rangle),\langle B^{\beta(\gamma_\alpha),*}_\rho:\rho \in w \rangle)$.
Clearly $\langle {\frak y}_\alpha:\alpha < \omega_{\ell(*)} \rangle$ is as
required.]
\sn
\item "{$(*)_3$}"   if ${\frak x}_1 \in Q_{n,w}$ and $\rho \in w$ and
$u = (w \backslash \{\rho\}) \cup \{\rho \char 94 <0>,\rho \char 94
<1>\}$, \ub{then} there is ${\frak x}_2 \in Q_{n+1,u}$ such that:
{\roster
\itemitem{ $(*)$ }  if ${\frak y} = (\langle a_\ell:\ell < n+1
\rangle,\langle B_\rho:\rho \in u \rangle)$ is an $\alpha$-witness of
${\frak x}_2$ and $\alpha < \omega_{\ell(*)}$ \ub{then} 
${\frak y}'_\alpha = (\langle a_\ell:\ell < n \rangle,
\langle B'_\rho:\rho \in w \rangle)$ is an $\alpha$-witness
for ${\frak x}_1$ where $B'_\rho$ is $B^{{\frak y}_2}_\rho$ if $\rho
\in u \cap w$ and $B'_\rho = B^{{\frak y}_2}_{\rho \char 94 <0>} \cup
B^{{\frak y}_2}_{\rho \char 94 <1>}$ if $\rho = \rho$.
\nl
[why?  similar to the proof of $(*)_2$ using clause (e) of the
assumption this time.]
\endroster}
\ermn
Together it is not hard to prove the non $\aleph_0$-unstability (as in
\cite{Sh:522}). .  \hfill$\square_{\scite{7.7}}$ 
\enddemo
\bigskip

\remark{\stag{7.8} Remark}  1) This 
claim can be generalized replacing $\aleph_0$
by $\mu$, strong limit singular of cofinality $\aleph_0$. 
\endremark
\bn
\centerline {$* \qquad * \qquad *$}
\bigskip

\definition{\stag{7.9} Definition}  1) tp$_\Delta(\bar a,A,{\frak A})
= \{\varphi(\bar x,\bar b):\varphi(\bar x,\bar y) \in \Delta$ and $b
\in  {}^{\ell g(\bar y)}(A)$ and 
${\frak A} \models \varphi[\bar a,\bar b]\}$. \nl
2) $\Phi_{{\frak A},\Delta,A}^{pr,m} = 
\{(\varphi_0(\bar x,\bar b),\varphi_1(\bar x,\bar b)):\varphi_0(\bar
x,\bar y),\varphi_1(\bar x,\bar y)$ belongs to $\Delta$ and $\bar b
\in {}^{\ell g(\bar y)} A$ and $\bar x = \langle x_\ell:\ell < m
\rangle $ and ${\frak A} \models \neg(\exists \bar x)[\varphi_0(\bar
x,\bar b) \and \varphi_1(\bar x,\bar b)]\}$ \nl
where $A \subseteq {\frak A},\Delta$ a set of ${\frak A}$-formulas,
and so
\nl
$\Phi^{pr,m}_{{\frak A},\Delta} 
= \{(\varphi_0(\bar x,\bar y),\varphi_1(\bar x,\bar y)):
\varphi_0,\varphi_1 \in \Delta,{\frak A} \models \neg \exists \bar
y \exists \bar x[\varphi_0(\bar x,\bar y) \and \varphi_1(\bar x,\bar y)]\}$.
\nl
3) ${\bold S}^m_\Delta(A,{\frak A}) = \{\text{tp}_\Delta(\bar a,A,{\frak
A}):\bar a \in {}^m {\frak A}\}$ where $A \subseteq {\frak A}$ and
$\Delta$ a set of $\Bbb L(\tau_{\frak A})$-formulas
\enddefinition
\bigskip

\definition{\stag{7.10} Definition}  1) We say $({\frak A},\Delta)$ is
$(\mu,\Delta,\lambda)$-unstable \ub{if} there are $M \subseteq
{\frak A},m < \omega$ and $\langle \bar a_\alpha:\alpha < \lambda \rangle$
such that:
\mr
\item "{$(a)$}"  $\bar a_\alpha \in {}^m{\frak A}$
\sn
\item "{$(b)$}"  if $\alpha \ne \beta$ are $< \lambda$ then for some
$(\varphi_0(\bar x,\bar b),\varphi_1(\bar x,\bar b)) \in
\Phi^{m,pr}_{{\frak A},\Delta,M}$ (see Definition \scite{7.9} below) 
we have $\varphi_0(\bar x,\bar b) \in 
\text{ tp}_\Delta(\bar a_\alpha,M,{\frak A})$ 
and $\varphi_1(\bar x,\bar b) \in \text{ tp}_\Delta(\bar a_\beta,M,{\frak A})$
\sn
\item "{$(c)$}"   $\|M\| \le \mu$.
\ermn
1A) Let ${\frak A}$ be $(\aleph_0,\Delta,\text{per})$-unstable mean
that $({\frak A},\Delta)$ is $\aleph_0$-unstable; here per stands for perfect.
\nl
2) We add ``weakly" if we weaken clause (b) to
\mr
\item "{$(b)^-$}"  tp$_\Delta(\bar a_\eta,M,{\frak A}) \ne \text{
tp}_\Delta(\bar a_\nu,M,{\frak A})$ for $\eta \ne \nu$ from $X$ \nl
(so if $\Delta$ is closed under negation there is no difference); in
part (1), $X = \lambda$ and in part (2), $X = {}^\omega 2$.
\ermn
3) We use $(\mu_0,\Delta,x,\Bbb Q)$ where $\Bbb Q$ is a forcing notion
\ub{if} the example is found in
$\bold V^{\Bbb Q}$ such that usually $M$ is in $\bold V$ and we add 
an additional possibility if $x = \text{ per}^{\bold V}$ then $M \in \bold V$
and $X = ({}^\omega 2)^{\bold V}$ (here per stands for perfect).
\nl
4) We may replace ``a forcing notion $\Bbb Q$" by a family ${\frak K}$
of forcing notions (e.g. the family of c.c.c. ones) meaning: for at
least one of them. \nl
5) We replace stable by unstable for the negation.
\enddefinition
\bn
\centerline {$* \qquad * \qquad *$}
\bn
\margintag{7.11}\ub{\stag{7.11} Observation}:  1) If $\Delta$ is closed under negation,
\ub{then} ${\frak A}$ is weakly $(\aleph_0,\Delta,\lambda)$-unstable
iff ${\frak A}$ is $(\aleph_0,\Delta,\lambda)$-unstable.
\bigskip

\definition{\stag{7.12} Definition}  Let $({\frak A},\Delta)$ be a
$\Sigma^1_{\ell(*)}$-candidate where $\ell(*) \in \{1,2\}$.  For $m <
\omega$ and $B \subseteq {}^m{\frak A}$ we define rk$^{\ell(*)}(B) =
\text{ rk}^{\ell(*)}(B,\Delta,{\frak A})$, an ordinal or infty or $-1$
by defining for any ordinal $\alpha$ when rk$^{\ell(*)}(B) \ge \alpha$
by induction on $\alpha$.
\enddefinition
\bn
\ub{Case 1}:  $\alpha = 0$.

rk$^{\ell(*)}(B) \ge \alpha$ iff $B \ne \emptyset$.
\bn
\ub{Case 2}:  $\alpha$ limit.

rk$^{\ell(*)}(B) \ge \alpha$ iff rk$^{\ell(*)}(B) \ge \beta$ for every
$\beta < \alpha$.
\bn
\ub{Case 3}:  $\alpha = \beta +1$.

rk$^{\ell(*)}(B) \ge \alpha$ iff (a) + (b) holds where
\mr
\item "{$(a)$}"  if $B = \cup\{B_i:i < \aleph_{\ell(*)-1}\}$ then for
some $i$ we have rk$^{\ell(*)}(B_i) \ge \beta$
\sn
\item "{$(b)$}"  we can find $\bar \varphi(\bar x,\bar b) \in
\Phi^m_{{\frak A},\Delta}$ and 
$B_0,B_1 \subseteq B$ such that rk$^{\ell(*)}(B_i) \ge \beta$ and
$\bar a \in B_\ell \Rightarrow {\frak A} \models \varphi_\ell(\bar
a,\bar b)$ for $\ell = 0,1$.
\endroster
\bn
\margintag{7.13}\ub{\stag{7.13} Observation}:  Assume $({\frak A},\Delta)$ is
$\aleph_{\ell(*)}$-candidate, $\ell(*) \in \{1,2\}$. \nl
1) If $\alpha \le \beta$ are ordinals and rk$^{\ell(*)}(B) \ge \beta$
then rk$^{\ell(*)}(B) \ge \alpha$. \nl
2) rk$^{\ell(*)}(B) \in \text{ Ord} \cup \{-1,\infty\}$ is well
defined (for $B \subseteq {}^m {\frak A})$. \nl
3) If $B_1 \subseteq B_2 \subseteq {\frak A}$ then rk$^{\ell(*)}(B_1)
\le \text{ rk}^{\ell(*)}(B_2)$.
\bigskip

\demo{Proof}  Trivial.
\enddemo
\bigskip

\proclaim{\stag{7.14} Claim}  The following are equivalent if
$2^{\aleph_0} \ge \aleph_{\ell(*)},({\frak A},\Delta)$ is a
$\Sigma^1_{\ell(*)}$-candidate:
\mr
\item "{$(a)$}"  {\rm rk}$^{\ell(*)}({}^m{\frak A}) \ge \omega_{\ell(*)}$
\sn
\item "{$(b)$}"  ${\frak A}$ is $(\aleph_0,\Delta)$-unstable
\sn
\item "{$(c)$}"  ${\frak A}$ is
$(\aleph_0,\Delta,\aleph_{\ell(*)})$-unstable
\sn
\item "{$(d)$}"  {\rm rk}$^{\ell(*)}({\frak A}) = \infty$.
\endroster
\endproclaim
\bigskip

\demo{Proof}  \ub{$(a) \Rightarrow (b)$}.

Let ${\Cal P}_\alpha = \{B \subseteq {}^m{\frak
A}:\text{rk}^{\ell(*)}(B) \ge \alpha\}$ and apply \scite{7.7A}.
\mn
\ub{$(b) \Rightarrow (c)$}:  Trivial.
\mn
\ub{$(c) \Rightarrow (d)$}:

Let $A \subseteq {\frak A}$ be countable and $\{\bar a_\alpha:\alpha <
\aleph_{\ell(*)}\} \subseteq {}^m {\frak A}$ exemplifies that ${\frak
A}$ is $(\aleph_0,\Delta,\aleph_{\ell(*)})$-unstable. \nl
Without loss of generality
\mr
\item "{$(*)$}"  if $\bar b \subseteq A,\varphi(\bar x,\bar y) \in
\Delta$ and $\{\alpha < \aleph_{\ell(*)}:{\frak A} \models
\varphi(\bar a_\alpha,\bar b)\}$ is bounded then it is empty.
\ermn
Now let ${\Cal P} = \{\{\bar a_\alpha:\alpha \in S\}:S
\subseteq \aleph_{\ell(*)}$ is unbounded.  Now we can prove by
induction on $\alpha$ that $B \in {\Cal P}
\Rightarrow \text{ rk}^{\ell(*)}(B) \ge \alpha$.
\hfil$\square_{\scite{7.14}}$
\mn
\ub{$(d) \Rightarrow (a)$}:  Trivial.
\enddemo
\bigskip

\definition{\stag{7.15} Definition}  If $p$ is a $(\Delta_1,m)$-type in
over $A$ in ${\frak A}$ (i.e. a set of formulas $\varphi(\bar x,\bar
a)$ with $\varphi(\bar x,\bar y) \in \Delta_1,\bar a \subseteq A$), we
let
$$
\align
\text{rk}^{\ell(*)}(p,\Delta,{\frak A}) = \text{
Min}\{\text{rk}^{\ell(*)} \dbca_{\ell < n} \varphi_\ell({}^m{\frak
A},\bar b_\ell),\Delta,{\frak A}):&n < \omega \\
  &\text{ and } \varphi_\ell(\bar x,\bar b_\ell) \in p \text{ for }
\ell < n\}.
\endalign
$$
\enddefinition
\bn
\margintag{7.16}\ub{\stag{7.16} Observation}  1) If $p \subseteq q$ (or just $q \vdash
p$) are $(\Delta,m)$-types in ${\frak A}$ \ub{then}
rk$^{\ell(*)}(q,\Delta,{\frak A}) \le \text{
rk}^{\ell(*)}(p,\Delta,A)$. \nl
2) If $q$ is a $(\Delta,m)$-type in ${\frak A}$ \ub{then} for some finite
$p \subseteq q$ we have

$$
\text{rk}^{\ell(*)}(q,\Delta,{\frak A}) = \text{
rk}^{\ell(*)}(p,\Delta,{\frak A})
$$
\mn
hence

$$
p \subseteq r \subseteq q \Rightarrow 
\text{ rk}^{\ell(*)}(r,\Delta,{\frak A}) =
\text{ rk}^{\ell(*)}(p,\Delta,{\frak A}).
$$
\bigskip

\proclaim{\stag{7.17} Claim}  1) In \scite{7.14} we can add
\mr
\item "{$(e)$}"  $({\frak A},\Delta)$ is not $\aleph_0$-stable
\sn
\item "{$(f)$}"  for some $\mu < \lambda$ the pair $({\frak A},\Delta)$ is
$(\mu,\Delta,\lambda)$-unstable and $\aleph_{\ell(*)} < \lambda$.
\endroster
\endproclaim
\bigskip

\demo{Proof}  \ub{$\neg(e) \Rightarrow \neg(c)$}.

Let $M \prec ({\Cal H}(\chi),\in,<^*_\chi)$ be countable such that $x
\in M$ for suitable $x$ and $m < \omega$.  For every $\bar a \in
{}^m{\frak A}$ there is a function $\bold c_{\bar a} \in M$ from
$\Phi^m_{({\frak A},\Delta)}$ to $\{0,1\}$ as in Definition
\scite{7.1}.  So if $\bar a_i \in {}^m{\frak A}$ for $i < \omega_{\ell(*)}$
then for some $i < j < \omega_{\ell(*)}$ we have 
$\bold c_{\bar a_i} = \bold c_{\bar a_j}$ because
$M$ is countable.  So clearly $(c)$ fails $\bar \varphi$.
\mn
\ub{$(e) \Rightarrow (c)$}.

Fix $({\Cal H}(\chi_0),\in,<^*_\chi)$ and let

$$
\align
{\Cal S}_0 = \{M \prec ({\Cal H}(\chi_0),\in,<^*_\chi):&{\frak A} \in M
\text{ and} \\
  &\|M\| = \aleph_{\ell(*)-1} \text{ and } \omega_{\ell(*)-1} +1
\subseteq M\}.
\endalign
$$
\mn
For $m < \omega$ and $\bold I \subseteq {}^m{\frak A}$ let ${\Cal
J}_{\bold I} = {\Cal J}[\bold I]$ be the family of ${\Cal S} \subseteq
{\Cal S}_0$ such that: we can find $\langle F_x,\bold c_x:x \in {\Cal
H}(\chi) \rangle$ (a witness) such that:
\mr
\item "{$(\alpha)$}"  $c_x:\Phi^m_{{\frak A},\Delta} \rightarrow
\{0,1\}$
\sn
\item "{$(\beta)$}"  $F_x:{}^{\omega >}({\Cal H}(\chi)) \rightarrow
{\Cal H}(\chi)$
\sn
\item "{$(\gamma)$}"  if $M \in S$ is closed under $F_x$ for $x \in M$
then for every $\bar a \in \bold I$ for some 
$y \in M,c_y$ is a witness for tp$(\bar a_M,M \cap {\frak A},{\frak A})$.
\ermn
Clearly ${\Cal J}_{\bold I}$ is a normal ideal on ${\Cal S}_0$.  Also
if $m < \omega \Rightarrow S_0 \in {\Cal J}[{}^m{\frak A}]$ then
increasing $\chi$ we get the desired result.  Toward contradiction
assume that $m < \omega$ and ${\Cal S} \notin {\Cal J}[{}^m{\frak A}]$
and let ${\Cal P}$ (i.e. ${\Cal P}_\alpha = {\Cal P}$ for $\alpha <
\omega_{\ell(*)}$) be the family of $\bold I \subseteq {}^m{\frak A}$
such that ${\Cal S}_0 \notin {\Cal J}_{\bold I}$.

We now finish by \scite{7.7A} once we prove
\mr
\item "{$\circledast$}"  if $\bold I \in {\Cal P}$ then for some $\bar
\varphi(\bar x,\bar b) \in \Phi^m_{{\frak A},\Delta}$ for each $\ell <
2$ the set $\bold I^\ell_{\bar \varphi(x,\bar b)}$ is $\{\bar a \in
\bold I:{\frak A} \models \varphi_\ell(\bar a,\bar b)\}$ belong to
${\Cal P}$.
\ermn
If not, for every $\bar \varphi(\bar x,\bar b) \in
\Phi^m_{{\frak A},\Delta}$ there is $\ell = \bold c[\bar \varphi(\bar
x,\bar b)] < 2$ and $\langle (F^{\bar \varphi(\bar x,\bar b)}_x,\bold
c^{\bar \varphi(\bar x,\bar b)}_x):x \in {\Cal H}(\chi) \rangle$
witnessing ${\Cal S}_0 \in \bold J[\bold I^\ell_{\bar \varphi}]$.

Define $(F_y,\bold c_y)$ for $y \in {\Cal H}(\chi)$ by: if $y =
\langle x,\bar \varphi(\bar x,\bar b) \rangle$ then $F_y = F^{\bar
\varphi(\bar x,\bar b)}_x,c_y = \bold c^{\bar \varphi(\bar x,b)}_x$,
otherwise $\bold c$.

Clearly we can find $M \in {\Cal S}_0$ such that
\mr
\item "{$\circledast_1$}"   if $\varphi(\bar x,\bar b) \in
\Phi^m_{{\frak A},\Delta} \cap N$ and $x \in M$ then $M$ is closed
under $F^{\bar \varphi(\bar x,\bar b)}_x$
\sn
\item "{$\circledast_2$}"   for some $\bar a \in {}^m{\frak A}$, no
$\bold c_y,y \in M$ defines tp$_\Delta(\bar a,M \cap {\frak A},{\frak
A})$.
\ermn
But $\bold c$ does it!  So we are done.
\mn
\ub{$(f) \Rightarrow (d)$}.

Like $(c) \Rightarrow (d)$.
\mn
\ub{$(c) \Rightarrow (f)$}.

Just use $\mu = \aleph_0$.  \hfill$\square_{\scite{7.17}}$
\enddemo
\bigskip

\proclaim{\stag{7.18} Claim}  Assume that $({\frak A},\Delta)$ is a
$\Sigma^1_{\ell(*)}$-candidate, $\ell(*) \in \{1,2\}$ and is
$\mu$-stable.

For some $\xi < \omega_1$ we have: if $\lambda \ge \mu,m < \omega,A
\subseteq {\frak A},|A| \le \lambda$ and $\bar a_\alpha \in {}^m{\frak
A}$ for $\alpha < \lambda^{+ \xi}$ \ub{then} for some $S \subseteq
\lambda^{+ \xi}$ of cardinality $\lambda$ the sequence $\langle \bar
a_\alpha:\alpha \in S \rangle$ is $\Delta$-indiscernible over $A$ in
${\frak A}$.
\endproclaim
\bigskip

\remark{Remark}  See more in \cite{Sh:F562}.
\endremark
\bigskip

\demo{Proof}  Assume not.  For $\xi < \omega_1$ let 

$$
\align
{\Cal P}_\xi = \{\{\bar a_\alpha:\alpha < \lambda^{+ \xi}\}:&\text{
for some } \lambda \ge \mu, \text{ for no } 
S \subseteq \lambda^{+ \xi} \text{ of cardinality} \\
  &\lambda \text{ is } \langle \bar a_\alpha:\alpha \in S \rangle
\text{ is } \Delta \text{-indiscernible over } A \text{ in } {\frak
A}\}.
\endalign
$$
\mn
The point is:
\mr
\item "{$\circledast$}"  if $\lambda^{+ \xi}$ is regular, $\xi > 0,A
\subseteq {\frak A},|A| \le \lambda,\bar a_\alpha \in {}^m {\frak A}$
for $\alpha < \lambda^{+ \xi}$ and $S \subseteq \lambda^{+ \xi}$ is
stationary then (a) or (b) where
{\roster
\itemitem{ $(a)$ }  for some club $E$ of $\lambda,\langle \bar
a_\alpha:\alpha \in S \cap E \rangle$ is $\Delta$-indiscernible over
$A$ in ${\frak A}$
\sn
\itemitem{ $(b)$ }  for some $m < \omega$ and club $E^*_n$ of
$\lambda^{+ \xi}$ we have
\sn
\itemitem{ $(b)_m(i)$ }  $\langle \bar a_\alpha:\alpha \in S \cap
E^*_m \rangle$ is $(\Delta,m)$-end extension indiscernible
\sn
\itemitem{ ${{}}(ii)$ }  for no club $E' \subseteq E^*_m$ of
$\lambda^{+ \xi}$ is $\langle \bar a_\alpha:\alpha \in S \cap E
\rangle$ a sequence which is $(\Delta,m+1)$-end extension indiscernible.
\endroster}
\ermn
Clearly clause (a) is impossible by our present assumptions so let
$E^*,m$ be as in clause (b).  By claim \scite{7.4}(4) there is a club
$E$ of $\lambda$ and $\langle f_n:n < \omega \rangle$ as there and let
$S^*_\gamma = \{\alpha \in S:f_{m+1}(\alpha) = \gamma\}$, so $\alpha >
\gamma,{\Cal P}_{m+1} = \{\gamma:S^*_\gamma$ is stationary.  Without loss of
generality $E^* \subseteq E$ and $\gamma \notin S_{m+1} \Rightarrow
S^*_\gamma = \emptyset$.  Without loss of generality $f_{m+1}$ is as
in claim \scite{7.19} below.

So by $(b)_m(ii)$ clearly $\Gamma_{m+1}$ is not a singleton (and it
cannot be empty), so we clearly have finished.
\hfill$\square_{\scite{7.18}}$
\enddemo
\bigskip

\proclaim{\stag{7.19} Claim}  Let $A,\langle \bar a_\alpha:\alpha <
\lambda \rangle,E,\langle f_n:n < \omega \rangle$ be as in
\scite{7.4}(4).  \ub{Then} without loss of generality 
(possibly shrinking $E$ and changing the $f_n$'s) we can add
\mr
\item "{$(iii)$}"  if $m < \omega$ and $\gamma_1 \ne \gamma_2$ are in
Rang$(f_{n+1})$ but $f_{n+1}(\alpha_1) = \gamma_1 \wedge
f_{n+1}(\alpha_2) = \gamma_2 \Rightarrow f_n(\gamma_1) =
f_n(\gamma_2)$ letting $S = \{\alpha:f_n(\alpha) = f_n(\gamma_1) =
f_n(\gamma_2)\}$ and $\beta = \text{ Min}(S \cap E \backslash
(\gamma_1 +1) \backslash (\gamma_2 + 1)$, 
\ermn
\ub{then} for some formula $\varphi(\bar x_0,\dotsc,\bar x_n)$ with
parameters from $A \cup \{a_\gamma:\gamma < \beta\}$ such that:
\mr
\item "{$(*)$}"  if $i<2,\alpha'_0 < \ldots < \alpha'_{n+1}$ are from
$S \cap E(\ell \le n)(\exists \alpha)(f_n(\alpha) = f_n(\alpha'_\ell)
\wedge f_{n+1}(\alpha) = \gamma_i)$ and $f(\alpha'_0) = \gamma_i$ then
${\frak A} \models \varphi[\bar a_{\alpha'_0},\dotsc,\bar
a_{\alpha'_n}] \Leftrightarrow i=0$.
\endroster
\endproclaim
\bigskip

\demo{Proof}  Easy.
\enddemo
\newpage

\head {Glossary} \endhead
 \spuriousreset
\bn
\S0 \ub{Introduction}
\mn
Theorem \scite{0.1}: No Polish group
\sn
Thesis \scite{0.2}: Polish algebras are large
\sn
Question \scite{0.3}:  What can be Aut$(\Bbb A),\Bbb A$ uncountable
\sn
Question \scite{0.4}:  Is there model theory of Polish algebras
\sn
Example \scite{0.5}:  Adding many Cohens
\sn
Example \scite{0.6}:  The complex field, the real field
\sn
Conjecture \scite{0.7}:  There is a dichotomy
\sn
Thesis \scite{0.8}:  Classification theory of such structures exists
\sn
Theorem \scite{0.9}:  There is a $F_\sigma$ abelian groups with complicated
categoricity behaviour
\sn
Conclusion \scite{0.10}: Categoricity can stop at $\aleph_n$
\sn
Theorem \scite{0.11}:  Indiscernibles exist
\sn
Definition \scite{0.12}:  Categoricity
\sn
Categoricity Question \scite{0.13}:  Is there such a classification
theory for equational theories
\sn
Notation \scite{0.21}:
\sn
Definition \scite{0.22}: group words
\bn
\S1 Metric groups and metric models
\sn
Definition \scite{1.1}: metric group, metric semigroups
\sn
Notation \scite{1.1A}:  For metric group $\bold M,\bold d_M$ is the
metric, $e_{\bold M}$ the unit, $G_{\bold M}$ the group
\sn
Definition \scite{1.2}: specially (metric group), specially$^+,\bar
\zeta$ is strongly O.K.
\sn
Observation \scite{1.3}: basic properties
\sn
Definition \scite{1.4}: automorphism of countable structures,
endomorphism semi group, monomorphism semi group
\sn
Claim \scite{1.5}: the above are separable metric groups semi groups
\sn
Definition \scite{1.6}: ${\frak a}$ is a metric algebra; unitary,
complete; specially$^{(+)}$ unitary; partial
\bn
\S2 \ub{Semi-metric groups: automorphism groups of uncountable
structures}
\mn
Definition \scite{1.7}: 1) $G$ is a complete/special metric group.
\nl

2) Similarly for semi-group
\sn
Discussion \scite{1.7A}: 1) Note that in \scite{1.7} we do not
necessarily have metric groups.
\nl

2)-4) Variants.
\sn
Definition \scite{1.8}: The sequence $\bar A$ is an 
$\omega$-representation of $\Bbb A$, and related
metrics
\sn
Claim \scite{1.9}: when $\Bbb A,\bar A$ gives a [semi]
complete/specially$^+$ metric group of 
\nl

\hskip10pt automorphisms (or semi group of
endomorphisms or semi group 
\nl

\hskip10pt of monomorphisms).
\sn
Discussion \scite{1.11A}: On variants of \scite{1.9}
\sn
Definition \scite{1.12}: $\Bbb A$-beautiful term and some distance
functions depending 
\nl

\hskip10pt on a representation
\sn
Claim \scite{1.12A}: beautiful terms induce operations on endomorphism
semi-group 
\sn
Claim \scite{1.13}: how nice is the derived metric algebra from
auto/endo/mono 
\nl

\hskip10pt semi-groups
\sn
Question \scite{1.20}: can an uncountable Polish algebra be free for some
variety?
\sn
Observation \scite{1.21}:  example answering the question
\sn
Remark \scite{1.22}:  another metric
\sn
Definition \scite{1.30}:  ${\frak a} = (M,{\frak d},\bold U)$ is a
metric topological algebra
\sn
Claim \scite{1.31}: sufficient condition for being complete metric
topological algebra
\sn
Discussion \scite{1.32}:  1) Replacing metric by a topology.
\nl

2) Replacing automorphisms by expansion to models of a universal Horn theory. 
\bn
\S3 Compactness of metric algebras
\sn
The completeness Lemma \scite{s.1}: give sufficient conditions for
solvability of a set of equations in a Polish algebra.
\sn
Remark \scite{s.2}:  Explaining \scite{s.1}.
\sn
Fact \scite{s.3}:  Recall free group is torsion free with no non-trivial
element divisible.
\sn
Fact \scite{s.4}:  Recall another consequence of freeness.
\sn
Claim \scite{s.5}: sufficient conditions for complete metric algebra
to be far from free.
\sn
Remark \scite{s.5A}:  On stable variants on the theorems (e.g. $\|M_{\frak
a}\| <$ cov(meagre) instead $M$ countable).
\sn
Remark \scite{3.5}:  On variants.
\bn
\S4 \ub{Conclusions}
\mn
Conclusion \scite{c.1}:  if $(G,{\frak d})$ is a complete metric space
of density $< |G|$ then $G$ is similar to free; semi-complete; is
enough; not discrete is enough.
\sn
Conclusion \scite{c.2}:  There is no free uncountable Polish group.
\sn
Claim \scite{c.3}:  Strengthening the ``non-free" replacing free.
\sn
Remark \scite{c.3A}:  On related ranks; this conclusion confirms the
complicatedness thesis.
\sn
Conclusion \scite{c.4}:  On Aut$({\Bbb A})$.
\sn
Claim \scite{c.6}: when the proof of \cite{Sh:744} works
\bn
\S5 \ub{Quite free but not free abelian groups}
\mn
Question \scite{6.0}:  1) Is the ``freeness of a (definable) abelian
group" absolute?
\nl

2),3) Variants.
\sn
Definition \scite{6.1}:  of $G_{k(*)}$
\sn
Claim \scite{6.2}: $G_{k(*)}$ is $\aleph_1$-free
\sn
Definition \scite{6.3}:  1),2) the subgroups $G_U,G_{U,u}$.
\nl

3) The set of equations $\xi_{U_1},\Xi_{U,u}$
\sn
Claim \scite{6.4}: How $G_U$ is generated
\sn
Main Claim \scite{6.5}: 1) $G_{U \cup u}/G_{U,u}$ is free if $|u| \le
k,|U| \le \aleph_k$.
\nl

2) $G_{U \cup u}/G_{U,u}$ is free if $|U| \le
\aleph_{k(*)-|u|}$
\sn
Claim \scite{6.6}: $G_U$ is not free if $|U| \ge \aleph_{k(*)+1}$
\sn
Claim \sciteu{6.6A}:  $G_{k(*)}$ is Borel (even $F_\sigma$) abelian
group
\sn
Conclusion \scite{6.6G}:  $G_{k(*)}$ is a Borel and even $F_\sigma$
abelian group, 
\nl

\hskip10pt $\aleph_{k(*)+1}$-free and if $2^{\aleph_0} \ge
\aleph_{k(*)+2}$, 
\nl

\hskip10pt not $\aleph_{k(*)+2}$-free.
\sn
Question \scite{6.10}:  Does appropriate $\Bbb Q$ necessarily add
reals
\sn
Fact \scite{6.11}:  Information from \cite{Sh:402}.
\sn
Corollary \scite{6.12}:  1) $G_{k(*)}$ purely embeddable into
${}^\omega \Bbb Z$.
\nl

2) Forcing making it free.
\bn
\S6 \ub{Beginning of stability theory}
\mn
Convention \scite{7.0}:  $\tau_{\frak A}$ countable, members of
${\frak A}$ are reals
\sn
Definition \scite{6a.deL}:  ${\frak A}$-formula, pairs of formulas and
set $\Delta$ of pairs
\sn
Definition \scite{7.1}:  $({\frak A},\Delta)$  a candidate, stability
\sn
Observation \scite{7.2}:  Basic facts.
\sn
Claim \scite{7.4}:  From stability to the existence of indiscernibles
\sn
Definition \scite{7.4A}:  ${\frak A}$ has $(\lambda,\Delta)$-order 
\sn
Claim \scite{7.7}: From order to unstability
\sn
Question \scite{7.7Q}:  What can be $\{\lambda:{\frak A}$ has
$(\Delta,\lambda)$-order$\}$?
\sn
Claim \scite{7.7A}: Sufficient conditions for being unstable
(i.e. having a perfect set of 
\nl

\hskip10pt pairwise explicitly contradictory type)
\sn 
Remark \scite{7.8}: Replacing $\aleph$ by, e.g. $\beth_\omega$
\sn
Comment: nonstable is unstable
\sn
Definition \scite{7.9}: tp$_\Delta(\bar a,A,{\frak A}),
\Phi^{pr,m}_{\Bbb A,\Delta,A}$ and $\bold S^m_\Delta(A,{\frak A})$
\sn
Definition \scite{7.10}: $(\Bbb A,\Delta)$ is
$(\mu,\Delta,\lambda)$-unstable
\sn
Observation \scite{7.11}: Weakly stable/unstable
\sn
Definition \scite{7.12}:  of rk$(B)$
\sn
Observation \scite{7.13}: Properties of rk
\sn  
Claim \scite{7.14}: Equivalences to rank being infinite
\sn
Definition \scite{7.15}:  rk in more cases
\sn
Subclaim \scite{7.16}:  properties of rk
\sn
Claim \scite{7.17}: More cases of equivalence in \scite{7.14}
\sn
Claim \scite{7.18}:  Existence of indiscernible
\sn
Claim \scite{7.19}: helping \scite{7.18}
\newpage


\nocite{ignore-this-bibtex-warning} 
\newpage
    
REFERENCES.  
\bibliographystyle{lit-plain}
\bibliography{lista,listb,listx,listf,liste}

\enddocument